\documentclass[12pt]{article}
\usepackage{graphicx, fullpage}
\usepackage{amsmath,amssymb,verbatim,amsthm,bbm}
\usepackage{natbib}
\usepackage[ruled,vlined]{algorithm2e}
\setcitestyle{authoryear,round}
\usepackage{xcolor}
\usepackage{caption}
\usepackage{dsfont}
\usepackage{listings}
\usepackage{todonotes}
\usepackage{subcaption}
\usepackage{booktabs}
\usepackage[normalem]{ulem}
\usepackage{listings}
\usepackage{xcolor}  
\usepackage{url}
\usepackage{appendix}

\newcommand{\Pcon}{\stackrel{\textrm{P}}{\rightarrow}}
\newcommand{\approxD}{\stackrel{\textrm{D}}{\approx}}
\newcommand{\Dcon}{\stackrel{\textrm{D}}{\rightarrow}}

\newtheorem{lemma}{Lemma}[section]
\newtheorem{definition}{Definition}[section]
\newtheorem{example}{Example}[section]
\newtheorem{assumption}{Assumption}[section]

\newtheorem{corollary}{Corollary}[section]

\newtheorem{theorem}{Theorem}[section]
\newtheorem{remark}{Remark}[section]

\newcommand{\Ex}{{\mathbb{E}}}
\newcommand{\E}{{\mathbb{E}}}

\newcommand{\var}{{\rm{var}}}

\DeclareMathOperator*{\argmin}{arg\,min}

\newcommand{\R}{\mathbb{R}}

\newcommand{\1}{\mathsf{1}}

\definecolor{cobalt}{rgb}{0.0, 0.28, 0.67}

\title{Statistical inference for boundary constrained parameters}
\author{Jaeseon Lee, Anirban  Bhattacharya and Suhasini Subba Rao}

\begin{document}

\maketitle

\begin{abstract}
Estimation and inference under boundary constraints arises in many areas of
statistics, but the available asymptotic theory remains limited and
fragmented relative to the classical case, where the parameter lies in the
interior of the parameter space. In this paper we derive a
boundary-constrained analogue of the standard asymptotic expansion for a
minimiser of a loss function, in which the deviation of the estimator from the
population parameter is decomposed over lower-dimensional subspaces determined
by the directional derivatives at the minimiser. Using this decomposition we
show that the limiting distribution of the estimator is a mixture of point
masses and truncated skew-normal distributions, with mixture probabilities
determined by the directional derivative conditions. A consequence of this
result is that the constrained and unconstrained coordinates do not separate:
the presence of parameters on or near the boundary affects the sampling
distribution of the unconstrained parameter estimators. We apply the expansion
to obtain the asymptotic distribution of the quasi log-likelihood ratio
statistic under a range of null hypotheses, and since the resulting mixture is
not pivotal in the quasi-likelihood setting, we propose a simulation-based
procedure for approximating critical values. The framework unifies and
simplifies several existing results. We illustrate the methodology by testing
for interactions in exponential family graphical models, where the parameters
are subject to compatibility constraints that ensure a joint distribution exists.

\noindent \emph{Some key words:} Convex optimisation, 
Directional Derivatives, First order expansion,
Quasi log-likelihood ratio test.
\end{abstract}


\section{Introduction}
\label{sec:intro}

Parameter estimation under boundary constraints is a classical problem in
statistics, with applications in fields ranging from biostatistics to
econometrics, dating back to
\cite{chernoff1954distribution,p:mor-71,p:cha-74}; for a book level treatment
see \cite{b:sil-sen-05}, and for a more recent review see
\cite{brazzale2024likelihood}.
The fundamental difference between this setting and the classical one, where
the parameter lies in the interior of the parameter space, is that standard
asymptotic theory no longer applies: if the estimator lies on the boundary, the
derivative of the loss need not be zero, and the expansions underlying
the classical theory are no longer available.
As a consequence, naively applying classical methods in the presence of
boundary constraints can lead to inaccurate inference. Perhaps more
surprisingly, the same is true of parameters which are not themselves
constrained, and for which classical inference might be thought to remain
valid; we return to this point below. Our interest in this area arose from testing interactions in graphical models constructed from
exponential family conditional distributions \citep{p:bes-74,p:yan-15,bhattacharya2026conditionally}, which offer a simple method for
modelling the joint distribution of random vectors. Depending on the
exponential families used, some of these distributions involve constraints on
the parameters in order to ensure existence of a joint distribution. 
Testing for interactions in such models therefore amounts to
a quasi-likelihood ratio test within a boundary-constrained parameter space; see Section \ref{sec:sims} for a concrete example. 


The asymptotic theory under nonstandard parameter constraints was developed in
\cite{chernoff1954distribution,p:sha-85,p:geyer-94,p:and-99}.
The bulk of the subsequent literature has focused on hypothesis testing, and in
particular on the behaviour of the log-likelihood ratio test when parameters of
interest may lie on the boundary, often in the presence of nuisance parameters.
Broadly speaking, this literature divides into two strands, according to
whether the nuisance parameters are themselves boundary restricted.

In the first strand, an early contribution is \cite{p:sel-lia-87}, who study
the log-likelihood ratio test under boundary constraints assuming that the
likelihood is correctly specified. By standardizing the score with the Fisher
information matrix, they orthogonalize the system, at the cost of transforming
the parameter space from a rectangular cone into a more general cone. With this
transformation they derive the limiting distribution of the log-likelihood
ratio test in certain specific examples (when two parameters lie on the
boundary), and show that the limit is a mixture of chi-squares of different
degrees of freedom.
Building on these results, in the setting where none of the nuisance parameters
are boundary restricted, \cite{p:susko-13} derives the distribution of the
log-likelihood ratio test (under a correctly specified likelihood)
\emph{conditioned} on the number in the activation set (that is, the number of
null hypothesis parameters estimated to lie on the boundary).
He shows that, under the null, the conditional test statistic asymptotically
follows a chi-square distribution whose degrees of freedom equal the number of
estimated null hypothesis parameters in the interior.
\cite{p:chen-18} generalize this conditional approach to the quasi-likelihood
setting by using a variant of the vertical scaling method of
\cite{p:chandler-07} to render the conditional distribution pivotal (in the
case with no nuisance parameters).
Recently, \cite{p:hua-20} applied a technique similar to that of
\cite{p:sel-lia-87} to the quasi log-likelihood ratio test, again focusing on
the case where two parameters lie on the boundary.
A related line of work, treated in \cite{b:sil-sen-05} and \cite{p:moh-20},
considers null hypotheses specified by a set rather than a single point. The
resulting distributions are variations of those appearing in \cite{p:susko-13},
and are also covered by the framework developed here.

In a separate strand of research,
\cite{p:and-01,p:kop-sin-11,cavaliere2017consistency,cavaliere2022bootstrap,p:gui-mc-25}
consider the case where the nuisance parameters are not identifiable under the
null, or the nuisance parameters themselves lie on the boundary. The results in
this setting are very different to those cited previously, and without an
appropriate correction the (quasi) likelihood ratio test can have an incorrect
level. Relatedly, \cite{p:gui-mc-25} show that inference for a parameter of
interest is affected by a boundary-constrained nuisance parameter in ways that
go beyond a change in variance or shape; our results give a precise description
of this effect.

To illustrate the gap between the boundary-constrained and classical settings,
suppose that $\mathcal{L}_{n}(\theta)$ is a loss based on a sample of
size $n$, normalized so that $\mathcal{L}_{n}(\theta) = O_{p}(1)$, that
$\widehat{\theta} = \arg\min_{\theta \in \Theta} \mathcal{L}_{n}(\theta)$, and
that $\theta_0$ is the population parameter, which lies in the \emph{interior}
of $\Theta$.
Under suitable regularity conditions it is well known (see, for example,
\cite{van2000asymptotic}) that $\widehat{\theta}$ has the first order expansion
\begin{align}
\label{eq:standard}
(\widehat{\theta}-\theta_0) = -W_{n}^{-1}s_{n} + o_{p}(n^{-1/2}),
\end{align}
where $W_{n} = \mathbb{E}[\nabla^{2} \mathcal{L}_{n}(\theta_0)]$ and
$s_{n} = \nabla \mathcal{L}_{n}(\theta_0)$.
This expansion is the core result underlying both the asymptotic normality of
$\widehat{\theta}$ and the asymptotic theory of the quasi log-likelihood ratio
(with or without nuisance parameters).
As far as we are aware, no analogue of \eqref{eq:standard} has been established
under boundary constraints.
Instead, prior work tends to focus on conditions which ensure that a local
quadratic approximation of the likelihood holds. In this paper we take a
different approach: we assume the quadratic approximation upfront, and use it as a
vehicle for obtaining the sampling properties of the parameter estimators. It
is this change of emphasis, from verifying the approximation to exploiting it,
that allows us to treat the boundary-constrained and interior cases within a
single expansion, and consequently to recover as special cases results which
have until now been obtained separately and by rather different arguments.


Our contributions are as follows. 
In Section \ref{sec:estimator}, 
using some simple tools from
optimisation theory, we derive a boundary-constrained analogue of
\eqref{eq:standard} in which $(\widehat{\theta}-\theta_0)$ is decomposed over
lower-dimensional subspaces of $\R^{p}$, where $p$ denotes the dimension of
$\theta$ and the partition is determined by the directional derivatives at the
minimiser.
Using this decomposition, we show that the limiting distribution of
$\sqrt{n}(\widehat{\theta}-\theta_0)$ is a mixture of point masses and
(truncated) skew-normal distributions, the point mass or skew-normal component
corresponding to whether a given coordinate lies on the boundary or in the
interior, weighted by mixture probabilities that are themselves determined by
the directional derivative conditions. It is interesting to note that, although
the skew-normal appears in many applied contexts, this appears to be the first
time it has arisen as the limiting distribution of an estimator in an
inferential problem of this kind.
In contrast to classical unconstrained inference, we show that the presence of
parameters close to the boundary impacts the distribution of the unconstrained
parameter estimators. This has a direct practical implication, since it means
that the boundary cannot be dealt with by restricting attention to the
constrained coordinates alone: the unconstrained estimators must be treated as
part of the same joint limit, and it is precisely this joint treatment that the
decomposition above makes possible.
In Section \ref{sec:LLRT}, we use this result to obtain the asymptotic distribution of the quasi
log-likelihood ratio statistic under various null hypotheses. We show that the
correctly-specified likelihood result of \cite{p:susko-13} is a special case of
our result; in addition, we obtain explicit expressions for the mixture
probabilities of each $\chi^{2}$ component, which were not derived in
\cite{p:susko-13}.
Since the resulting mixture is not pivotal 
in the quasi-likelihood
setting,  in Section \ref{sec:simtest}
we propose a simulation-based procedure for approximating the critical
values, both for the quasi-likelihood and the conditional quasi-likelihood
versions of the test.
Section \ref{sec:sims} presents simulations
for the graphical model example described above. All proofs are in the
supplementary material.

Finally, we mention that our approach is transparent: the boundary-constrained expansion leads
to a considerably simplified proof of the results in \cite{p:sha-85} and
\cite{p:susko-13}, as well as of the special cases treated in
\cite{p:sel-lia-87}, \cite{p:hua-20} and \cite{p:moh-20}, thus unifying these
previously separate results within a single framework.
We believe the expansions developed here may be useful beyond the settings
treated explicitly, including other constrained M-estimation problems and tests under
inequality restrictions; we leave these extensions to future work.

\section{The sampling properties of estimators that can lie on the
  boundary of the parameter space}\label{sec:estimator}

\subsection{Notation and preliminary results}\label{sec:estimator1}
In this section we introduce notation and some preliminary results; this section can be skipped on first reading.

Let $\mathcal{L}_{n}(\cdot)$ denote a loss function defined on $\mathcal{C} := \R_+^b \otimes \R^{p-b}$ for some 
$b \in [p] :\, = \{1, \ldots, p\}$. We let $B =\{1,\ldots,b\}\subseteq [p]$ denote the set of parameters {\it restricted} to be nonnegative, and 
$U=\{b+1,\ldots,p\}=[p]\backslash B$ denote the {\it unrestricted} coordinates. We assume $\mathcal{L}_{n}$ is smooth, and possesses a 
continuous second derivative on some open set in $\R^p$ containing $\mathcal{C}$. 

 Without loss of generality, we will use this {\it rectangular cone} $\mathcal{C}$ throughout the article as a representative parameter space encapsulating $b$ independent affine constraints.
We assume the same population parameter $\theta_0$ 
uniquely minimises the expected loss for all $n$, i.e., $\theta_0 =
(\theta_{0,B}^\top,\theta_{0,U}^\top)^\top=
\arg\min_{\theta\in \mathcal{C}} \Ex[\mathcal{L}_n(\theta)]$. Uniqueness avoids the identifiability issues considered in \cite{p:and-01}.
Crucially, 
 we allow for the possibility that $\theta_0 \in \partial \mathcal{C}$, where
 $\partial \mathcal{C}$, the {\it boundary} of $\mathcal{C}$, consists of vectors such that at least one entry of 
 $\theta_{B}$ equals zero.

\begin{remark}
More generally, suppose the parameter space is $\{\theta \in \R^p : c_j^\top \theta + f_j \ge 0,\ j=1,\ldots,b\}$, where 
$\{c_j\}_{j=1}^b$ are linearly independent in $\R^p$. Extend $\{c_j\}_{j=1}^b$ to a linearly independent set $\{c_j\}_{j=1}^p$, and define the transformed parameters $\eta_j = c_j^\top \theta + f_j$ for $j \le b$ and $\eta_j = c_j^\top \theta$ for $j > b$. The mapping $\theta \mapsto \eta$ is bijective, and $\eta \in \mathcal{C}$. Other constraint types can be transformed to this form. For example, in isotonic regression with $\theta_1 \le \theta_2 \le \cdots \le \theta_p$, let $\eta_1 = \theta_1$ and $\eta_j = \theta_j - \theta_{j-1}$ for $2 \le j \le p$ (see \cite{b:sil-sen-05}, Chapter 2 and \cite{p:moh-20}).
\end{remark}

We let 
 $\widehat{\theta} = \arg\min_{\theta \in \mathcal{C}} \mathcal{L}_n(\theta)$. Clearly, $\widehat{\theta}$ can lie on $\partial \mathcal{C}$, and this is likely if $\theta_0 \in \partial \mathcal{C}$. One of our goals is to obtain the asymptotic distribution of $\widehat{\theta}$.
  When $\theta_0 \in \partial \mathcal{C}$, standard asymptotic techniques that assume an interior population parameter are not applicable. To address this, we will connect the constrained estimator to its unconstrained counterpart 
  $\widetilde{\theta}\in \arg\min_{\theta\in \widetilde{\Theta}}\mathcal{L}_{n}(\theta)$
  where $\widetilde{\Theta}$ is an open set that contains $\mathcal{C}$ and $\nabla \mathcal{L}_{n}(\widetilde{\theta})=0$.
  Under standard assumptions, the asymptotic sampling properties of the unconstrained estimator $\widetilde{\theta}$ are well known. Therefore proper understanding on this connection can allow us to derive the 
  sampling properties of 
  $\widehat{\theta}$.
The following lemma we make this connection, but in the specific situation that the loss is convex.
\begin{lemma}
\label{lemma:boundary}
Suppose $\mathcal{L}_{n}(\cdot)$ is a strictly convex loss function and let $\widehat{\theta} = \arg\min_{\theta\in \mathcal{C}}\mathcal{L}_{n}(\theta)$ and $\widetilde{\theta} = \arg\min_{\theta\in \R^{p}}\mathcal{L}_{n}(\theta)$.
If there exists $1 \le j \le b$ such that $\widetilde{\theta}_{j} < 0$, then $\widehat{\theta}$ lies on the boundary of $\mathcal{C}$; that is, there exists at least one $1 \le k \le b$ with $\widehat{\theta}_{k} = 0$.
\end{lemma}
 The lemma shows that if $\widetilde{\theta}$ lies outside the feasible region, then $\widehat{\theta}$ must lie on the boundary and the problem reduces to the appropriate subspace. However, the lemma does not quantify how $\widetilde{\theta}$ relates to $\widehat{\theta}$, and the strict convexity requirement on $\mathcal{L}_n$ may not hold in applications. We therefore approximate $\mathcal{L}_n$ locally by a strictly convex quadratic function, which later yields an explicit link between constrained and unconstrained estimators over suitable subspaces (see Section~\ref{sec:explicitexpression}).
 Under the stated conditions, we  
define the score and expected Hessian as
\begin{align}
\label{eq:ahat}
s_{n} = \nabla\mathcal{L}_n(\theta_0), \qquad
W = \Ex[\nabla^2 \mathcal{L}_{n}(\theta_0)].
\end{align}
where $\E[s_n] = 0$ and $W$ is strictly positive definite (these two conditions will be assumed throughout the paper). 
We assume that $\Ex[\nabla^2 \mathcal{L}_{n}(\theta_0)]$ does not depend on $n$ however this condition can be relaxed. 
Using $s_{n}$ and $W$, we define the quadratic approximation
\begin{align}
\label{eq:Qgeneral}
Q_n(\delta) \equiv Q_n(\delta; s_n, W) = \mathcal{L}_n(\theta_0) + \delta^\top s_{n} +
\frac{1}{2} \delta^\top W \delta,
\end{align}
where we set $\delta = \theta-\theta_0$.
Since $\theta\in \mathcal{C}$, we have $\delta \in \mathcal{C}_{\xi}$ where
\begin{align}
\label{eq:parameterspace}
\mathcal{C}_{\xi}  = \bigg[ \prod_{j=1}^b [\xi_j, \infty) \bigg]\otimes \R^{p-b}, \qquad \xi_{j}=-\theta_{0,j} \le 0 \textrm{ for }1\leq j \leq b.
\end{align}
We will view $\xi$ as a $p$-dimensional vector $\xi = (\xi_B^\top, 0_U^\top)^\top$ with $\xi_B = - \theta_{0, B}$. Define
\begin{align}
\label{eq:LQ}
\widehat{\theta}_{L} = \arg\min_{\theta\in \Theta}
\mathcal{L}_{n}(\theta), \qquad
\widehat{\delta}_{Q} = \arg\min_{\delta\in \mathcal{C}_{\xi}}
Q_{n}(\delta) 
\quad\textrm{ and }\quad 
\widehat{\theta}_{Q} = \widehat{\delta}_{Q}+\theta_0,
 \end{align}
where $\Theta = \prod_{j=1}^b [0, M_j] \otimes \prod_{\ell=b+1}^p [M_{1, \ell}, M_{2, \ell}]$ is a compact subset of $\mathcal{C}$ containing $\theta_0$ which is introduced for technical reasons, standard in asymptotic analysis of extremum estimators. 

Central to our results is that the loss function
$\mathcal{L}_{n}(\cdot)$ and corresponding estimators 
defined in (\ref{eq:LQ}) satisfy the following rates
\begin{align}
\label{eq:rate2}
\sup_{\theta\in \mathcal{I}_{n}}|\mathcal{L}_n(\theta) -Q_{n}(\theta-\theta_0)| = o_{p}(n^{-1})
\quad \text{and} \quad
\widehat{\theta}_{L} = \widehat{\theta}_{Q}+o_{p}(n^{-1/2}),
\end{align}
where $\mathcal{I}_{n} = \{\theta:\|\theta-\theta_0\|_2 = O(n^{-1/2})\}$ (implicit in the above is that $\widehat{\theta}_{L} = \theta_{0}+O_{p}(n^{-1/2})$). 
We use (\ref{eq:rate2}) as our fundamental
starting point.  
There exists several articles which 
prove similar conditions to  (\ref{eq:rate2}); see \cite{p:sel-lia-87,p:and-99,p:hua-20,p:gui-mc-25}. 
For completeness, we provide conditions 
for (\ref{eq:rate2}) to hold in the supplement. 

 To reduce notation, from now on
 we drop the suffix $n$ and set
 $s=s_{n}$ and $Q=Q_{n}$. 
 In what follows we consider the minimiser $\widehat{\delta} \equiv \widehat{\delta}_Q$ of $Q(\delta)$ over $\delta \in \mathcal{C}_{\xi}$ and later use (\ref{eq:rate2}) to transfer the results to $\widehat{\theta}_L$. 
 Further, to make the discussion more transparent we will usually refer to the loss function being minimised over $\mathcal{C}$ 
 rather than $\Theta$.

\subsection{Closed form expression for minimiser of the quadratic loss over
  the rectangular cone set}\label{sec:explicitexpression}

In this section we give an explicit expression for the minimiser of the quadratic approximation over the rectangular cone, first in one dimension to build intuition, and then in arbitrary dimensions. We also identify the affine subspace of $\R^{p}$
on which the solution lies via directional derivative conditions, and show how this yields a partition of $\R^p$. 
The only assumption we make in this section
is that $W$ is a strictly positive definite matrix.

Suppose that $\mathcal{L}_{n}(\theta)$ is a univariate loss function that is minimised over the parameter space $[0,\infty)$. The quadratic approximation centred at the population parameter $\theta_0$ is
\[
Q(\delta) = \mathcal{L}_{n}(\theta_0)+ s\,\delta + \tfrac{1}{2}w\,\delta^{2},
\]
where $\delta, s$ and $w$ are univariate quantities, with $w>0$. After centralisation, $Q(\delta)$ is minimised over $\mathcal{C}_{\xi} = [\xi,\infty)$ where $\xi=-\theta_0\leq 0$. The unconstrained minimiser is
\[
\widetilde{\delta} = \argmin_{\delta\in \R}Q(\delta) = -w^{-1}s.
\]
If $\widetilde{\delta}>\xi$, then $\widetilde{\delta}$ also minimises $Q(\delta)$ over $\mathcal{C}_{\xi}$. Otherwise, if $\widetilde{\delta}<\xi$, convexity implies
\[
\argmin_{\delta\in [\xi,\infty)}Q(\delta)=\xi.
\]
Thus, if the unrestricted estimator lies outside the boundary, the restricted minimiser lies on the boundary.

We now extend to arbitrary dimension. For $S \subset [p]$ and $a \in \R^p$, write $a_S = (a_j : j \in S)$ for the sub-vector of $a$ with indices in $S$. Also, for $K \subset S$, $(a_S)_K$ denotes the sub-vector of $a_S$ corresponding to $K$ (so $(a_S)_K = a_K$). We will also write $a = (a_S^\top, a_{S^c}^\top)^\top$ by convention, even if $S$ is not contiguous. When clear from the context, with a slight abuse of notation, we sometimes also write $a = (a_S, a_{S^c})$ to reduce notation. 
For $a, a'\in \R^q$, we write $a<a'$ if $a_j<a'_j$ for all $j\in [q]$. We use $|A|$ to denote the cardinality of the set $A$.


Recall the unconstrained and constrained quadratic problems
\begin{align}
\label{eq:tildedelta}
\widetilde{\delta} = \arg\min_{\delta \in \R^p} Q(\delta) = -W^{-1}s, 
\qquad
\widehat{\delta} = \arg\min_{\delta \in \mathcal{C}_\xi} Q(\delta).
\end{align}
In light of Lemma \ref{lemma:boundary}, when the unrestricted minimiser $\widetilde{\delta}$ lies outside $\mathcal{C}_\xi$, the constrained minimiser $\widehat{\delta}$ lies on a lower-dimensional affine subspace of $\R^p$.
We now characterise $\widehat{\delta}$.

For each $C \subseteq B$, define
\begin{align}
\label{eq:AF-defin}
A := A_{C}= B \setminus C, 
\qquad
F:=F_{C} = C \cup U,
\end{align}
so that $A\cup F=[p]$. Express
$s$ and $W$ in the 
corresponding blocked vector/matrix form
\[
s = \begin{pmatrix} s_{A} \\ s_{F} \end{pmatrix},
\qquad
W = \begin{pmatrix} W_{AA} & W_{AF} \\ W_{AF}^{\top} & W_{FF} \end{pmatrix}.
\]
Using the above, we define $\delta^{(C)}$ (a candidate solution for $\widehat{\delta}$) to be the unconstrained minimizer of $Q$ subject to the entries corresponding to $A$ fixed at $\xi_A$, that is, 
\begin{align}
\label{eq:deltaA}
\delta^{(C)} = \begin{pmatrix} \xi_{A} \\ \delta^{F} \end{pmatrix},
\qquad
\delta^{F} = \argmin_{\delta_{F} \in \R^{|F|}} Q(\xi_{A}, \delta_{F})
= -W_{FF}^{-1}\big(s_{F}+ W_{FA}\xi_{A}\big).
\end{align}
Note the use of superscripts $(C)$ and $F$ on $\delta$ indicates constructed vectors (not subvectors of a larger vector). Borrowing optimisation terminology,
$C$ and $A$ are respectively the 
\emph{inactive} and \emph{active} subsets of the boundary constraints. In particular, $C = B$ corresponds to $A = \emptyset, F = [p]$, and $\delta^F = \widetilde{\delta}$. 
In Remark 2.2 below, we connect $\delta^F$ for any $F$ with $\widetilde{\delta}$.

\begin{remark}\label{remark:SWinv}
Since $s = (s_{A}^{\top},s_{F}^{\top})^\top
= -W (\widetilde{\delta}_{A}^{\top},
\widetilde{\delta}_{F}^{\top})^\top$, we have
$s_{F} = -\big(
W_{FA}\widetilde{\delta}_{A}+W_{FF}\widetilde{\delta}_{F}\big)$. Substituting
this into $\delta^{F}$ gives
\begin{align}
\label{eq:deltaproject}
\delta^{F} =
\widetilde{\delta}_{F}+(W_{FF})^{-1} W_{FA}\,(\widetilde{\delta}_{A} -
\xi_{A}),
\end{align}
linking the unconstrained least squares estimators
$\widetilde{\delta}$ and $\delta^{F}$.
\end{remark}

In the following theorem we show
that $\widehat{\delta}$ equals $\delta^{(C)}$ for some subset $C \subseteq B$, whose identity is determined by $(s, W)$.
We briefly sketch the proof below, which is based on the {\it first-order optimality condition} for constrained minimization of smooth convex functions over convex sub-domains. 
For a smooth function $f$ on a closed set $\mathcal{C}$ and a feasible direction $\nu$ at $x\in \mathcal{C}$ (i.e., $x+t\nu\in \mathcal{C}$ for all sufficiently small $t\ge 0$), the directional derivative of $f$ at $x$ along $\nu$ is
$D_{\nu} f(x) \;=\; \lim_{t\downarrow 0}(f(x+t\nu)-f(x))/t \;=\; \nu^\top \nabla f(x)$.
If $x$ minimises $f$ over $\mathcal{C}$ and $\mathcal{C}$ is convex, then $D_{\nu} f(x)\ge 0$ for all feasible $\nu$;
conversely, for convex $f$, this first-order condition is sufficient for $x$ to be the global minimum.
See Supplement~\ref{sec:optim_results} for details.

Since $\mathcal{C}_{\xi}$ is a convex set, $\delta^{(C)}$ minimises $Q$ over $\mathcal{C}_\xi$ if and only if $D_{\nu}Q(\delta^{(C)}) \ge 0$ for all feasible directions $\nu$ at $\delta^{(C)}$.
At $\delta^{(C)} = (\xi_A, \delta^F)$,
the set of all feasible directions in $\mathcal{C}_{\xi}$ is $\{\nu=(\nu_A, \nu_F):\ \nu_A\in \R^{|A|}_+, \nu_F \in \R^{|F|} \}$.  
Further,
substituting $\delta^{(C)}$ into $\nabla Q$ and using (\ref{eq:deltaproject}) gives
\begin{align}
  \label{eq:schur}
\nabla Q(\delta^{(C)})=
\begin{pmatrix}
-\,W_{A|F}\,(\widetilde{\delta}_A-\xi_A)\\[2pt]
0_F \\
\end{pmatrix},
\quad \textrm{ where }W_{A|F}=W_{AA}-W_{AF}W_{FF}^{-1}W_{FA}.
\end{align}
These two facts lead to the conditions in
Theorem~\ref{theorem:blocks} below. 

\begin{theorem}
\label{theorem:blocks}
Let $Q(\cdot)$ be defined as in (\ref{eq:Qgeneral}),  $\mathcal{C}_\xi$ be as in
(\ref{eq:parameterspace}), and $\mathcal{P}(B)$ denote the power set of $B$. Recall $\widetilde{\delta}$ and $\widehat{\delta}$ from \eqref{eq:tildedelta}. Let $C,A_{C},F_{C}$ be as in \eqref{eq:AF-defin}, and recall from \eqref{eq:deltaA} that $\delta^{(C)} = (\xi_{A_C}^\top, (\delta^{F_C})^\top)^\top$, with $\delta^{F_C}=-W_{F_CF_C}^{-1}(s_{F_C}+W_{F_{C}A_{C}}\xi_{A_C})$. Then
\begin{align}
\label{eq:consequence}
\widehat{\delta}
= \sum_{C\in \mathcal{P}(B)}
\begin{pmatrix}
\xi_{A_{C}} \\
\delta^{F_{C}}
\end{pmatrix}
\underbrace{\,
I\!\left[\delta^{F_C}_{C} > \xi_{C} \right]
\,
I\!\left[
W_{A_{C}|F_{C}} \,(\widetilde{\delta}_{A_{C}}-\xi_{A_C}) \leq 0\right]}_{=I(s\in R_{C})},
\end{align}
where the sets
\[
R_{C}
=\Big\{a\in \R^{p} \,:\, \big(-W_{F_{C}F_{C}}^{-1}\,(a_{F_{C}} +W_{F_CA_C} \xi_{A_C}) \big)_{C}>\xi_{C}, \ 
  -W_{A_{C}|F_{C}}\big[(W^{-1}a)_{A_{C}}+\xi_{A_C}\big]\leq 0\Big\}
\]
form a disjoint partition of $\R^{p}$.
\end{theorem}

To see why $\{R_{C}\}_{C\in \mathcal{P}(B)}$ forms a disjoint partition of $\R^{p}$, define disjoint sets based on the active set of $\widehat{\delta}(a) = \argmin_{\delta \in \mathcal{C}_\xi} Q(\delta; a, W)$ as $a$ varies over $\R^p$. For each $C\in \mathcal{P}(B)$ set
\[
\bar{R}_{C} = \{a\in \R^{p}:\ \widehat{\delta}_{A_C}(a) =\xi_{A_C} \ \text{and} \
\widehat{\delta}_{C}(a) \in \xi_C + (0,\infty)^{|C|}\}.
\]
By construction, $\{\bar{R}_{C}\}_{C\in \mathcal{P}(B)}$ are disjoint, and since $\widehat{\delta}(a)$ is defined for every $a \in \R^p$, we have $\cup_{C\in \mathcal{P}(B)} \bar{R}_{C} = \R^{p}$. In the proof of Theorem \ref{theorem:blocks} we show that $\bar{R}_{C}=R_{C}$.

\begin{example}\label{example:p1}
The simplest case is $p=1$, $B=\{1\}$ and $\xi=0$ (the population parameter lies on the boundary). Then $\mathcal{P}(B) = \{\emptyset,\{1\}\}$, $\mathcal{C}_0 = \R_{+}$ and $U=\emptyset$. If
$\widetilde{\delta}>0$, then $C=\{1\}$ and $\widehat{\delta}=\widetilde{\delta}$. Else, if
$W \widetilde{\delta} < 0$, the positivity of $W$ implies $\widetilde{\delta}<0$, so $C=\emptyset$ and
$\widehat{\delta}=0$.
\end{example}

\begin{example}\label{example:p2}
Consider $p=2$, $B=\{1,2\}$, $U=\emptyset$ and $\xi=0$. Then
$\mathcal{P}(B) = \{\emptyset,\{1\},\{2\},B\}$. Define $\widetilde{\delta}=-W^{-1}s$, and
\[
\delta^{1} = -W_{11}^{-1}s_{1} =
\widetilde{\delta}_{1}+W_{11}^{-1}W_{12}\,\widetilde{\delta}_{2},
\qquad
\delta^{2} = -W_{22}^{-1}s_{2} =
\widetilde{\delta}_{2}+W_{22}^{-1}W_{21}\,\widetilde{\delta}_{1}.
\]
Hence
\[
\delta^{(B)} = \widetilde{\delta}, \quad
\delta^{(1)} =
\begin{pmatrix}
\delta^{1} \\ 0
\end{pmatrix}, \quad
\delta^{(2)} =
\begin{pmatrix}
0 \\ \delta^{2}
\end{pmatrix}, \quad
\delta^{(\emptyset)} =
\begin{pmatrix}
0 \\ 0
\end{pmatrix}.
\]
Since $W_{1|2}$ and $W_{2|1}$ are scalars and $W$ is positive definite, both $W_{1|2}$ and $W_{2|1}$ are strictly positive. Thus: (i) for $C=\{2\}$ and $A=\{1\}$, $W_{1|2}\,\widetilde{\delta}_{1}<0$ is equivalent to $\widetilde{\delta}_{1}<0$; (ii) for $C=\{1\}$ and $A=\{2\}$, $W_{2|1}\,\widetilde{\delta}_{2}<0$ is equivalent to $\widetilde{\delta}_{2}<0$. Finally, when $C=\emptyset$ and $A=\{1,2\}$, $W \widetilde{\delta} = -s$, so $W\widetilde{\delta}<0$ is equivalent to $s_{1}>0$ and $s_{2}>0$, i.e. $\delta^{1}<0$ and $\delta^{2}<0$. Altogether, Theorem \ref{theorem:blocks} yields
\begin{align}\label{eq:casep2}
\begin{pmatrix}
\widehat{\delta}_{1} \\
\widehat{\delta}_{2}
\end{pmatrix}
=
\begin{cases}
\widetilde{\delta} & \text{if }\ \widetilde{\delta}_1>0,\ \widetilde{\delta}_2>0,\\[2pt]
\delta^{(1)} & \text{if }\ \delta^{1}>0,\ \widetilde{\delta}_2<0,\\[2pt]
\delta^{(2)} & \text{if }\ \widetilde{\delta}_1<0,\ \delta^{2}>0,\\[2pt]
(0,0) & \text{if }\ \delta^{1}<0,\ \delta^{2}<0.
\end{cases}
\end{align}
A related expansion is given in \cite{p:hua-20}.
\end{example}

We next rewrite Theorem \ref{theorem:blocks} in terms of $\Sigma = W^{-1}$ using block matrix identities. By the block inversion formula, $W_{A_C|F_C} = \Sigma_{A_C A_C}^{-1}$, and $(W_{F_C F_C})^{-1} W_{F_C A_C} = -\Sigma_{F_CA_C} (\Sigma_{A_CA_C})^{-1}$.

\begin{corollary}\label{cor:blocks}
Under the conditions of Theorem \ref{theorem:blocks}, with $\Sigma = W^{-1}$, we can write
\[
\widehat{\delta}
= \sum_{C\in \mathcal{P}(B)}
\begin{pmatrix}
\xi_{A_{C}} \\
\delta^{F_{C}}
\end{pmatrix}
\,
I\!\left[\big(\delta^{F_C}\big)_{C} > \xi_{C} \right]
\,
I\!\left[
\Sigma_{A_{C}A_{C}}^{-1} \,(\widetilde{\delta}_{A_{C}}-\xi_{A_C}) \leq 0\right],
\]
where
\[
\delta^{F_{C}} =
\widetilde{\delta}_{F_{C}}-\Sigma_{F_CA_C} (\Sigma_{A_CA_C})^{-1}
(\widetilde{\delta}_{A_C}-\xi_{A_C}).
\]
\end{corollary}

\begin{example}\label{ex:SL87_case2and3}
We show that Case 2 in \cite{p:sel-lia-87}, with
$p=2$, $\xi=0$ and $\mathcal{C}_0 = \R^{2}_{+}$, is a special case of
Corollary \ref{cor:blocks}. If $\Sigma = W^{-1}$, then
$\widetilde{\delta} = -\Sigma s$, and
\[
\delta^{1} = \widetilde{\delta}_{1} -
\Sigma_{12}\Sigma_{22}^{-1}\,\widetilde{\delta}_{2},
\qquad
\delta^{2} = \widetilde{\delta}_{2} -
\Sigma_{21}\Sigma_{11}^{-1}\,\widetilde{\delta}_{1}.
\]
Substituting into (\ref{eq:casep2}) recovers Case 2 in \cite{p:sel-lia-87}.

It is interesting to note that Case 3 in \cite{p:sel-lia-87} appears to conjecture
that there exist matrices $\{L_{C}\}_{C\in \mathcal{P}(B)}$ which give the decomposition 
\begin{align*}
  \widehat{\delta}  =
  \sum_{C\in \mathcal{P}(B)}
  \left(
  \begin{matrix}
    0_{A_{C}} \\
   \widetilde{\delta}_{F_{C}} - \Sigma_{F_{C}A_{C}}\Sigma_{A_{C}A_{C}}^{-1}\widetilde{\delta}_{A_{C}}
   \end{matrix} 
\right)I\left[L_{C}\widetilde{\delta}>0\right].
\end{align*}
This is similar in spirit to Corollary \ref{cor:blocks}, when $\xi = 0$, though no
expression for $L_{C}$ is given.
\end{example}

\subsection{Asymptotic sampling properties of the boundary estimators}\label{sec:distribution}

We now return to the minimiser $\widehat{\theta}_{L}$ of the loss function
$\mathcal{L}_{n}(\theta)$ over
$\Theta$. Recall the  population parameter is $\theta_{0} = (\theta_{0,B},\theta_{0,U})$. Also recall that $\widehat{\delta} = \argmin_{\delta\in
    \mathcal{C}_{\xi}}Q(\delta)$ where 
$\xi = (\xi_{B},0_{U})$ and   $\xi_{B} =
-\theta_{0,B}$.

In the following theorem we obtain an analogous expansion to
(\ref{eq:standard}) under boundary constraints by 
combining Theorem \ref{theorem:blocks} with (\ref{eq:rate2}). 
\begin{theorem}\label{theorem:nonstandard}
 Suppose the rates in (\ref{eq:rate2}) hold. Recall the 
 notation from Theorem \ref{theorem:blocks}.
Then we have the first order expansion 
\begin{align}
\label{eq:Lest}
    \widehat{\theta}_{L}-\theta_0
     &= \sum_{C\in \mathcal{P}(B)}
         \left(
       \begin{matrix}
         \xi_{A_{C}} \\
         \delta^{F_{C}}\\
       \end{matrix}  
       \right)
  I\left[\delta^{F_C}_{C} > \xi_{C} \right]
        I\left[
        W_{A_{C}|F_{C}}(\widetilde{\delta}_{A_{C}}-\xi_{A_{C}}) \leq
  0\right] + o_{p}(n^{-1/2})
\end{align}
where $\delta^{F_C}
=-W_{F_{C}F_{C}}^{-1}(s_{F_C}+W_{F_CA_{C}}\xi_{A_{C}})$.
\end{theorem}

The above result is key to obtaining the the asymptotic sampling
distribution of $\widehat{\theta}_{L}-\theta_0$. To that end, we first introduce the following
definition of a skewed normal and truncated skewed normal distribution. 

\begin{definition}[The generalised multivariate skewed and
  truncated skewed normal]
  \label{def:SN_Variants}
  Let $\phi_{d}(\cdot;\Sigma)$ and $\Phi_{d}(\cdot; \Sigma)$ 
  denote the probability density function and cumulative distribution function of a $N_d(0, \Sigma)$ distribution. For $Y \sim N_d(0, \Sigma)$ and $H \subseteq \R^d$, we write $\operatorname{pr}(Y \in H; \Sigma) = \int_H \phi_d(y; \Sigma) dy$. 
  Thus, $\Phi_d(y; \Sigma) = \operatorname{pr}(Y \le y; \Sigma)$; also denote $\widetilde{\Phi}_{d}(y;\Sigma) = \operatorname{pr}(Y>y;\Sigma)$. 
  Further, if $(Y_{1}^{\top},Y_{2}^{\top})\sim N_{d}(0,\Sigma)$, then 
  we denote $\Phi^{*}_{d}(y_{1},y_{2};\Sigma) = \operatorname{pr}(Y_{1}\leq y_{1},Y_{2}>y_{2};\Sigma)$.

  Let $d_1, d_2, d_3, d_4$ be positive integers with $d_1 + d_2 = d_4$. 
  Let $E_i \in \R^{d_i \times d_i}$ for $i = 1, \ldots, 4$ be positive definite matrices. Also, 
  let $B_{12} \in \R^{d_1 \times d_2}, B_{13} \in \R^{d_1 \times d_3}, B_{23} \in \R^{d_2 \times d_3}$ be real matrices. Finally, let $\zeta = (\zeta_1, \zeta_2) \in \R^{d_4}$ with $\zeta_i \in \R^{d_i}$ for $i = 1, 2$.
\\[1ex]
\noindent {\bf One-sided Skewed Normal distribution.} 
A random variable $Y \in \R^{d_2}$ follows a one-sided skewed normal distribution, denoted $Y \sim \widetilde{SN}_{d_{1},d_{2}}(\zeta_1,B_{12},E_2, E_1)$, if its density at any $y \in \R^{d_2}$ is given by 
\begin{align*}
\widetilde{SN}_{d_{1},d_{2}}(y; \zeta_1, B_{12},E_2, E_1) = 
\mathcal{Z}^{-1} \, \widetilde{\Phi}_{d_1}(\zeta_1 -B_{12}y;E_1) \, \phi_{d_2}(y;E_2),
\end{align*}
with the normalizing constant $\mathcal{Z} = \int_{\R^{d_2}} \widetilde{\Phi}_{d_1}(\zeta_1-B_{12}x;E_1) \, \phi_{d_2}(x;E_2) dx$. Observe that if $B_{12} = 0$, this reduces to a $N_{d_2}(0, E_2)$ distribution.
\\[1ex]
\noindent {\bf Two-sided Skewed Normal distribution.} 
A random variable $Y \in \R^{d_3}$ follows a two-sided skewed normal distribution, denoted $Y \sim SN_{d_4,d_3}(\zeta_1,\zeta_2,B_{13},B_{23}, E_3, E_4)$, if its density at any $y \in \R^{d_3}$ 
is given by 
\begin{align*}
SN_{d_4, d_3}(y; \zeta_{1},\zeta_{2},B_{13},B_{23},E_3, E_4) :\, = 
\mathcal{Z}^{-1} \, \Phi^{*}_{d_4}(\zeta_{1}-B_{13}y, \zeta_{2}-B_{23}y; E_4) \, \phi_{d_3}(y;E_3),
\end{align*}
with the normalizing constant $\mathcal{Z} = \int_{\R^{d_3}} \Phi^{*}_{d_4}(\zeta_{1}-B_{13}x, \zeta_{2}-B_{23}x; E_4) \, \phi_{d_3}(x; E_3) \, dx$. 
\\[1ex]
\noindent {\bf Truncated Skewed Normal distribution.} A random variable $Y \in \R^{d_2}$ follows a truncated skewed normal distribution, denoted 
$Y \sim tSN_{d_1,d_2}(\zeta_1,\zeta_2,B_{12},E_2,E_1)$, if its density at any $y \in \R^{d_2}$ is given by 
\begin{align*}
tSN_{d_1,d_2}(y; \zeta_1,\zeta_2,B_{12},E_2,E_1) :\,= \mathcal{Z}^{-1} \, \Phi_{d_1}(\zeta_1-B_{12}y;E_1) \, \phi_{d_2}(y;E_2) \, I(y > \zeta_2),
\end{align*}
with the normalizing constant $\mathcal{Z} = \int_{\R^{d_2}} \Phi_{d_1}(\zeta_1-B_{12}x;E_1) \, \phi_{d_2}(x;E_2) \, I(x > \zeta_2) \,dx$. If $B_{12} = 0$, this reduces to a truncated multivariate normal distribution with density proportional to $\phi_{d_2}(y; E_2) I(y > \zeta_2)$. Moreover, if the indicator is removed, then this reduces to a multivariate skewed normal distribution. 
\end{definition}

\begin{remark}
\label{remark:SN_defn}
    Several formulations of multivariate skewed normal distributions appear in the literature; see \cite{b:azzalini-13} Chapters~5 and 7 for a comprehensive review. Their connection with the distributions appearing in Definition \ref{def:SN_Variants} is discussed in detail in Remark \ref{remark:supSN} of the supplement. We introduce these specific distributions above because of their appearance inside our subsequent characterization of the limiting distribution of $\sqrt{n}(\widehat{\theta}_L - \theta_0)$. 
\end{remark}


The following lemma connects the distributions introduced in Definition \ref{def:SN_Variants} to conditional distributions of multivariate Gaussians subject to linear inequality constraints; related constrictions of multivariate skewed normal distributions appear in \cite{azzalini1996multivariate, arellano2006unification}. 
\begin{lemma} \label{lemma:SN_TSN}
Suppose $Y_{1},Y_{2},Y_{3}$ are multivariate random vectors 
of dimension $d_1,d_2,d_3$ respectively. Let $d = d_1 + d_2 + d_3$, and assume 
\begin{align*}
\left(
\begin{matrix}
Y_{1} \\
Y_{2} \\
Y_{3} \\
\end{matrix}
\right)\sim N_{d}\left(
0,V\right), \textrm{ where }
V = \left(
\begin{matrix}
V_{11} & V_{12} & V_{13} \\
V_{21} & V_{22} & V_{23} \\
V_{31} & V_{32} & V_{33} \\
\end{matrix}
\right).
\end{align*}
Let $V_{1|2} = V_{11} - V_{12} V_{22}^{-1}V_{21}$ 
and $V_{(1,2)|3}=
V_{(1,2),(1,2)} - (V_{13}^\top, V_{23}^\top)^\top V_{33}^{-1} (V_{31},
V_{32})$, where $V_{(1,2),(1,2)}$ denote the $(1,2)$ sub-diagonal of the matrix of $V$.
Then
\begin{itemize}
    \item[(i)] $Y_{2} \mid (Y_{1}<\zeta_1,Y_{2}>\zeta_2)\sim 
tSN_{d_1,d_2}(\zeta_1,\zeta_2,V_{12}V_{22}^{-1}, V_{22}, V_{1|2})$.

    In the case $V_{12} = 0$, then
$Y_{2}\mid (Y_{1}<\zeta_1,Y_{2}>\zeta_2)\sim 
tN_{d_2}(\zeta_2, V_{22})$ 
    \item[(ii)] $Y_{3}\mid (Y_{1}<\zeta_1,Y_{2}>\zeta_2)\sim
      SN_{d_1+d_2,d_3}( \zeta_1,\zeta_2,
      V_{13}V_{33}^{-1}, V_{23}V_{33}^{-1}, V_{33}, V_{(1,2)|3})$.

In the case $V_{12}=V_{13}=0$, then $Y_3\mid(Y_1 < \zeta_1, Y_2>\zeta_2) \sim \widetilde{SN}_{d_2,d_3}(\zeta_2, V_{23}V_{33}^{-1},V_{33}, V_{2|3})$.
\end{itemize}
\end{lemma}
Lemma \ref{lemma:SN_TSN} is instrumental in showing that the asymptotic distribution of $\sqrt{n}(\widehat{\theta}_L - \theta_0)$ is a mixture of skewed normals and truncated skewed normals in Theorem \ref{theorem:CLT} below. For that, we first state an assumption regarding the unconstrained minimizer $\widetilde{\delta} = - W^{-1} s$ in \eqref{eq:tildedelta}.  

\begin{assumption}\label{assum:CLT}
Assume that
(a) $\Ex[s] = 0$ and $\|\var[s]\|_{\operatorname{spec}} = O(n^{-1})$ (where $\|\cdot\|_{\operatorname{spec}}$ denotes the spectral norm) and 
(b) $\sqrt{n} \widetilde{\delta} \Dcon N(0,\Sigma) \, \text{as }n\to\infty$
    with $\Sigma$ positive definite. 
\end{assumption}

Analogous to the classical least squares estimator, 
we show that Assumption \ref{assum:CLT}(a) yields a root-$n$ rate of convergence for the constrained minimizer
$\widehat{\delta}$ (this together with (\ref{eq:rate2}) gives a root-$n$ rate for $\widehat{\theta}_{L}$).
\begin{lemma}\label{lemma:LSrateG}
Suppose Assumption \ref{assum:CLT}(a) holds. Then $\widehat{\delta} = O_p(n^{-1/2})$.
\end{lemma}
In the following theorem we 
obtain the asymptotic distribution of 
$\widehat{\theta}_{L}$ in the case that
$\xi = 0$ (the general case is considered in 
Section \ref{sec:closetoboundary} and 
Lemma \ref{lemma:finite_approx} of the supplement), 
and 
show that it is a mixture of skewed normals and truncated skewed
normal distributions. Specifically, we establish that the limiting distribution $P$ of $\sqrt{n}(\widehat{\theta}_L - \theta_0)$ is a mixture distribution of the form 
\begin{align}\label{eq:asymp_distr_generic}
P = \sum_{C \in \mathcal{P}(B)} \pi_C \, P^C, \quad P^C = \left(\begin{matrix}
      m_{A_C}(0)\\
    f^{F_C}
\end{matrix}\right), \quad f^{F_C} :\, = \left(\begin{matrix} f^{F_C}_C \\ f^{F_C}_U  \end{matrix}\right),
\end{align}
where the $\pi_C$ are mixture weights (i.e. $\pi_C \ge 0$ and $\sum_{C \in \mathcal{P}(B)} \pi_C = 1$), $m_{A_{C}}(0)$ denotes the distribution on $\R^{|A_C|}$ degenerate at zero, and $f^{F_C}$ denotes a joint distribution on $\R^{|F_C|}$ of $(Y_C, Y_U)$ with $Y_C \sim f^{F_C}_C$ and $Y_U \sim f^{F_C}_U$. We will use this notational convention, which suppresses the dependence between $Y_C$ and $Y_U$ for brevity, through the rest of this paper.
We offer precise characterisations of the constituent pieces $\pi_C$ and $f^{F_C}$ below. As far as we are aware this is the first characterisation of 
the general asymptotic distribution of boundary restricted estimators.
\begin{theorem}\label{theorem:CLT}
  Suppose (\ref{eq:rate2}) and Assumption \ref{assum:CLT} hold, and recall the notation in Theorem \ref{theorem:blocks}. Assume $\xi = 0$ (in other words, $\theta_{0,B}=0$). 
    For each $C\in \mathcal{P}(B)$ we define the matrix
    $V_{C} =  T_{C} \Sigma T_{C}^\top$ where 
    \begin{align*}
        T_{C} &= 
        \left( \begin{matrix}
            W_{A_C|F_C} & 0 & 0\\
            W_{C|U}^{-1}W_{CA_C|U} & I_{|C|} & 0 \\
         W_{U|C}^{-1}W_{UA_C|C}   & 0 & I_{|U|}
        \end{matrix}\right),
    \end{align*}
with    $W_{CA_C|U} = W_{CA_C}-W_{CU}W_{UU}^{-1}W_{UA_C}$,
$W_{UA_C|C}=W_{UA_C}-W_{UC}W_{CC}^{-1}W_{CA_C}$; and $\Sigma$ as in Assumption \ref{assum:CLT}. 
Then, we have that 
\begin{align}\label{eq:TCdel}
\sqrt{n} ((W_{A_C|F_C}\widetilde{\delta}_{A_C})^{\top}, (\delta_{C}^{F_C})^{\top},(\delta_{U}^{F_C})^{\top})^\top = \sqrt{n} \, T_C \widetilde{\delta} \Dcon N(0,V_{C}).
\end{align}
Using this, we have $\sqrt{n}(\widehat{\theta}_{L}-\theta_0)\Dcon P$ with $P$ as in \eqref{eq:asymp_distr_generic}, where 
$\pi_C = \operatorname{pr} \big( \delta^{F_C}_C > 0, W_{A_{C}|F_{C}} \widetilde{\delta}_{A_C} \leq 0  ; \Sigma \big)$ 
and $f^{F_C}$ is the limiting law of $\sqrt{n} \, ( (\delta_C^{F_C})^\top, (\delta_U^{F_C})^\top)^\top \mid \delta_C^{F_C} > 0, W_{A_C|F_C} \widetilde{\delta}_{A_C} \le 0$, induced from \eqref{eq:TCdel}. In particular, recalling the tSN and SN distributions defined in Definition \ref{def:SN_Variants},
\begin{align*}
& f^{F_C}_C \equiv tSN_{|A_C|,|C|}(0,0,  V_{A_C C}V_{CC}^{-1}, V_{C C}, V_{A_{C}|C}), 
\\
& f^{F_C}_{U} \equiv SN_{b, p-b}\left(0, 0,
        V_{A_C U}V_{UU}^{-1}, V_{CU}V_{UU}^{-1}, V_{UU}, V_{B|U} \right).
\end{align*}


In the case that $U=\emptyset$, 
$T_{C}$ is a $2\times 2$ block matrix
and $f^{F_C} = f^{F_C}_C$.
\end{theorem}  
We provide some intuition behind the theorem. From Theorem \ref{theorem:nonstandard}, $\sqrt{n}(\widehat{\theta}_{L}-\theta_0) = \sqrt{n} \widehat{\delta} + o_p(1)$. When $\xi = 0$, from Theorem \ref{theorem:blocks}, $\widehat{\delta} = (0^\top, (\delta^{F_C}_C)^\top)^\top$ if $\delta_C^{F_C} > 0$ and $W_{A_C|F_C} \widetilde{\delta}_{A_C} \le 0$. Also, when $\xi = 0$, $\delta^F = \widetilde{\delta}_F + W_{F_C F_C}^{-1} W_{F_C A_C} \widetilde{\delta}_A$ from \eqref{eq:consequence}. Thus, under the limiting $N(0, \Sigma)$ law of $\sqrt{n} \widetilde{\delta}$, the mixture weight $\pi_C = \operatorname{pr} \big( \delta^{F_C}_C > 0, W_{A_{C}|F_{C}} \widetilde{\delta}_{A_C} \leq 0; \Sigma \big)$, and $f^{F_C}$ is the induced law of $\sqrt{n} \, \delta^{F_C} \mid \delta_C^{F_C} > 0, W_{A_C|F_C} \widetilde{\delta}_{A_C} \le 0$, which from \eqref{eq:TCdel} is recognized as a {\it lower-dimensional marginal} of a truncated multivariate normal distribution, subject to linear inequality constraints. Such a distribution is generally intractable; see \cite{p:cartinhour1990one} and \cite[Proposition 1]{p:zhou2024mass} for some special cases. The marginals $f^{F_C}_C$ and $f^{F_C}_U$ are simpler to describe, with $f^{F_C}_C$ a truncated skewed normal and $f^{F_C}_U$ a two-sided skewed normal as in Definition \ref{def:SN_Variants}. 

In particular, the limiting distribution of $\sqrt{n} \, \widehat{\delta}_U$ is $\sum_{C} \pi_C f^{F_C}_U$. Thus, an interesting aspect of the above theorem is that the {\it marginal distribution} of the unconstrained 
parameter estimators 
$\widehat{\delta}_{U}$ depends on the constrained parameter estimators -- not just through a variance matrix (as in classical asymptotic theory) -- but through mixtures of skewed normal distributions which depend on the number of boundary constrained parameters. We illustrate this in Figure \ref{fig:marginal} below.

\begin{remark}[Sanity Check]
To check the validity of the 
partition probabilities in 
Theorem \ref{theorem:CLT}
we run a numerical simulation in the case $p=3$, $B=\{1,2,3\}$ 
and $U=\emptyset$ (an explicit expression for the  asymptotic distribution of  $\widehat{\theta}_{L} - \theta_{0}$ in the case $p=3$, $B=\{1,2,3\}$ is given in the Supplement). 
We first obtain the empirical distribution of 
$\widehat{\delta}=\arg\min_{\delta \in \mathcal{C}}Q(\delta)$
(evaluated with an optimization routine). With this we calculate the proportion
of times $\widehat{\delta}$ lies within each $\R^{q}$-dimensional
subspace of $\R^{3}$.  This proportion is compared to the probabilities 
$\operatorname{pr}(\delta^{F}_{C}>0, W_{A|F}\widetilde{\delta}_{A}\leq 0)$ (as stated in the  theorem). Despite no explicit expression existing for these probabilities,  since 
$\widetilde{\delta}$ are normal, we 
approximate it with simulation. In the
simulations we set 
\begin{align*}
    W = \left( \begin{matrix}
        1 & 0.5 & 0.3 \\
        0.5 & 1 & 0.4 \\
        0.3 & 0.4 & 1
    \end{matrix}\right),
    \qquad
    \Sigma = \left( \begin{matrix}
        1 & 0.6 & 0.2 \\
        0.6 & 1 & 0.5 \\
        0.2 & 0.5 & 1
    \end{matrix}\right).
\end{align*}
For each replication, we draw
from $\widetilde{\delta} \sim N(0, \Sigma)$. 
 Table~\ref{tab:mixture_prob} gives the result based on $10,000$ replications.
We observe a close alignment 
between the proportion (the difference is due to sampling variability)
in the simulations and the probabilities in the theorem. Thus giving a
numerical justification to Theorem \ref{theorem:CLT}.

\begin{table}[ht]
\centering
\caption{Comparison of simulated frequencies and calculated mixture probabilities.}
\label{tab:mixture_prob}
\begin{tabular}{|c|ccc||c|ccc|}
\hline
$C$ & Simulation & Probability & Difference & $C$ & Simulation & Probability & Difference\\
\hline
$\emptyset$     & 0.3952 & 0.3878 &  0.0074 & $\{1,2\}$       & 0.0860 & 0.0917 & -0.0057 \\
$\{1\}$         & 0.0675 & 0.0685 & -0.0010 &$\{1,3\}$       & 0.0311 & 0.0328 & -0.0017 \\
$\{2\}$         & 0.0241 & 0.0239 &  0.0002 &$\{2,3\}$       & 0.0703 & 0.0707 & -0.0004 \\
$\{3\}$         & 0.0883 & 0.0908 & -0.0025 &$\{1,2,3\}$     & 0.2375 & 0.2339 &  0.0036 \\
\hline
\end{tabular}
\end{table}

\end{remark}

We now apply the above result to the special case that $W = \Sigma^{-1}$, which corresponds to the case the loss is equal to the negative log-likelihood.
\begin{corollary}\label{theorem:CLT3}
  Suppose the conditions in Theorem \ref{theorem:CLT} hold and
  additionally assume that $W = \Sigma^{-1}$ with $\Sigma$ as in Assumption \ref{assum:CLT}. Then $V_C$ in Theorem \ref{theorem:CLT} simplifies to 
 \begin{align*}
      V_{C} &=
\left(
  \begin{matrix}
    \Sigma_{A_C A_C}^{-1} & 0 & 0\\
 0 &   \Sigma_{C|A_C} & \Sigma_{CU|A_C} \\
 0 &  \Sigma_{CU|A_C}^\top  & \Sigma_{U|A_C} \\
   \end{matrix} 
             \right),
    \end{align*}
    with
    $\Sigma_{C|A_C} = \Sigma_{CC} -
    \Sigma_{CA_C}\Sigma_{A_CA_C}^{-1}\Sigma_{A_CC}$,
    $\Sigma_{U|A_C} = \Sigma_{UU} -
    \Sigma_{UA_C}\Sigma_{A_CA_C}^{-1}\Sigma_{A_CU}$ and
    $\Sigma_{CU|A_C} = \Sigma_{CU} -
\Sigma_{CA_C}\Sigma_{A_CA_C}^{-1}\Sigma_{A_CU}$.
    Further, 
    $\pi_C = \operatorname{pr}\big( \delta^{F_C}_C > 0; \Sigma\big)
    \operatorname{pr}\big( \Sigma_{A_C A_{C}}^{-1}\widetilde{\delta}_{A_C} \leq 0; \Sigma\big)$, 
\begin{align*}
f_C \equiv tN_{|C|}(0, \Sigma_{C|A_C}), \quad f_U \equiv \widetilde{SN}_{|C|, p-b}\big(0, \Sigma_{CU|A_C} \Sigma_{U|A_C}^{-1},
    \Sigma_{U|A_C}, \Sigma_{C|U,A_{C}} \big),
\end{align*}
where $\Sigma_{C|U,A_{C}}=\Sigma_{C|A_C}- \Sigma_{CU|A_C}  \Sigma_{U|A_C}^{-1} \Sigma_{UC|A_C}$.

\end{corollary}

\begin{remark}
Comparing Theorem \ref{theorem:CLT} and Corollary \ref{theorem:CLT3}, certain simplifications ensue when $W = \Sigma^{-1}$. 
These are due to the fact that for a given set $C$,
$\widetilde{\delta}_{A_{C}}$ is asymptotically independent of 
$\delta^{F_C}$. Consequently,  the distribution of $\delta^{F_{C}}_{C}$ is not a truncated skewed normal, just a truncated normal with variance $\Sigma_{C|A_{C}}$.  
This simplifies the distributions of some of the log-likelihood ratio tests considered in Section \ref{sec:LLRT}.
\end{remark}
In Supplement \ref{sup:examples} we illustrate the results in this section 
with examples for the cases $p=1,2,3$.

\subsection{Case that some of the parameters do not lie on the boundary} 
\label{sec:closetoboundary}
Theorem \ref{theorem:CLT} and Corollary \ref{theorem:CLT3} assume all the boundary restricted population parameters,
$\theta_{0,B}$, do lie on the boundary. We will use these in Section \ref{sec:LLRT} 
when we derive the distribution of the quasi-likelihood ratio tests under the null hypothesis that the parameters of interest lie on the boundary. However, in the case that some of the
boundary restricted parameters do not lie on the boundary, i.e.
if $\xi_{E}< 0$ 
($\theta_{0,E}>0$)
for some 
$E\subseteq B$,
then for large sample sizes the distribution of  
$\widehat{\delta}_{0,E}$
will behave as if it were unrestricted. However, if the population parameter is not on the boundary but close to it (of order $O(n^{-1/2})$), the boundary will have an impact.
To illustrate this, consider the case $p=1$ and 
$B=\{1\}$, where $\xi=-\theta_{0}<0$. If $\widetilde{\delta}<\xi$ then
$\widehat{\delta} = \xi$, else $\widehat{\delta}=\widetilde{\delta}$.
Consequently, if 
$\widetilde{\delta}\approx N(0,\sigma^2/n)$ then the distribution of $\widehat{\theta}_{L}-\theta_{0}$ is 
a mixture of a point mass and a truncated normal density
\begin{align*}
\widehat{\theta}_{L}-\theta_{0}\approxD 
 m(\xi)\Phi\left(\xi;\frac{\sigma^{2}}{n}\right)
+tN\left(\xi;\frac{\sigma^{2}}{n}\right)
\left[1-\Phi\left(\xi;\frac{\sigma^{2}}{n}\right)\right]
\end{align*}
where $m(\xi)$ is the point mass centered at $\xi$.
As the sample size grows the probability 
$\Phi\left(\xi;\frac{\sigma^{2}}{n}\right)\rightarrow 0$
and we obtain the classical large sample approximation
$\sqrt{n}(\widehat{\theta}_{L}-\theta_0)\Dcon N(0,\sigma^2)$. We now extend this argument to the case $p>1$. 

To simplify notation and make the discussion more transparent, we consider the case that none of the boundary constrained parameters
$\theta_{0,B}$ actually lie on the boundary but in the interior.
Thus $\xi_{j}=-\theta_{0,j}<0$ for all $j\in B$. 
Asymptotically, the distribution of  $\sqrt{n}(\widehat{\theta}_L-\theta_{0})$ will be Gaussian. However, 
for small sample sizes the boundary will have an impact not just on the distribution of 
$\widehat{\theta}_{B}$ but also on 
the unconstrained parameter estimators 
$\widehat{\theta}_{U}$. 
We state the more refined finite sample approximation below  
\begin{align} \label{eq:finite_approx}
   \widehat{\theta}_{L}-\theta_0 &\approxD 
   \sum_{C\in \mathcal{P}(B)} 
 \left[\left(
    \begin{matrix}
      m_{A_{C}}(\xi_{A_C})\\
        tSN_{|A_C|,|C|}(\zeta_{A_C},\zeta_{C}, V_{A_C C}V_{CC}^{-1}, V_{C C}/n,
        V_{A_C|C}/n)  \\
     SN_{b, p-b}\left(\zeta_{A_C},\zeta_{C},
     V_{A_CU}V_{UU}^{-1}, V_{CU}V_{UU}^{-1}, V_{UU}/n, V_{B|U}/n \right)  
      \end{matrix}
  \right) -
  \left(
  \begin{matrix} 
  0 \\
  \zeta_C-\xi_C\\
  \zeta_U\\
  \end{matrix}
  \right)\right]\nonumber\\
&\times \operatorname{pr}\big( \delta^{F_C}_{C} > \xi_{C},W_{A_{C}|F_{C}} (\widetilde{\delta}_{A_C}-\xi_{A_C}) \leq 0 \big)
\end{align} 
where 
$\zeta_{A_{C}} = (T_{C}\xi)_{A_C},\quad 
\zeta_{C} = (T_{C}\xi)_{C}, \quad
\zeta_U = (T_C \xi)_U$, 
and $m_{A_C}(\xi_{A_C})$ is the point mass centered at $\xi_{A_C}$
 (a proof of can be found in Lemma \ref{lemma:finite_approx} of the supplement). 
We observe that the distribution of 
$\widehat{\theta}_{U}$ is a mixture of shifted skew-normals, each determined by $\theta_{B}$ (the distance of the population parameter from the boundary). In contrast, 
if there is no boundary constraint, then 
$\theta_{B}$ only impacts the
distribution of $\widehat{\theta}_{U}$ through the variance of $\widehat{\theta}_{U}$. 
This closely aligns with the observations of
\cite{p:gui-mc-25} who shows that the 
log-likelihood ratio test in the presence of  nuisance parameters that lie on or close to the boundary impacts the critical values of the test statistic. 

To illustrate how 
$\xi$ affects the marginal distribution of $\widehat{\delta}_U$, 
we consider the case $p=2$, $B=\{1\}$, $U=\{2\}$ 
and $\widetilde{\delta}\sim N(0,\Sigma)$ 
(the actual values of $W$ and $\Sigma$ are given in Example 2 
Section \ref{sec:Boundary} in the supplement). 
In Figure \ref{fig:marginal}
we plot the marginal distribution of 
$\widehat{\delta}_{2}$ 
when
$\xi_{1} \in \{0,-0.1,-0.2,-0.5,-1, -2\}$.
We observe that the {\it shape} of 
the distribution changes with $\xi$; the closer $\xi$ to zero, the more skewed the marginal of $\widehat{\delta}_{2}$, with skewness ranging from $0.170$ to $0.619$. There is also bimodality, as it is a mixture of two shifted skewed normals. This example clearly illustrates that using a normal distribution to construct confidence intervals for estimators where one estimator can lie on or close to the boundary can lead to incorrect coverage.
Curiously, the different distributions all appear to be very close in the tails. This is because the tails of a skewed normal resemble those of a
normal. 

\begin{figure}[h!]
\begin{center}
        \includegraphics[height=6cm]{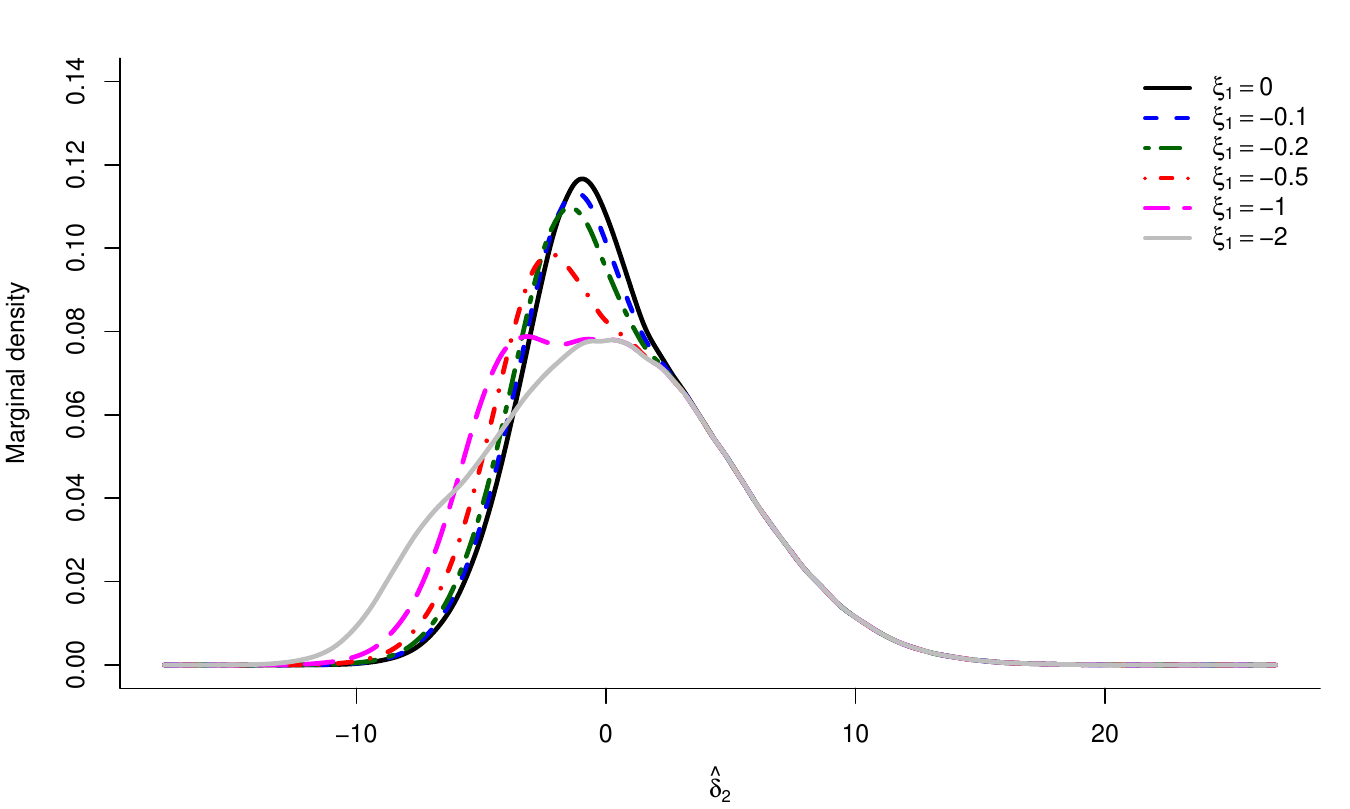}
        \caption{\label{fig:marginal}Marginal distribution of $\widehat{\delta}_{2}$ for different $\xi_1$}

        \smallskip
        \noindent \textbf{Alt text:}
        Six curves show the marginal density of $\widehat{\delta}_2$ for different values of $\xi_1$.
        As $\xi_1$ is close to zero, the density becomes more sharply peaked and skewed. The left tails differ noticeably across the distributions, while the right tails nearly overlap.
\end{center}        
       \end{figure}

\section{The quasi log-likelihood ratio test in the presence of
  nuisance parameters}\label{sec:LLRT}
In this section we use the results in Section \ref{sec:estimator} to
derive the distribution of the quasi log-likelihood ratio test (quasi-LR) in the presence of nuisance
parameters. Our results  generalise and unify several 
disparate results in the literature. Further, the explicit expressions that we
derive for the quasi-LR under the different 
nulls can be used to obtain simulation based 
critical values (see Section \ref{sec:simtest}).

Let $\theta = (\theta_B, \theta_U)$ where $\theta_B \in [0,\infty)^{b}$
and $\theta_U \in \mathbb{R}^{p-b}$.
We consider several testing problems, each defined by its own null and alternative hypothesis, but all of the general form
\begin{align}\label{eq:generic_test}
    H_0: \theta \in \mathcal{H}_{0} \quad \textrm{vs.} \quad
    H_1: \theta \in \mathcal{H}_{1}\backslash \mathcal{H}_0
\end{align}
where $\mathcal{H}_{0}\subset \mathcal{H}_{1}$. 
Let $\mathcal{L}_{n}(\theta)$ denote the loss function
considered in  Section \ref{sec:estimator1} and
$-n\mathcal{L}_{n}(\theta)$ the corresponding 
 quasi log-likelihood.
The quasi-LR corresponding to \eqref{eq:generic_test} is
\begin{align}
\label{eq:LR}  
LR = 2n\left[\min_{\theta^{} \in\mathcal{H}_{0}}
  \mathcal{L}_{n}(\theta^{}) - \min_{\theta \in \mathcal{H}_{1}}\mathcal{L}_{n}(\theta)\right].
\end{align}
Let $\theta_0$ denote the population parameter and $\theta^{*}$ denote
the best fitting parameter under the null hypothesis. As the aim is to
derive the asymptotic distribution of $LR$ under the null hypothesis, 
we set $\theta^{*}=\theta_{0}$ and use the quadratic loss $Q_n$ defined in 
(\ref{eq:Qgeneral}) 
to obtain the asymptotic distribution of $LR$. Let 
$\mathcal{H}_{j}-\theta_0$ denote the parameter space 
$\mathcal{H}_{j}$ shifted by $\theta_0$.

It can be shown that if the null hypothesis holds with population parameter $\theta_0$ and
  condition (\ref{eq:rate2}) holds under both the null and alternative, then
\begin{align}
\label{eq:LRclassicalquad}
LR =& 2n\inf_{\delta\in \mathcal{H}_{0}-\theta_0}     
     Q_n(\delta) - 2n\inf_{\delta\in \mathcal{H}_{1}-\theta_0}Q_{n}(\delta) + o_p(1) \\
=&n\inf_{\delta\in \mathcal{H}_{0}-\theta_0} 
(\delta-\widetilde{\delta})^\top W(\delta-\widetilde{\delta}) - 
n\inf_{\delta\in \mathcal{H}_{1}-\theta_0} 
(\delta-\widetilde{\delta})^\top W(\delta-\widetilde{\delta}) + o_{p}(1),\nonumber
\end{align}
where we recall $\sqrt{n}\widetilde{\delta} \Dcon N(0,\Sigma)$ (under Assumption \ref{assum:CLT}).
\cite{p:sel-lia-87,p:and-01, p:hua-20, p:gui-mc-25} to name but a few authors, use the latter approximation as a means of deriving the limiting distribution of $LR$ under the null hypothesis. The focus has mainly been on
$W=\Sigma^{-1}$, where the distribution of $LR$ has been derived for certain hypothesis. The distribution in the 
general case $W\neq \Sigma^{-1}$ and for more general hypothesis has received less attention. 
We use the results in 
Sections \ref{sec:explicitexpression} and  \ref{sec:distribution}
to obtain an expression for 
(\ref{eq:LRclassicalquad}). This allows us to easily
deduce the asymptotic distribution of $LR$ under the null for arbitrary $W$ and general hypothesis.  
It has the additional benefit that the results 
when $\Sigma = W^{-1}$ are simply  a corollary of our general framework.


\subsection{Asymptotic distribution of the quasi-LR under 
different nulls}

We now derive the distribution of 
$LR$ under various different null. As a comparison we recall that in the case that $H_0:\theta_{B}=0$ vs $H_{1}:\theta_{B}\neq 0$ 
(a two sided test, hence no parameter constraints), 
then by Wilks' theorem
\begin{align}
\label{eq:classicalLR}
LR = n\delta^{\top}_{B}W_{B|U}\delta_{B}+o_{p}(1).
\end{align}
Asymptotically, this is a 
weighted chi-square with $|B|$-degrees of freedom. 

\vspace{1.5mm}
\noindent {\bf Hypothesis 1}
We first consider the one-sided hypothesis
\begin{align*}
    H_0: \theta_B = 0 \quad \textrm{vs.} \quad
    H_1: \theta_B \geq 0 \textrm{ and }\theta_{B}\neq 0.
\end{align*}
 We derive the asymptotic distribution of $LR$ under the null
hypothesis. An important feature of the null hypothesis is that
none of the parameters are restricted and can lie anywhere on
$\R^{p-b}$. It is this feature, together with Theorem \ref{theorem:CLT} that will give a relatively simple
distribution for $LR$.
In this case the quasi-LR in \eqref{eq:LR} is 
\begin{align}
\label{eq:LR1}
LR_{1} &= 2n\left[\inf_{\theta_{U} \in \R^{p-b}}\mathcal{L}_{n}(0_{B},\theta_{U})  -  \inf_{\theta \in \R_{+}^{b}\times\R^{p-b}}\mathcal{L}_{n}(\theta)\right].
\end{align}
In the theorem below we obtain the distribution of 
$LR_1$ under the null hypothesis. 
Using the 
quadratic approximation in Lemma \ref{lemma:QuadLLRT}, under the null hypothesis we have  
\begin{align*}
LR_1 &  =   2n[Q_n(0_B,\delta^{U}) - Q_n(\widehat{\delta})] + o_p(1),
\end{align*}
where 
$\delta^{U}= -W_{UU}^{-1}s_{U}$ (see (\ref{eq:deltaA})) and
$\widehat{\delta} = \argmin_{\delta\in \mathcal{C}}Q_{n}(\delta)$, whose explicit form is given in Theorem \ref{theorem:blocks}.

\begin{theorem}\label{theorem:LLRexpansion}
Let $LR_{1}$ be defined as in (\ref{eq:LR1}) and suppose 
Assumption \ref{assum:CLT} holds.
Suppose the null is true with population parameter $\theta_{B}=0$ and
  condition (\ref{eq:rate2}) is satisfied for 
   $\widehat{\theta}_{j,L}$
  and 
  $\widehat{\theta}_{j,Q}$ (evaluated under the null and the alternative;
  $j=0$ and $1$).
Then,  we have
\begin{align*}
 LR_1 &=  n\sum_{C\in \mathcal{P}(B)\backslash \emptyset}
(\delta^{F_{C}}_{C})^{\top} W_{C|U}
                  \delta^{F_{C}}_{C} 
  I\left[\delta^{F_{C}}_{C} > 0\right]
        I\left[W_{A_{C}|F_{C}}\widetilde{\delta}_{A_{C}} \leq
      0\right]+ 0 \cdot I\left[W_{B|U}\widetilde{\delta}_{B} \leq
      0\right]+    o_{p}(1),
\end{align*}
where 
$\sqrt{n}\delta^{F_C}_C \mid (\delta^{F_C}_C > 0 , W_{A_C|F_C} \widetilde{\delta}_{A_C}\leq 0) \Dcon  tSN_{|A_C|,|C|}(0,0,V_{A_CC}V_{CC}^{-1},V_{CC},V_{A_C|C})$
with the matrix $V_C$ defined in Theorem~\ref{theorem:CLT} and
$W_{C|U} = W_{CC}-W_{UC}^{\top}W_{UU}^{-1}W_{UC}$.

Thus the distribution of $LR_{1}$ is a mixture of 
quadratic forms of truncated skewed normals,
weighted by the probabilities
$\operatorname{pr} (s\in R_{C})$. 
\end{theorem}

Theorem \ref{theorem:LLRexpansion} is analogous to (\ref{eq:classicalLR}), but $LR_1$ is decomposed in terms of 
disjoint mixtures of weighted chi-square distributions, whose probabilities are determined by the probability of the event
$\operatorname{pr}(s\in R_{C})$. 
The reason for the appearance of the sum in Theorem \ref{theorem:LLRexpansion}, is that under the null
$R_{C} = \R^{p}$, whereas under the alternative 
$\cup_{C\in \mathcal{P}(B)}R_{C}=\R^{p}$ (hence the alternative is a refinement of the partition under the null). Interestingly, the distribution of the nuisance parameters do not play a role, except through its variance. 

The distribution of 
$LR_1$ is not pivotal as it depends on $\Sigma$ and $W$, however,  given estimates of $\Sigma$ and $W$ it can be simulated.
In Section \ref{sec:simtest} we describe how one can 
approximate the critical values of the distribution of $LR_1$ under the null hypothesis. 

We now consider the special case that $W=\Sigma^{-1}$, and show the distribution is pivotal. 

\begin{corollary}[The case $\Sigma = W^{-1}$]\label{cor:pivotal}
Suppose the conditions in Theorem~\ref{theorem:LLRexpansion} hold and $W= \Sigma^{-1}$. Then $W_{C|U} = \Sigma_{C|A_C}^{-1}$ and 
$\sqrt{n}\delta^{F_C}_C\mid(\delta^{F_C}_C > 0, \Sigma_{A_C A_C}^{-1}\widetilde{\delta}_{A_C} \leq 0) 
\Dcon  tN_{|C|}(0,\Sigma_{C|A_C})$. 
Thus 
\begin{align*}
    LR_1 \Dcon  
    \sum_{C \in \mathcal{P}(B)\backslash \emptyset} 
   \chi^{2}_{|C|}
    \operatorname{pr}[\delta^{F_C}_C > 0 ] 
    \operatorname{pr}[ \Sigma_{A_C A_C}^{-1}\widetilde{\delta}_{A_C} \leq 0 ]
    + 0 \cdot \operatorname{pr}[\Sigma_{BB}^{-1} \widetilde{\delta}_B  \leq 0  ].
\end{align*}
The above is sometimes called the chi-bar squared distribution.  Note that in the case $\Sigma = \sigma^{2}I$, then the weights $\operatorname{pr}[\delta^{F_C}_C > 0 ] 
    \operatorname{pr}[ \Sigma_{A_C A_C}^{-1}\widetilde{\delta}_{A_C} \leq 0 ]=(1/2)^{b}$. 
\end{corollary}
\cite{p:susko-13}, Theorem 1, directly follows from this result, specifically 
$(LR_{1} \mid |\widehat{C}| = c) \Dcon \chi^{2}_{c}$.
We note the reason that 
this test is pivotal (with the exception of the weights) is 
that the weight matrices in the partition $\Sigma_{C|A_{C}}^{-1}$ coincide with the precision matrix of the truncated normal. 

\vspace{1.5mm}
\noindent {\bf Hypothesis 2}
We now consider the mixed hypothesis, where 
some of the parameters we are testing do not lie on the boundary 
\begin{align*}
    H_0: \theta_B = 0 
    \textrm{ and }\theta_{U_{1}}=0
    \quad \textrm{vs.} \quad
    H_1: \theta_B \geq 0 \textrm{ and }\theta_{B}\neq 0
    \textrm{ or }\theta_{U_{1}}\neq 0,
\end{align*}
where $U_{1}\subseteq U$. A similar set-up is considered in
\cite{p:susko-13} and \cite{p:hua-20}.

We introduce the additional notation $U_{2} = U\backslash U_{1}$. Again, using that under the 
null hypothesis the parameter space is not restricted and
$LR$ is defined as  
\begin{align}
\label{eq:LR2}
LR_{2} &= 2n\left[\inf_{\theta \in \R^{|U_{2}|}}\mathcal{L}_{n}(0_{B},0_{U_{1}},\theta)  -  \inf_{\theta \in \R_{+}^{b}\times\R^{p-b}}\mathcal{L}_{n}(\theta)\right].  
\end{align}
Below we obtain the sampling distribution of $LR_2$.
\begin{lemma}\label{cor:LLRexpansion}
Let $LR_{2}$ be defined as in (\ref{eq:LR2}).
Suppose the assumptions in Theorem \ref{theorem:LLRexpansion} hold but now under 
the null hypothesis 
$H_0: \theta_B = 0 
    \textrm{ and }\theta_{U_{1}}=0$
    ($U_{1}\neq \emptyset$). Then 
we have
\begin{align*}
 LR_2 &=  n\sum_{\stackrel{C\in \mathcal{P}(B)}{C_{1}=C\cup U_1}}
(\delta^{F_{C}}_{C_1})^{\top} W_{C_1|U_2}
                  \delta^{F_{C}}_{C_1} 
  I\left[\delta^{F_{C}}_{C} > 0\right]
        I\left[W_{A_{C}|F_{C}}\widetilde{\delta}_{A_{C}} \leq
      0\right]+ o_{p}(1).
\end{align*}
In the 
case that $\Sigma=W^{-1}$, then
\begin{align*}
    LR_2 \Dcon  
   \sum_{\stackrel{C\in \mathcal{P}(B)}{C_{1}=C\cup U_1}}
\chi^{2}_{|C_{1}|} 
  \operatorname{pr}\left[\delta^{F_{C}}_{C} > 0\right]
    \operatorname{pr}\left[\Sigma_{A_{C} A_{C}}^{-1} \widetilde{\delta}_{A_{C}} \leq
      0\right]
\end{align*}
\end{lemma}

In hypothesis 1 and 2, the null is fixed and the 
nuisance parameters do not play a role. 
In the following two hypothesis we consider below 
either the null is not a point but a set or the 
nuisance parameters play a role in the 
limiting distribution. In such cases, it is not
possible to derive the distribution of 
$LR$ under the null. One needs to find critical values 
that are robust to the type of null or the 
nuisance parameters. For this we use the following 
simple fact that $Q_{n}(\delta_{\xi})$ is monotonically decreasing in $\xi$ where $\delta_{\xi} = \arg\min_{\delta\in \mathcal{C}_{\xi}}Q_{n}(\delta)$
(precisely, $Q_{n}(\delta_{\xi_1})\leq Q_{n}(\delta_{\xi_2})$ if 
$\xi_{1}<\xi_{2}$) (see Supplement \ref{sup:LLRT}). 

\vspace{1.5mm}
\noindent {\bf Hypothesis 3}
We consider the boundary 
constrained testing problems of type B, considered in \cite{b:sil-sen-05} and 
\cite{p:moh-20}. We focus on the special case 
\begin{align*}
    H_0: \theta_B \geq  0 
    \quad \textrm{vs.} \quad
    H_1: \theta_B \not\geq 0.
\end{align*}
Note that \cite{p:moh-20} consider the more general case that $H_{0}:A\mu\geq 0$ where the 
$\mu$ is the mean of a Gaussian random vector. 
Motivated by Section 3.8, \cite{b:sil-sen-05}, we define the quasi-likelihood  as
\begin{align}
\label{eq:LR3}
LR_3 = 2n\left[ \min_{\theta\in \mathcal{C}}\mathcal{L}_{n}(\theta) - 
\min_{\theta\in \R^{p}}\mathcal{L}_{n}(\theta)
\right]
\end{align}
where $\mathcal{C} = \R_{+}^{b}\times \R^{p-b}$. 
Below we obtain the asymptotic distribution of 
$LR_3$ under the null hypothesis and show how it connects to the least favourable null distribution defined in \cite{b:sil-sen-05} (Section 3.8) and the 
adaptive test described in \cite{p:moh-20}. 
Both \cite{b:sil-sen-05} and \cite{p:moh-20} consider the specific case of the Gaussian likelihood with boundary constraints on the mean and $W=\Sigma^{-1}$. 

\begin{theorem}\label{theorem:typeB}
Let $LR_{3}$ be defined as in (\ref{eq:LR3}).
Suppose the assumptions in Theorem \ref{theorem:LLRexpansion} hold but under the  null hypothesis 
$H_0: \theta_B \geq 0$. Then, for a given 
$\theta_B \in \R_{+}^{b}$
we have 
\begin{align}
 (LR_3\mid \theta_B) &= 
2n\min_{\delta\in\mathcal{C}_{\xi}}Q_{n}(\delta)+n\widetilde{\delta}^{\top}W\widetilde{\delta}+
 o_{p}(1)  \nonumber\\
 &=n \sum_{C\in \mathcal{P}(B)\backslash B}
 (\xi_{A_C} - \tilde{\delta}_{A_C})^\top W_{A_C|F_C} 
(\xi_{A_C} - \tilde{\delta}_{A_C}) 
  I\left[\delta^{F_{C}}_{C} > \xi_{C}\right]
        I\left[W_{A_{C}|F_{C}}(\widetilde{\delta}_{A_{C}}-\xi_{A_{C}}) \leq
      0\right] \nonumber\\
      &\qquad +0\cdot I\left[\delta^{F_{B}}_{B} > \xi_{B}\right]+o_{p}(1), \label{eq:normexpansion}
\end{align}
where $\xi_{B} = -\theta_{B}$, 
$\widetilde{\delta} = -W^{-1}s$.
Thus, asymptotically 
the least favourable distribution under the null is at
$\theta_B=0$ and for any  
$\theta_{B}\in \mathbb{R}_{+}^{b}$ we have 
$(LR_{3}\mid \theta_{B})\leq (LR_{3} \mid \theta_{B}=0)+o_{p}(1)$. 
If $W = \Sigma^{-1}$, then 
\begin{align}
(LR_{3}\mid \theta_B  = 0) \Dcon 
 \sum_{C\in \mathcal{P}(B)\backslash B}
\chi^{2}_{|A_{C}|}
  \operatorname{pr}\left[\delta^{F_{C}}_{C} > 0\right]
        \operatorname{pr}\left[\Sigma_{A_C A_C}^{-1} \widetilde{\delta}_{A_{C}} \leq
      0\right]. \label{eq:chitype2}
\end{align} 
\end{theorem}
We now connect Theorem \ref{theorem:typeB} to the results of 
\cite{b:sil-sen-05} and \cite{p:moh-20} who consider the case that 
$W=\Sigma^{-1}$.  The distribution of  
$LR_{3}$ is unknown since the null  
$\mathcal{C}$ is a set and not a specific value. However, our results mirror Theorem 3.8.1, \cite{b:sil-sen-05} who shows that the 
$p$-value over all values in the null space is maximised at the boundary $\theta_{B} = 0$. With this 
\cite{b:sil-sen-05} suggest using the mixture of chi-squares given in 
(\ref{eq:chitype2}) to obtain the least favourable p-values under null
(obtained using simulations
as an exact expression for the weights is not given). On the other hand, \cite{p:moh-20} 
suggest using a data adaptive test that they show has more power. 
Putting into the framework of Theorem \ref{theorem:typeB}, this corresponds to using  $\theta_B=0$ in the conditioning set
\begin{align*}
(LR_{3}\mid \theta_{B} = 0, |\widehat{C}|=c)
\Dcon \chi^{2}_{b-c}
\end{align*}
where we recall that 
$\widehat{C}=(j\in B\textrm{ where }\widehat{\theta}_{j,L}>0)$ 
and $\widehat{\theta} = \arg\min_{\theta\in \mathcal{C}}\mathcal{L}_{n}(\theta)$. 

Away from the case that $W=\Sigma^{-1}$, the 
form for the least favourable distribution 
$LR_{3} \mid \theta_{B}=0$ given in (\ref{eq:normexpansion}) 
has a curious asymmetry to 
Theorem \ref{theorem:LLRexpansion}
(and correspondingly 
(\ref{eq:chitype2}) and Corollary \ref{cor:pivotal}). The partition sets are the 
same, the difference is the quadratic form 
$ (\tilde{\delta}_{A_C})^\top W_{A_C|F_C} 
\tilde{\delta}_{A_C}$ (in Theorem \ref{theorem:typeB}) vs $(\delta^{F_{C}}_{C})^{\top} W_{C|U}
\delta^{F_{C}}_{C}$ (in Theorem \ref{theorem:LLRexpansion}).

In the quasi-likelihood case,
Theorem \ref{theorem:typeB} together with
a simulation based method (similar to that described in Section \ref{sec:simtest}) can be used to obtain the critical values of the
least favourable null distribution 
$LR_{3}\mid \theta_{B}=0$ and its adaptive counterpart $LR_{3}\mid (\theta_{B}=0,|\widehat{C}|=c)$. 
 
So far we have considered
cases where nuisance parameters play no role in the asymptotic distribution. Below we show that boundary constrained nuisance parameters do play a role. 

\vspace{1.5mm}
\noindent {\bf Hypothesis 4}
\emph{The case that some of the nuisance parameters are boundary constrained} The hypothesis is 
\begin{align*}
    H_0: \theta_{B_{1}} =  0 
    \quad \textrm{vs.} \quad
    H_1: \theta_{B_1} \geq 0 \textrm{ and } \theta_{B_1} \neq 0
\end{align*}
where $B_{1}\subset B$ (one can also include a subset of $\theta_{U}$ in the above too).
Let $B_{2} = B\backslash B_{1}$, such that 
$\theta_{B_{2}}$ are the boundary constrained nuisance parameters, $b_{1}=|B_{1}|$ and
$b_2=|B_{2}|$. Let $\xi_{B_2} = -\theta_{B_2} = 
(\xi_{2,1},\ldots,\xi_{2,b_2})$.
It can be shown that LR under the null is 
\begin{align*}
LR_4 \mid \theta_{B_2} = 2n\left[
\inf_{\delta\in \mathcal{C}_{0,\xi_{B_2}}}Q_{n}(0_{B_1},\delta) - \inf_{\delta\in \mathcal{C}_{1,\xi_{B_2}}}Q_{n}(\delta)
\right] + o_{p}(1),
\end{align*}
where
$\mathcal{C}_{0,\xi_{B_2}} = 
\prod_{j=1}^{b_2} [\xi_{2,j}, \infty)\otimes \R^{p-b}$ 
and 
$\mathcal{C}_{1,\xi_{B_2}} = 
\R_{+}^{b_1}\otimes \mathcal{C}_{0,\xi_{B_2}}$.
When there are nuisance parameters on the boundary, then the
partition sets of $\R^{p}$ under the alternative is not a refinement of that under the null.  This is in contrast to 
hypothesis 1-3 considered above. Consequently, the nuisance parameters play a role in the distribution of LR. In Supplement 
\ref{sup:LLRT} we obtain an exact expression for 
$LR_4$ which clearly illustrates this fact. 

This property is well known in the econometric literature (\cite{p:and-gug-09,p:mcc-17,p:gui-mc-25}). 
For example, the first paragraph of 
\cite{p:gui-mc-25} states ``As is now well-known, the asymptotic distribution of a likelihood ratio (LR) statistic is discontinuous under the null hypothesis being tested at points for which a vector nuisance
parameter lies on the boundary of its parameter space''. However,  it is not immediately obvious 
what is meant by this. The expansion of LR gives a clear meaning;  nuisance parameters that are close to the boundary impact the distribution under the null (which is also apparent in Section \ref{sec:closetoboundary}). 
There exist several interesting solutions; \cite{p:gui-mc-25} is one such example. Motivated by \cite{p:mcc-17}, \cite{p:gui-mc-25} 
exploit the monotonicity of 
$\inf_{\delta\in \mathcal{C}_{\epsilon}}Q(\delta)$ to replace 
the $\theta_{B_2}$ in 
$\mathcal{C}_{j,\xi_{B_2}}$ with the bounds of the rectangular confidence intervals based on the unconstrained estimator $\widetilde{\theta}_{B}$. 
They show that the critical values are computationally fast to implement, uniformly holds the level and has more power than previous methods. 



\section{Simulation based approach to Critical Values}\label{sec:simtest}

In this section, we apply the results from the previous sections to 
obtain the distribution of 
$LR$ under certain null hypotheses.
We focus on Hypothesis 1 considered in 
Section \ref{sec:LLRT}.
The results in Section \ref{sec:distribution}
demonstrate that the distribution of 
$\widehat{\theta}_{L}-\theta_0$ under boundary constraints 
does not have a classical textbook distribution that standard software can provide. 
Nevertheless, 
the results readily lend themselves to a simulation based approach for obtaining the critical values of the test statistic.   
In this section, we consider two testing procedures, one based on the marginal distribution of 
$LR$ under the null hypothesis and the other based on 
the conditional likelihood $LR\mid \widehat{C}=C$  (the latter method is motivated by \cite{p:susko-13}). 

In the first procedure, 
we simulate the distribution of $LR$ under the null $H_{0}:\theta_{B} = 0$. 
Theorem \ref{theorem:blocks} forms the foundation of this procedure. 
The pseudo code is given in the Algorithm 1.

Given the simulated conditional $LR^{(1)}, \dots, LR^{(M)}$,
we approximate the critical value by the empirical $(1-\alpha)$-quantile $\widehat{q}_{1-\alpha}$.
The null hypothesis is rejected when $LR_{obs} > \widehat{q}_{1-\alpha}$.

\begin{algorithm}[h!]
\caption{Simulation of the marginal distribution of $LR$ under $H_{0}:\theta_{B}=0$}
\label{alg:full_simulation}



\KwIn{$\Sigma$, $W$, sample size $n$, index sets $B$ and $U$,
and number of iterations $M$.}

\BlankLine
\ForEach{$C\in\mathcal{P}(B)$}{
    Set $A_C=B\setminus C$ and $F_C=C\cup U$\;
    Compute $W_{A_C|F_C}$ and, if $C\neq\emptyset$, $W_{C|U}$\;
}

\BlankLine
\For{$m\leftarrow 1$ \KwTo $M$}{
    Draw $\widetilde{\delta}^{(m)}
    \sim N(0,\Sigma/n)$\;
    Set $LR^{(m)}=0$\;

    \ForEach{$C\in\mathcal{P}(B)$}{
        Compute $\delta^{F_C}
        =\widetilde{\delta}^{(m)}_{F_C}
        +W_{F_CF_C}^{-1}W_{F_CA_C}
        \widetilde{\delta}^{(m)}_{A_C}$\;

        \If{$\delta_C^{F_C}>0$ \textbf{and}
        $W_{A_C|F_C}\widetilde{\delta}^{(m)}_{A_C}\leq0$}{
            \If{$C\neq\emptyset$}{
                Set $LR^{(m)}
                =n(\delta_C^{F_C})^\top W_{C|U}\delta_C^{F_C}$\;
            }
            \textbf{break}\;
        }
    }
}

\KwOut{$\{LR^{(1)},\ldots,LR^{(M)}\}$}
\end{algorithm}


\vspace{1.5mm}

\noindent{\bf Approximating the conditional $LR$}

\vspace{1.5mm}

\noindent Approximating the empirical 
critical values of 
$LR$ can be computationally intensive
if $B$ is large, as it requires moving through all the sets in the power set $\mathcal{P}(B)$. 
An alternative and computationally faster method is to obtain the conditional likelihood of $LR$ conditioned on the set $C$, where 
$\widehat{\theta}_{A_C,L}=0$ and 
$\widehat{\theta}_{F_C,L}$
lies in the interior of $\R_{+}^{|C|} \otimes \R^{p-b}$.
This approach is similar in spirit to the conditional testing procedure in \cite{p:susko-13}.
Let $\widehat{C} =\{j \in B: \widehat{\theta}_{j,L} > 0\}$.
Then $\widehat{C}$ identifies the boundary restricted parameters that are estimated to lie in the interior of the parameter space.
For a fixed $C \in \mathcal{P}(B)$, by applying Theorem~\ref{theorem:LLRexpansion} 
under the null hypothesis 
$H_{0}$ and $\xi_{B} = 0$, we have
\begin{align*}
    LR \mid (\widehat{C}=C)
    =
    \left\{
    \begin{array}{cc}
    n (\delta^{F_C}_{C})^\top W_{C|U} \delta^{F_C}_{C} \mid (W_{A_C|F_C} \widetilde{\delta}_{A_C}\leq 0, \delta^{F_C}_C>0)+ o_{p}(1) & C\neq \emptyset \\
    O_{p}(n^{-1/2})& C = \emptyset \\
    \end{array}
    \right.
\end{align*}
From Theorem~\ref{theorem:CLT}, 
$\delta^{F_C}_C \mid (W_{A_C|F_C} \widetilde{\delta}_{A_C} \leq 0, \delta^{F_C}_C>0)$ has a truncated skew-normal distribution, 
whose distribution can be simulated without searching over all sets in $\mathcal{P}(B)$.
The pseudo code for this procedure 
can be found in Supplement 
\ref{sup:hyp}.

\begin{remark}\label{remark:LRinvSigmaW}
In the special case $\Sigma = W^{-1}$ we have
    $LR  \mid (\widehat{C}=C) \Dcon \chi^2_{|C|}$,
where the above follows from 
Corollary \ref{cor:pivotal}.
Moreover, the limiting conditional distribution depends on $C$ only through its cardinality.
Thus, one can condition on 
all the sets where $|C|=c$ to give 
$LR  \mid (|\widehat{C}|=c) \Dcon \chi^2_{c}$.
The latter procedure corresponds to 
the conditional test proposed in 
\cite{p:susko-13}. 
However, in the case that 
$W\neq \Sigma^{-1}$ this does not hold.
\end{remark}

Given the simulated conditional $LR^{(1)}, \dots, LR^{(M)}$,
we approximate the critical value by the empirical $(1-\alpha)$-quantile $\widehat{q}_{1-\alpha, C}$.
For the observed set $\widehat{C}$, the null hypothesis is rejected when $LR_{obs} > \widehat{q}_{1-\alpha, \widehat{C}}$.


\section{Simulations}\label{sec:sims}

Motivated by \cite{p:bes-74,p:yan-15} who model multivariate count data in terms of univariate conditional Poissons, we test for dependence between the variables (conditional interactions). 
Define $(X,Y,Z)$ to follow the conditional trivariate Poisson distribution if all its univariate conditionals are Poisson distributed:
\begin{align*}
X\mid Y,Z\sim Pois(\lambda_{x}(Y,Z)), \quad
Y\mid X,Z\sim Pois(\lambda_{y}(X,Z)), 
\quad 
Z\mid X,Y\sim Pois(\lambda_{z}(X,Y)),
\end{align*}
where the conditional intensities are modeled as
\begin{align*}
\log \lambda_{x}(Y,Z) &= 
\psi_{x} - \phi_{xy}Y - \phi_{xz}Z,\quad
\log \lambda_{y}(X,Z) = 
\psi_{y} - \phi_{xy}X - \phi_{yz}Z, \\
\log \lambda_{z}(X,Y) &= 
\psi_{z} - \phi_{xz}X - \phi_{yz}Y.
\end{align*}
In order that a joint distribution exists, certain
{\it compatibility conditions} \citep{p:bes-74} need to hold. Specifically, 
$(\psi_{x},\psi_{y},\psi_{z})\in \R^{3}$ and 
$(\phi_{xy},\phi_{yz},\phi_{xz})\in \R_{+}^{3}$. In the case that $\phi_{xy}=0$, 
$X$ and $Y$ are conditionally independent given $Z$ (similar for the other terms). The case 
$\phi_{xy}=0,\phi_{yz}=0,\phi_{xz}=0$ corresponds to marginal independence between the variables. We  conduct a test for marginal independence, i.e., 
$H_{0}:(\phi_{xy},\phi_{yz},\phi_{xz}) =(0, 0, 0)$ vs
$H_{1}:(\phi_{xy},\phi_{yz},\phi_{xz}) \geq 0$ (at least one parameter is non-zero). 

As the joint probability mass function of $(X,Y,Z)$ lacks an analytic form due to the presence of an intractable nonrandomising constant,  we use the negative {\it partial likelihood} based on the conditional distributions to estimate the parameters. Specifically, define
\begin{align*}
\mathcal{L}_{n}(\theta) = 
\frac{1}{n}\sum_{i=1}^{n}\ell_{i}(\theta)
\end{align*}
where $\theta = (\psi_x, \psi_y, \psi_z, \phi_{xy}, \phi_{yz}, \phi_{xz})^\top$,
$\ell_{i}(\theta) = -
(\ell_{x,i}(\theta)+\ell_{y,i}(\theta)+\ell_{z,i}(\theta))$ and
\begin{align*}
\ell_{x,i}(\theta) &= X_{i}(\psi_{x} - \phi_{xy}Y_{i} - \phi_{xz}Z_{i})
-\exp(\psi_{x} - \phi_{xy}Y_{i} - \phi_{xz}Z_{i}),\\
\ell_{y,i}(\theta) &= Y_{i}(\psi_{y} - \phi_{xy}X_{i} -\phi_{yz}Z_{i})
-\exp(\psi_{y} - \phi_{xy}X_{i} - \phi_{yz}Z_{i}),\\
\ell_{z,i}(\theta) &= Z_{i}(\psi_{z} - \phi_{xz}X_{i} - \phi_{yz}Y_{i})
-\exp(\psi_{z} - \phi_{xz}X_{i} - \phi_{yz}Y_{i}).
\end{align*}
Thus $\mathcal{C}=\R^{3}\times \R_{+}^3$. Using the above we evaluate 
$LR = 2n[\mathcal{L}_{n}(\widehat{\theta}^{*}, 0) - \mathcal{L}_{n}(\widehat{\theta})]$, where 
$\widehat{\theta}^{*} = (\widehat{\psi}_{x},\widehat{\psi}_{y},\widehat{\psi}_{y})
=\argmin_{\R^{3}} \mathcal{L}_{n}(\psi_{x},\psi_{y},\psi_{z}, 0)$ and 
$\widehat{\theta} = \argmin_{\mathcal{C}} \mathcal{L}_{n}(\theta)$. 
To obtain 
the critical values using Algorithms 1 and 2 we require 
$W$ and $\Sigma$ which we estimate with
\begin{align*}
\widehat{W} &=\frac{1}{n}\sum_{i=1}^{n}\nabla^{2}\ell_{i}(\widetilde{\theta}), \quad \widehat{\Sigma} = \widehat{W}^{-1}\widehat{V}\widehat{W}^{-1}
\textrm{ where } 
\widehat{V} =\frac{1}{n}\sum_{i=1}^{n}(\nabla\ell_{i}(\widetilde{\theta})) (\nabla\ell_{i}(\widetilde{\theta}))^{\top},
\end{align*}
where $\widetilde{\theta}=\argmin_{\R^6} \mathcal{L}_{n}(\theta)$ is an unconstrained minimiser over $\R^6$.

In the simulations we set 
$(\psi_x, \psi_y, \psi_z) = (1,1,1)$. 
Under the null 
$(\phi_{xy},\phi_{yz}, \phi_{xz})$ 
is set to $(0,0,0)$  and under the alternative 
it is set to $(0.1,0.1,0.1)$.
We evaluate the empirical Type I error and power of the global and conditional 
LR test evaluated over $R=1000$ replications
for sample size $n \in \{10, 50, 100, 200\}$ and level $\alpha \in \{0.01, 0.05, 0.1\}$. The critical values are obtained by setting $M=1000$. The results are given in Tables \ref{tab:global_lr_A1} and \ref{tab:conditional_lr_A2}. The performance of both tests appear to be quite similar -- though the simulations suggest 
the conditional test may have slightly less power. 

\begin{table}[htbp]
\centering
\caption{Empirical Type I error and power of the global LR test from Algorithm~\ref{alg:full_simulation}.}
\label{tab:global_lr_A1}

\begin{tabular}{r rr rr rr}
\toprule
& \multicolumn{2}{c}{$\alpha = 0.01$}
& \multicolumn{2}{c}{$\alpha = 0.05$}
& \multicolumn{2}{c}{$\alpha = 0.10$} \\
\cmidrule(lr){2-3}
\cmidrule(lr){4-5}
\cmidrule(lr){6-7}
$n$
& Type I error & Power
& Type I error & Power
& Type I error & Power \\
\midrule
 10 & 0.002 & 0.016 & 0.051 & 0.179 & 0.133 & 0.368 \\
 50 & 0.010 & 0.305 & 0.066 & 0.575 & 0.125 & 0.716 \\
100 & 0.011 & 0.631 & 0.066 & 0.846 & 0.117 & 0.906 \\
200 & 0.012 & 0.947 & 0.061 & 0.991 & 0.107 & 0.998 \\
\bottomrule
\end{tabular}
\end{table}
\begin{table}[htbp!]
\centering
\caption{Empirical Type I error and power of the conditional LR test from Algorithm~\ref{alg:conditional_LR}.}
\label{tab:conditional_lr_A2}

\begin{tabular}{r rr rr rr}
\toprule
& \multicolumn{2}{c}{$\alpha = 0.01$}
& \multicolumn{2}{c}{$\alpha = 0.05$}
& \multicolumn{2}{c}{$\alpha = 0.10$} \\
\cmidrule(lr){2-3}
\cmidrule(lr){4-5}
\cmidrule(lr){6-7}
$n$
& Type I error & Power
& Type I error & Power
& Type I error & Power \\
\midrule
 10 & 0.001 & 0.006 & 0.031 & 0.089 & 0.091 & 0.214 \\
 50 & 0.011 & 0.186 & 0.057 & 0.425 & 0.115 & 0.536 \\
100 & 0.014 & 0.486 & 0.047 & 0.740 & 0.091 & 0.830 \\
200 & 0.010 & 0.884 & 0.046 & 0.971 & 0.099 & 0.988 \\
\bottomrule
\end{tabular}
\end{table}



\subsection*{Acknowledgements}
The second and third
authors were supported in part by NSF Grants DMS-2610612 and DMS-2210726 respectively.

\subsection*{Supplementary Material}\label{supplementary_material}
The Supplementary Material includes a basic summary of 
the 
results from optimisation that we use, 
the proofs of 
all the results in this paper and 
several examples of the boundary constrained estimator for the case $p=1$, 
$p=2$ and $p=3$.
The R code for reproducing the simulations and examples is available at 
\url{https://github.com/jaeseonlee-stat/Boundary-Inference}.

\bibliographystyle{abbrvnat}
\bibliography{literature}

\newpage

\appendix

\section{Background results from optimisation}
\label{sec:optim_results}
Many of the results in this paper use arguments from optimisation theory, in
particular regarding directional derivatives. We formally define the
directional derivative below and describe how it relates to the minimum of a loss
function over a set (a comprehensive review is given in \cite{b:bor-lew-00}).

\begin{definition}[Feasible direction]\label{def:feasible}
Let $C \subset \R^p$ be a closed nonempty set.
A vector $\nu \in \R^p$ is said to be a feasible direction at $x \in C$ if
there exists a scalar 
$\tau > 0$ such that  $x + t \nu \in C$ for all $0 \leq t \leq \tau$.
\end{definition}

\begin{definition}[Directional Derivative]
Let $f: C \to \R$ be a smooth function, where $C$ is as in Definition \ref{def:feasible}. The directional derivative of the function $f$ at $x\in C$ in the
feasible direction $\nu$ is defined as $D_\nu f(x) = \lim_{t \to 0} \frac{f(x + t \nu) - f(x)}{t}$, where the right hand side is well-defined since $\nu$ is a feasible direction. If $\nabla f(x)$ denotes the vector of partial derivatives of $f$ at $x$, then one has $D_{\nu}f(x)=\nu^{\top}\nabla f(x)$.
\end{definition}  

\begin{lemma}[Directional Derivative at the minimum point]\label{lemma:direction}
   Let $L$ be a smooth loss function defined on a convex set $C$, and $\theta_0 = \argmin_{\theta\in    C} L(\theta)$.
   Then $D_{\nu}L(\theta_0)\geq 0$ for all feasible directions $\nu$ at
  $\theta_0$ in $C$.
  Further, if $\theta_0$ lies within the interior of $C$ then 
  $\nabla L(\theta_0) = 0$. 

  Conversely, suppose $L$ is a convex function. If there exists
a $\theta_0\in C$ where for all  $\theta\in
 C$ we have 
    $(\theta-\theta_0)^\top \nabla L(\theta_0) \geq 0$, then  
    $\theta_0 = \argmin_{\theta \in C}L(\theta)$. 
\end{lemma}
\begin{proof}
This result is well known, but for 
 completeness we give a proof. Suppose $h \in (0,1)$ and
  $\theta\in C$, then by convexity of $C$ we have 
 $\theta_0 + h (\theta-\theta_0) \in C$.
Since $\theta_0$ is a minimiser of $L$ over $C$, then
\begin{align}
\label{eq:Ltheta0}
    L(\theta_0) \leq L\big( \theta_0 + h(\theta-\theta_0) \big).
\end{align}
Set $\nu = \theta-\theta_0$ and fixing $\nu$ 
we define the function $g_{\nu}(h) = L\big( \theta_0 + h\nu
\big)$. Then (\ref{eq:Ltheta0}) can be written as     $g_{\nu}(0) \leq g_{\nu}(h)$.
Consequently, for all $h\in (0,1)$ we have
\begin{align*}
    \frac{g_{\nu}(h) - g_{\nu}(0)}{h} \geq 0.
\end{align*}
Thus 
\begin{align*}
  \lim_{h \downarrow 0} \frac{g_{\nu}(h) - g_{\nu}(0)}{h} 
    =\nu^{\top}\nabla L(\theta_0
)=D_{\nu}(\theta_0)
    \geq 0.
\end{align*}
As the above holds for all
feasible directions $\nu$ at $\theta_0$ in $C$, then
$D_{\nu}(\theta_0)\geq 0$. Thus we have proved  the first claim.


Now suppose $\theta_0$ lies within the interior of $C$. 
Then both $\nu$ and $-\nu$ are feasible directions at $\theta_0$, so
\begin{align*}
    D_{\nu}(\theta_0) = \nu^\top \nabla L(\theta_0) \geq 0 \quad \textrm{and} \quad
D_{-\nu}(\theta_0) =  -\nu^\top \nabla L(\theta_0) \geq 0.
\end{align*}
Since $D_{\nu}(\theta_0)=-D_{-\nu}(\theta_0)$ for all feasible
directions at $\theta_0$, 
and both are greater than or equal to zero, then $\nabla L(\theta_0) = 0$.

We now prove the last claim. 
Suppose $L$ is convex, then for all $\theta,\theta_0\in C$ we have
\begin{align*}
    L(\theta) \geq L(\theta_0) + (\theta-\theta_0)^\top \nabla L(\theta_0).
\end{align*}
If there exists a $\theta_0\in C$ where $(\theta-\theta_0)^\top \nabla L(\theta_0)\geq 0$ for all  $\theta\in C$, the above
 display immediately implies $L(\theta)\geq L(\theta_0)$, thus proving that
$\theta_0 = \argmin_{\theta \in C}L(\theta)$.
\end{proof}

We now state an expanded version of Lemma \ref{lemma:boundary} below. 
\begin{lemma}
\label{lemma:boundary2}
Suppose $L$ is a strictly convex loss function and let $\widehat{\theta} = \arg\min_{\theta\in \mathcal{C}}L(\theta)$ and $\widetilde{\theta} = \arg\min_{\theta\in \R^{p}}L(\theta)$.
\begin{itemize}
\item[(i)] If there exists $1 \le j \le b$ such that $\widetilde{\theta}_{j} < 0$, then $\widehat{\theta}$ lies on the boundary of $\mathcal{C}$; that is, there exists at least one $1 \le k \le b$ with $\widehat{\theta}_{k} = 0$.
\item[(ii)] Suppose $\widetilde{\theta}_{1} < 0$ while $\widetilde{\theta}_{j} \ge 0$ for $j \in \{2,\ldots,b\}$. Then $\widehat{\theta}_{1} = 0$ and
\[
\widehat{\theta}_{-1} = \argmin_{\theta_{-1} \in \R^{b-1}_{+}\times \R^{p-b}} L(0,\theta_{-1}).
\]
\end{itemize}
\end{lemma}
\begin{proof}
 To  prove (i)
    we prove the contrapositive, that is if $ \widehat{\theta}_{j} > 0$ for all  $1
    \leq j
    \leq b$ then $\widetilde{\theta}_{j}>0$ for all $1\leq j\leq b$.
    Since
    $\widehat{\theta}$ is an interior point of the $\mathcal{C}$
    by Lemma \ref{lemma:direction}
          $\nabla L(\widehat{\theta}) = 0$.
     Strictly convex of $L$, means  
     the solution  $\nabla L(\theta) = 0$
     is unique, implying 
         $\widehat{\theta} = \widetilde{\theta}$.
     Thus $\widetilde{\theta}_{j}>0$ for 
     $1\leq j \leq b$. This proves (i).

   To prove (ii), suppose by contradiction, that 
        $\widehat{\theta}_{1} > 0$.
    Since we have $\widetilde{\theta}_{1} < 0$ and $\widehat{\theta}_{1} > 0$, there exists $\lambda \in (0,1)$ such that
    \begin{align*}
        \lambda \widehat{\theta}_{1} + (1-\lambda) \widetilde{\theta}_{1} = 0.
    \end{align*}
    With this $\lambda$ we define
    \begin{align*}
        \theta^{(\lambda)} = \lambda \widehat{\theta} + (1-\lambda) \widetilde{\theta},
    \end{align*}
    where by construction
    $\theta_{1}^{(\lambda)} = \lambda \widehat{\theta}_{1} + (1-\lambda)
    \widetilde{\theta}_{1} = 0$. 
    Under the stated conditions
$(0,\widetilde{\theta}_{2},\ldots,\widetilde{\theta}_{p})\in \mathcal{C}$, thus 
$\theta^{(\lambda)} = \lambda \widehat{\theta} + (1-\lambda) \widetilde{\theta}\in \mathcal{C}$. We use this fact to obtain the contradiction.

   Returning to the loss function $L$.  By strict convexity of $L$
    \begin{align*}
        L(\theta^{(\lambda)}) < \lambda L(\widehat{\theta}) + (1-\lambda) L(\widetilde{\theta}) \leq L(\widehat{\theta})
    \end{align*}
    where the last inequality is due to the fact that
    $L(\widetilde{\theta}) \leq L(\widehat{\theta})$ (by definition
    of $\widetilde{\theta}$).
   Since $\theta^{(\lambda)} \in \mathcal{C}$, this  
  contradicts that $\widehat{\theta}$ is the minimiser over
  $\mathcal{C}$. 
Consequently, $\widehat{\theta}_{1}=0$ and
$$\widehat{\theta}_{-1} = \argmin_{\theta_{-1} \in \R^{b-1}_{+}\times\R^{p-b}}
L(0,\theta_{-1}).$$
Thus proving the result. 
\end{proof}

\section{Auxillary lemmas}

The following lemma is used in the proof of 
Lemma \ref{lemma:LSrateG}, Lemma \ref{lemma:rateLS}, and Lemma \ref{lemma:finite_approx}. The lemmas in this section may also be of independent interest. 
\begin{lemma}\label{lemma:convergence}
Suppose that $\{\widetilde{\delta}_{n}\}_{n\geq 1}$ are a sequence of random vectors where  
$\widetilde{\delta}_{n}\Pcon 0$ as $n\rightarrow \infty$. If 
$W$ is a positive definite matrix and $\alpha$
 is a vector where 
 $\alpha\geq 0$ and $\alpha\neq 0$, then 
\begin{align*}
\operatorname{pr} \left( W(\widetilde{\delta}_{n}+\alpha)\leq 0\right)\rightarrow 0 \quad
\textrm{ as }
n\rightarrow \infty.
\end{align*}
\end{lemma}
\begin{proof}
 To prove this, we use the following fact. Suppose $W$ 
 is a $p\times p$
 positive definite matrix and $\alpha$ is a $p$-dimensional vector where $\alpha\geq 0$ but $\alpha\neq 0$, then 
for at least one entry $j$ we have $(W\alpha)_{j}>0$. 
For completeness we now prove this fact. 
By positive definiteness of $W$ then for all $\alpha\neq 0$ we have $\alpha^{\top}W\alpha=\sum_{i=1}^{p}\alpha_{i}(W\alpha)_{i}>0$. If $\alpha\geq 0$, then there must exist at least one $j$ where $(W\alpha)_{j}>0$, else $\alpha^{\top}W\alpha\leq 0$, which would 
contradict the assumption that
$W$ is positive definite or $\alpha \neq 0$. We now show that the result follows from this fact. 

By the above there exists a $j$ where $(W\alpha)_{j}>0$. Hence, since 
$\widetilde{\delta}_{n}\rightarrow \infty$
we have
\begin{align*}
(W(\widetilde{\delta}_{n}+\alpha))_{j}=
(W\widetilde{\delta}_{n})_{j}+
(W\alpha)_{j}\Pcon (W\alpha)_{j}>0
\end{align*}   
This implies $\operatorname{pr} \left( W(\widetilde{\delta}_{n}+\alpha)\leq 0\right)\rightarrow 0$ as $n\rightarrow \infty$.
\end{proof}

The following lemma is used in the proof of 
Corollary \ref{cor:pivotal}.
\begin{lemma} \label{lemma:trun_chi}
    Let $Z=(Z_1, \dots, Z_q)^\top \sim N_q(0, I_q)$.
    Then
    \begin{align*}
        Z^\top Z \mid (Z>0) \sim \chi^2_q.
    \end{align*}
    Further, let $A$ be a $q \times q$ positive definite matrix. Then
    \begin{align*}
        Z^\top Z \mid (A Z>0) \sim \chi^2_q.
    \end{align*}
\end{lemma}
\begin{proof}
    Since $Z_i$ follows a univariate standard normal distribution,
    \begin{align*}
        Z_i^2 \mid (Z_i>0) \sim \chi^2_1, 
        \quad 
        i=1,\dots, q.
    \end{align*}
    Since $Z_1, \dots, Z_q$ are independent, conditioning on $(Z>0)$ preserves the independence.
    Thus
    \begin{align*}
        Z^\top Z \mid (Z>0) \sim \chi^{2}_{q}.
    \end{align*}

    To prove the second claim, we use the following two facts
    \begin{align}
    \label{eq:fact1}
    (AZ > 0) = (AZ /\|Z\|_2> 0)
    \end{align}
    almost surely since $\|Z\|_2>0$ almost surely. The second fact is that if 
    $Z\sim N_{q}(0,I_{q})$ then 
    \begin{align}
    \label{eq:fact2}
    Z^{\top}Z=\|Z\|_2^2 \textrm{ and } 
     Z/\|Z\|_2 
    \end{align}
    are independent 
     (see \cite{b:bil-99}, Corollary~4.3). Using this we have 
    \begin{align*}
        Z^\top Z \mid (A Z>0) 
        \stackrel{\mathcal{D}}{=} 
         \|Z\|_{2}^{2}
        \mid \left( A \frac{Z}{\|Z\|_2}> 0 \right)\stackrel{\mathcal{D}}{=}
         Z^\top Z \sim \chi^2_{q},
    \end{align*}
    the second equality in the display follows from (\ref{eq:fact1}) and third from 
    (\ref{eq:fact2}).
\end{proof}

The following lemma is used in the proof of 
Lemma \ref{cor:LLRexpansion}.
\begin{lemma} \label{lemma:trun_chi2}
$X_{1}$ and $X_{2}$ are $q_{1}$ and 
    $q_{2}$-dimensional random vectors respectively.
    Suppose $X=(X_1^\top, X_2^\top)^\top \sim N_{q_1+q_2}(0, \Sigma)$.
    Then
    \begin{align*}
        X^\top \Sigma^{-1} X \mid (X_1>0) \sim \chi^{2}_{q_1+q_2}.
    \end{align*}
\end{lemma}
\begin{proof}
    Using the Schur complement,
    \begin{align*}
        \Sigma^{-1}
        &= \left( \begin{matrix}
            \Sigma_{11}^{-1} + \Sigma_{11}^{-1}\Sigma_{12}\Sigma_{2|1}^{-1}\Sigma_{21}\Sigma_{11}^{-1} & -\Sigma_{11}^{-1}\Sigma_{12}\Sigma_{2|1}^{-1} \\
            -\Sigma_{2|1}^{-1}\Sigma_{21}\Sigma_{11}^{-1} & \Sigma_{2|1}^{-1}
        \end{matrix} \right)  \\
        &= \left( \begin{matrix}
            I_{q_1} & -\Sigma_{11}^{-1} \Sigma_{12} \\
            0 & I_{q_2}
        \end{matrix} \right)
        \left( \begin{matrix}
            \Sigma_{11}^{-1} & 0 \\
            0 & \Sigma_{2|1}^{-1}
        \end{matrix} \right)
        \left( \begin{matrix}
           I_{q_1} & 0 \\
            -\Sigma_{21}\Sigma_{11}^{-1} & I_{q_2}
        \end{matrix} \right).
    \end{align*}
    Using this, we can decompose
    \begin{align*}
        X^\top \Sigma^{-1} X
        =
        X_1^\top \Sigma_{11}^{-1} X_1
        + (X_2 - \Sigma_{21}\Sigma_{11}^{-1}X_1)^\top \Sigma_{2|1}^{-1} (X_2 - \Sigma_{21}\Sigma_{11}^{-1}X_1).
    \end{align*}
    Here we have $(X_2 - \Sigma_{21}\Sigma_{11}^{-1}X_1)$ follows $N_{q_2}(0, \Sigma_{2|1})$ and is independent of $X_1$.
    Consequently, the two quadratic terms are independent conditional on $(X_1>0)$ and the conditional distribution of the second term does not change upon conditioning on $(X_1>0)$. Hence
\begin{align*}
X^\top \Sigma^{-1} X \mid (X_{1}>0) \stackrel{D}{=} 
X_1^\top \Sigma_{11}^{-1} X_1 \mid (X_{1}>0)
+(X_2 - \Sigma_{21}\Sigma_{11}^{-1}X_1)^\top \Sigma_{2|1}^{-1} (X_2 - \Sigma_{21}\Sigma_{11}^{-1}X_1).
\end{align*}    
By independence we can focus on the individual
    distributions in the above.
    First, since $Z_{1} = \Sigma_{11}^{-1/2}X_1 \sim N_{q_1}(0, I_{q_1})$ and $\Sigma_{11}^{1/2}$ is positive definite. 
by the second claim of Lemma~\ref{lemma:trun_chi} we have    
    $$X_1^\top \Sigma_{11}^{-1} X_1 \mid (X_1>0) =
    (\Sigma_{11}^{-1/2}X_{1})^{\top} (\Sigma_{11}^{-1/2}X_1)\mid (X_1>0)
    =Z_{1}^{\top}Z_{1}\mid \Sigma_{11}^{1/2}Z_{1}
    \sim \chi^2_{q_1} 
    $$ 
    Second, $(X_2 - \Sigma_{21}\Sigma_{11}^{-1}X_1) \sim N_{q_2}(0, \Sigma_{2|1})$ we have
    $$(X_2 - \Sigma_{21}\Sigma_{11}^{-1}X_1)^\top \Sigma_{2|1}^{-1} (X_2 - \Sigma_{21}\Sigma_{11}^{-1}X_1) \sim \chi^2_{q_2}.$$
    Thus, conditional on $(X_1>0)$, the two quadratic terms are independent and follow $\chi^2_{q_1}$ and $\chi^2_{q_2}$, respectively. Hence,
    \begin{align*}
        X^\top \Sigma^{-1} X \mid (X_1>0) \sim \chi^2_{q_1+q_2}.
    \end{align*}
    This gives the desired result. 
\end{proof}

\section{Sufficient conditions for validity of
 (\ref{eq:rate2})}

We now state conditions that ensure 
 (\ref{eq:rate2}) holds. These
conditions are only sufficient conditions and can be relaxed.
Let $\|\cdot\|_{2}$  
denote the Euclidean norm of a vector, $\|\cdot\|_{\operatorname{spec}}$ 
and $\textrm{vec}(\cdot)$ denote the spectral norm of a matrix and 
vectorised matrix.
We recall that 
$\Theta = [0,M_1)\times \cdots \times [0,M_b)\times(-M_{1,b+1},M_{2,b+1})\times \cdots \times (M_{1,p},M_{2,p})$ is a compact subset where the population parameter $\theta_{0}\in \Theta$ and can lie
on the boundary. Our important assumption is that there exists an open set $\widetilde{\Theta}$, such that $\Theta \subset \widetilde{\Theta}$ ($\theta_{0}$ lies within the interior of $\widetilde{\Theta}$), where the loss function $\mathcal{L}_{n}(\theta)$ satisfies the classical inferential assumptions. For example, if 
$\widetilde{\theta}_{L} = \arg\min_{\theta \in \widetilde{\Theta}}\mathcal{L}_{n}(\theta)$ then under the standard assumptions we have
that
$\widetilde{\theta}_{L} - \theta_0 = O_{p}(n^{-1/2})$.

\begin{assumption}\label{assum:main}
  \begin{itemize}
 \item[(i)] The parameter space $\widetilde{\Theta}$ is compact
   and convex.
\item[(ii)] For all $n$,
  $\mathcal{L}(\theta)=\E[\mathcal{L}_{n}(\theta)]$ and let 
  $\theta_0=\arg\min_{\theta\in\widetilde{\Theta}} \mathcal{L}(\theta)$, where
  $\theta_0$ lies in the interior of $\widetilde{\Theta}$ and is unique.  Both $\nabla \mathcal{L}(\theta)$ and
  $\nabla^{2}\mathcal{L}(\theta)$ are well defined at $\theta_0$
\item[(iii)] 
  (a) $\nabla \mathcal{L}(\theta_0)=0$ 
  ($\Ex[s_{n}] = \nabla \mathcal{L}(\theta_0) = 0$)
  and (b)
  $W = 
  \nabla^{2}\mathcal{L}(\theta_0)$ is a strictly positive definite
  matrix.

\item[(iv)] $\|\var[\nabla\mathcal{L}_n(\theta_0))]\|_{\operatorname{spec}} = O(n^{-1})$  
\item[(v)]  The second order partial derivatives 
$\nabla^{2} \mathcal{L}_{n}(\theta)$ are continuous on 
$\widetilde{\Theta}$, and at $\theta_0$ we have 
$\|\nabla^{2}\mathcal{L}_n(\theta_0)
  -\nabla^{2}\mathcal{L}(\theta_0)\|_{\operatorname{spec}}=O_{p}(n^{-1/2})$,
  and
  $\|\var[ \operatorname{vec}( \nabla^{2}\mathcal{L}_n(\theta_0))]\|_{\operatorname{spec}} = O(n^{-1})$. 
\item[(vi)]   $\sup_{\theta\in \widetilde{\Theta}}|\mathcal{L}_n(\theta) -
 \mathcal{L}(\theta)|\Pcon 0$.
\item[(vii)] 
For each $j_1,j_2,j_3\in \{1,\ldots,p\}$, the third order derivative 
$\frac{\partial^{3}
\mathcal{L}_{n}(\theta)}{\partial \theta_{j_1}\partial \theta_{j_2}\partial \theta_{j_3}}$ is continuous on 
$\widetilde{\Theta}$. Further
there exists a 
random variable $K_{n,j}$ where 
  \begin{eqnarray*}
    \sup_{\theta\in \widetilde{\Theta}}\Big\| \frac{\partial
    \nabla^{2}\mathcal{L}_{n}(\theta)}{\partial \theta_{j}}
    \Big\|_{\operatorname{spec}} \leq K_{n,j}
  \end{eqnarray*}
  such that $\sup_{n}\E[K_{n,j}]<\infty$ and $\theta =
  (\theta_{1},\ldots,\theta_{p})^\top$.
  \end{itemize}  
\end{assumption}

In the proofs we use the notation that $\lambda_{\min}(W)$ denotes the smallest eigenvalue of $W$.

\begin{lemma}\label{lemma:BS}
Let $Q_{n}(\cdot)$ be defined as in (\ref{eq:Qgeneral}). 
Suppose Assumption \ref{assum:main}(i,ii,iii,v,vii)  holds and  $\mathcal{I}_{\varepsilon} = 
\{\theta: \|\theta-\theta_0\|_{2}\leq \varepsilon\}\subset \widetilde{\Theta}$.  Then 
  \begin{align*}
\sup_{\theta\in \mathcal{I}_{\varepsilon}} 
|\mathcal{L}_{n}(\theta)-Q_{n}(\theta-\theta_0)|
    =O_{p}\left(\varepsilon^{2}n^{-1/2}+\varepsilon^{3}\right).
\end{align*}    
\end{lemma}
\begin{proof}
Assumption \ref{assum:main}(v) ensures that a second order mean value 
expansion of $\mathcal{L}_{n}(\theta)$ about $\theta_0$ exists. This facilitates the decomposition of $\mathcal{L}_{n}(\theta)$  as 
\begin{align}\label{eq:Q_approx}
  \mathcal{L}_n(\theta) =Q_{n}(\delta) +R_{n}(\delta,\theta_0),
\end{align}
where $\delta=\theta-\theta_0$ and 
\begin{align*}
  R_{n}(\delta,\theta_0) &=
\frac{1}{2} \delta^\top\left[\nabla^2 \mathcal{L}_n(\overline{\theta}) -
                           \nabla^{2}\mathcal{L}(\theta_0)\right]\delta
  =\frac{1}{2}(B_{n}(\delta,\theta_0)+S_{n}(\delta,\theta_0))
\end{align*}
with 
\begin{align*}
B_{n}(\delta,\theta_0)
  &=\delta^{\top}[\nabla^{2}\mathcal{L}_{n}(\theta_0)-
    \nabla^{2}\mathcal{L}(\theta_0)]\delta \textrm{ and }
  S_{n}(\delta,\theta_0) 
  = \delta^{\top}
  [\nabla^{2}\mathcal{L}_{n}(\overline{\theta})-
    \nabla^{2}\mathcal{L}_n(\theta_0)]\delta
\end{align*}
where $\bar{\theta}\in \widetilde{\Theta}$ is a convex combination of $\theta$ and $\theta_0$. Again, 
under Assumption \ref{assum:main}(v) we have
\begin{align}\label{eq:Bbound}
\sup_{\|\theta-\theta_0\|_2 \leq \varepsilon} |B_{n}(\delta,\theta_0)| &\leq 
\varepsilon^{2}\|\nabla^{2}\mathcal{L}_{n}(\theta_0)-\nabla^{2}\mathcal{L}(\theta_0)\|_{\operatorname{spec}}
\leq \varepsilon^{2}O_{p}(n^{-1/2}).
\end{align}
Finally, we consider the stochastic term $S_{n}(\delta, \theta_0)$. 
Under Assumption \ref{assum:main}(vii) we can apply 
the mean value theorem to each component of the second derivative
$\nabla^{2}\mathcal{L}_{n}(\theta_0)$, this gives the bound 
\begin{align*}
    |S_n(\delta, \theta_0)| 
     &\leq  C \|\delta\|_{2}^{2}\sum_{j=1}^{p}
                           \left\|\frac{\partial \nabla^2
                                 \mathcal{L}_n(\bar{\bar{\theta}}_j)}{\partial
                           \theta_{j}}\right\|_{\operatorname{spec}}
                           |\overline{\theta}_j - \theta_{0,j}|,
\end{align*}
where $C>0$ is a generic constant and
$\bar{\bar{\theta}}_j$ is a convex combinations of $\overline{\theta}$ and $\theta_0$. Thus by 
applying Assumption \ref{assum:main}(vii) we have 
\begin{align}
  \label{eq:Sbound}
\sup_{\|\theta-\theta_0\|_2\leq \varepsilon}|S_{n}(\delta,\theta_0)|
  &\leq  O_{p}(\varepsilon^{3}).
\end{align}
(\ref{eq:Bbound}) and (\ref{eq:Sbound}) give the result.
\end{proof}  

\vspace{2mm}
Let 
\begin{align}
\label{eq:LQdefin}
\widehat{\theta}_{L} = \arg\min_{\theta\in \Theta}
  \mathcal{L}_{n}(\theta) \textrm{ and }
\widehat{\delta}_{Q} = \arg\min_{\delta\in \mathcal{C}_{\xi}},
  Q_{n}(\delta)
\end{align}
with the above we define $\widehat{\theta}_{Q} = \widehat{\delta}_{Q}+\theta_0$.
In the following lemmas we obtain the classical $n^{1/2}$-convergence rates for the boundary constrained parameter estimators. 

\begin{lemma}\label{lemma:rateLS}
Let $s_{n} = \nabla \mathcal{L}_{n}(\theta_0)$ and 
$W = \nabla^{2}\mathcal{L}(\theta_0)$. 
    Suppose $\Ex[s_n]=0$,  $\| \var[s_n]\|_{\operatorname{spec}} = O(n^{-1})$ and $W$ is strictly positive definite. Then
    \begin{align*}
        \| \widehat{\theta}_Q-\theta_0\|_2 = O_p(n^{-1/2}).
    \end{align*}
\end{lemma}
\begin{proof}
Under the stated conditions we have $s_{n} = O_{p}(n^{-1/2})$, further to simplify notation we let $s=s_{n}$. 

To prove $\widehat{\theta}_{Q}\Pcon \theta_{0}$ we use Theorem
\ref{theorem:blocks}, from Section \ref{sec:explicitexpression}. This decomposition is completely independent of the conditions on 
$\mathcal{L}_{n}(\cdot)$ besides that $W$ is a positive definite matrix. From this decomposition we have
\begin{align*}
\widehat{\delta}  =   \widehat{\theta}_{Q} -\theta_{0}
     &= \sum_{C\in \mathcal{P}(B)}
         \left(
       \begin{matrix}
         \xi_{A_C} \\
         \delta^{F_{C}}
       \end{matrix}  
       \right)
I\left[\delta^{F_C}_{C} > \xi_{C} \right]
        I\left[
        W_{A_{C}|F_{C}} (\widetilde{\delta}_{A_{C}} - \xi_{A_C}) \leq 0\right],
  \end{align*}
where $\delta^{F_C} = - W_{F_CF_C}^{-1} (s_{F_{C}} + W_{F_C A_C} \xi_{A_C})$.
When $\xi_B = 0$, then 
$\delta^{F_C} = - W_{F_CF_C}^{-1}s_{F_{C}}$ and $\|\delta^{F_C}\|_2 = O_p(n^{-1/2})$ (since $\Ex[s]=0$ under Assumption \ref{assum:main}(iii)).
Consequently using the above decomposition we have 
\begin{align*}
    \|\widehat{\theta}_Q-\theta_0\|_2
    \leq \sum_{C \in \mathcal{P}(B)} \|\delta^{F_C}\|_2
    = O_p(n^{-1/2}).
\end{align*}
On the other hand, when $\xi_B \neq 0$,
$\delta^{F_C} = - W_{F_CF_C}^{-1} (s_{F_{C}} + W_{F_C A_C} \xi_{A_C})$ 
need not converge to zero.
Hence we focus on the indicators.
Under Assumption~\ref{assum:main}(iii), $W_{A_C|F_C}$ is positive definite. 
Also the unrestricted estimator $\widetilde{\delta}=-W^{-1}s = O_p(n^{-1/2})$.
Therefore, applying Lemma~\ref{lemma:convergence} we have
\begin{align*}
    \operatorname{pr} [W_{A_C|F_C}(\widetilde{\delta}_{A_C}-\xi_{A_C}) \leq 0] \Pcon 0.
\end{align*}
Using, again, the decomposition in the above
\begin{align*}
    \sqrt{n}
    \left(
       \begin{matrix}
         \xi_{A_C} \\
         \delta^{F_{C}}
       \end{matrix}  
       \right)
I\left[\delta^{F_C}_{C} > \xi_{C} \right]
        I\left[
        W_{A_{C}|F_{C}} (\widetilde{\delta}_{A_{C}} - \xi_{A_C}) \leq 0\right]
        \Pcon 0,
\end{align*}
implying $\|\widehat{\theta}_Q-\theta_0\|_2 = O_p(n^{-1/2})$.
Note that this proof utilizes that an explicit expression for
$\widehat{\theta}_{Q}$ exists. Hence we  do not require that the
sample space $\mathcal{C}$ is compact. 
\end{proof}

\begin{lemma}\label{lemma:consistency}
  Suppose Assumption \ref{assum:main} holds. 
  Then
  \begin{align}
  \label{eq:consistency}  
    \widehat{\theta}_{L}\Pcon \theta_{0} 
  \end{align}
  as $n\rightarrow \infty$ and 
  \begin{align}
  \label{eq:rates}  
    \|\widehat{\theta}_{L} - \theta_{0}\|_2 = O_p(n^{-1/2}).
  \end{align}
\end{lemma}
\begin{proof}
Under Assumption
\ref{assum:main}(i,ii,vi) 
$\widehat{\theta}_{L}\Pcon \theta_{0}$ follows from \cite{van2000asymptotic} Theorem~5.7. This gives convergence in probability, next we require the rate. 

To prove (\ref{eq:rates}) we use an argument similar to 
Lemma 1, \cite{chernoff1954distribution}. 
The proof uses the following equations 
\begin{align}
Q_{n}(0) &= \mathcal{L}_{n}(\theta_0),
 \qquad Q_{n}(\widehat{\delta}_{L}) = 
 \mathcal{L}_{n}(\theta_0)+\widehat{\delta}_{L}^{\top}s + 
 \frac{1}{2}\widehat{\delta}_{L}^{\top}W\widehat{\delta}_{L}, \nonumber\\
 \mathcal{L}_{n}(\widehat{\theta}_{L}) &=
Q_{n}(\widehat{\delta}_{L})+R_{n}(\widehat{\delta}_{L},\theta_{0})
\label{eq:LQexpansions}
\end{align}
and the fact that under Assumption \ref{assum:main}(iii,iv), $\|s\|_{2} = O_{p}(n^{1/2})$.

Using (\ref{eq:LQexpansions}) and that 
$\widehat{\theta}_{L}$ is the minimiser of
$\mathcal{L}_{n}$ over $\Theta$ we have 
\begin{align*}
0\leq \mathcal{L}_{n}(\theta_0)
-\mathcal{L}_{n}(\widehat{\theta}_{L})
= - \widehat{\delta}_{L}^{\top}s - 
 \frac{1}{2}\widehat{\delta}_{L}^{\top}W\widehat{\delta}_{L} - 
R_{n}(\widehat{\delta}_{L},\theta_{0}).
\end{align*}
The difference $\mathcal{L}_{n}(\theta_0)
-\mathcal{L}_{n}(\widehat{\theta}_{L})$ is just an intermediate step, what is important in the above  is the lower bounded by zero;
\begin{align*}
0\leq 
- \widehat{\delta}_{L}^{\top}s - 
 \frac{1}{2}\underbrace{\widehat{\delta}_{L}^{\top}W\widehat{\delta}_{L}}_{\geq 0} - 
R_{n}(\widehat{\delta}_{L},\theta_{0})
\leq\|\widehat{\delta}_{L}\|_2\|s\|_{2} - 
 \lambda_{\min}(W)\frac{1}{2}\|\widehat{\delta}_{L}\|_{2}^{2} - 
|R_{n}(\widehat{\delta}_{L},\theta_{0})|,
\end{align*}
where $\lambda_{\min}(W)$ denotes the smallest eigenvalue of $W$.
Thus 
\begin{align*}
\|\widehat{\delta}_{L}\|_2\|s\|_{2} \geq&
 \lambda_{\min}(W)\frac{1}{2}\|\widehat{\delta}_{L}\|_{2}^{2} +
|R_{n}(\widehat{\delta}_{L},\theta_{0})| \\
\geq& \|\widehat{\delta}_{L}\|_{2}^{2}
\left(
\frac{1}{2}\underbrace{\lambda_{\min}(W)}_{>0}+
C(n^{-1/2}+\|\widehat{\delta}_{L}\|_{2})
\right)
\end{align*}
where $C$ is a positive random variable with 
$C=O_{p}(1)$; the above follows from
$|R_{n}(\widehat{\delta}_{L},\theta_{0})| = 
\|\widehat{\delta}_{L}\|_{2}^{2}O_{P}(n^{-1/2}+\|\widehat{\delta}_{L}\|_{2})$
(see Lemma~\ref{lemma:BS}).
Finally, from (\ref{eq:consistency}) we have
$\|\widehat{\delta}_{L}\|_{2}\Pcon 0$ and 
we recall that 
$\|s\|_{2} = O_{p}(n^{-1/2})$. Thus the above inequality forces 
$\|\widehat{\delta}_{L}\|_{2} = O_p(n^{-1/2})$, which proves the assertion. 
\end{proof}

Having derived individual rates for $\widehat{\theta}_{L}$ and
$\widehat{\theta}_{Q}$ we show that $\widehat{\theta}_{Q}$ is an
approximation of $\widehat{\theta}_{L}$ with an approximation error
smaller than $O_{p}(n^{-1/2})$.

\begin{lemma}\label{lemma:rate}
Let $\widehat{\theta}_{L}$ and $\widehat{\theta}_{Q}$ be defined as in (\ref{eq:LQdefin}), where
 $W=\nabla^{2}\mathcal{L}(\theta_0)$ is a positive definite matrix. 
 Suppose the following conditions 
  \begin{align}\label{eq:uniform_approx}
      \sup_{\|\theta-\theta_0\|_2 \leq Cn^{-1/2}} | \mathcal{L}_{n}(\theta) - Q_n(\theta-\theta_0)|
      = O_{p}(n^{-3/2})
  \end{align}
$\|\widehat{\theta}_{L} -\theta_{0}\|_{2} = O_{p}(n^{-1/2})$  and
$\|\widehat{\delta}_{Q}\|_{2} = O_{p}(n^{-1/2})$
  hold. Then 
  \begin{align*}
  \widehat{\theta}_{L} = \widehat{\theta}_{Q}+O_{p}(n^{-3/4}).
\end{align*} 
\end{lemma}
\begin{proof}
Our aim is to show that 
\begin{align*}
0\leq
  \frac{1}{2}\|\widehat{\delta}_{L}-\widehat{\delta}_{Q}\|_{2}^{2} \lambda_{\min}(W)\leq
  Q_{n}(\widehat{\delta}_{L})-Q_{n}(\widehat{\delta}_{Q})=O_{p}(n^{-3/2}),
\end{align*}  
where $\lambda_{\min}(W)$ denotes the smallest eigenvalue of $W$. From which the result will follow.

Algebraic manipulations gives 
\begin{align}
 Q_{n}(\widehat{\delta}_{L})-Q_{n}(\widehat{\delta}_{Q})
    &= (\widehat{\delta}_{L}-\widehat{\delta}_{Q})^\top \left( s +
    W(\widehat{\delta}_{Q}-\delta_0) \right) + \frac{1}{2} (\widehat{\delta}_{L}-\widehat{\delta}_{Q})^\top W (\widehat{\delta}_{L}-\widehat{\delta}_{Q})  \nonumber\\
    &= (\widehat{\delta}_{L}-\widehat{\delta}_{Q})^\top \nabla
      Q_n(\widehat{\delta}_{Q}) + \frac{1}{2}
      (\widehat{\delta}_{L}-\widehat{\delta}_{Q})^\top W
      (\widehat{\delta}_{L}-\widehat{\delta}_{Q}). \label{eq:Qdiff}
\end{align}

We now study the above. 
Suppose that $E_{\widehat{\theta}_{Q}}$ denotes all 
feasible directions at $\widehat{\theta}_{Q}$ in 
$\mathcal{C}$. Since  
$\widehat{\theta}_{L}\in \mathcal{C}$ then $\nu = 
\widehat{\delta}_{L}-\widehat{\delta}_{Q}=
\widehat{\theta}_{L}-\widehat{\theta}_{Q} \in E_{\widehat{\theta}_{Q}}$ and
 $D_{\nu}  Q_n(\widehat{\delta}_{Q}) \geq 0$ (since $\widehat{\theta}_{Q}$ is a minimum point); see 
Lemma \ref{lemma:direction}. Substituting this bound into
(\ref{eq:Qdiff}) we have the lower bound
\begin{align}
    Q_{n}(\widehat{\delta}_{L})-Q_{n}(\widehat{\delta}_{Q}) 
  &= D_{\nu} Q_n(\widehat{\delta}_{Q})+ \frac{1}{2} (\widehat{\delta}_{L}-\widehat{\delta}_{Q})^\top W
      (\widehat{\delta}_{L}-\widehat{\delta}_{Q}) \nonumber\\
    &\geq \frac{1}{2} (\widehat{\delta}_{L}-\widehat{\delta}_{Q})^\top W (\widehat{\delta}_{L}-\widehat{\delta}_{Q}) \nonumber\\
    &\geq \frac{1}{2} \| \widehat{\delta}_{L}-\widehat{\delta}_{Q}\|_2^2 \lambda_{\min}(W).
    \label{eq:rateQ1}
\end{align}
On the other hand, since 
$\|\widehat{\theta}_{L}-\theta_0\| = O_p(n^{-1/2})$
and $\|\widehat{\theta}_{Q}-\theta_0\| = O_p(n^{-1/2})$, then 
\begin{align}
 0\leq   Q_{n}(\widehat{\delta}_{L})-Q_{n}(\widehat{\delta}_{Q})
    &= Q_{n}(\widehat{\delta}_{L}) - \mathcal{L}_n(\widehat{\theta}_{L}) 
    + \underbrace{\mathcal{L}_n(\widehat{\theta}_{L}) - \mathcal{L}_n(\widehat{\theta}_{Q})}_{\leq 0}
    + \mathcal{L}_n(\widehat{\theta}_{Q}) -Q_{n}(\widehat{\delta}_{Q}) \nonumber\\
    &\leq Q_{n}(\widehat{\delta}_{L}) - \mathcal{L}_n(\widehat{\theta}_{L}) 
    + \mathcal{L}_n(\widehat{\theta}_{Q}) -Q_{n}(\widehat{\delta}_{Q}) 
    \nonumber\\
    &\leq |Q_{n}(\widehat{\delta}_{L}) - \mathcal{L}_n(\widehat{\theta}_{L}) | + |\mathcal{L}_n(\widehat{\theta}_{Q}) -Q_{n}(\widehat{\delta}_{Q})| \nonumber\\
    & = O_p(n^{-3/2}),\label{eq:rateQ2}
\end{align}
where the second line in the display uses that
$\mathcal{L}_n(\widehat{\theta}_{L}) -
\mathcal{L}_n(\widehat{\theta}_{Q}) \leq 0$ and
the last line uses \eqref{eq:rates} in Lemma~\ref{lemma:consistency} and \eqref{eq:uniform_approx}.
Since 
(\ref{eq:rateQ1}) and (\ref{eq:rateQ2}) give
lower and upper bounds for 
$Q_{n}(\widehat{\delta}_{L})-Q_{n}(\widehat{\delta}_{Q})$ these can be combined to give 
\begin{align*}
    \frac{1}{2} \| \widehat{\theta}_{L}-\widehat{\theta}_{Q}\|_2^2 \lambda_{\min}(W)
    &\leq O_p(n^{-3/2}).
\end{align*}
Finally, since $\lambda_{\min}(W)>0$,  we have 
$\| \widehat{\theta}_{L}-\widehat{\theta}_{Q}\|_2 = O_p(n^{-3/4})$,
thus proving the result. 
\end{proof}

\emph{Proof of (\ref{eq:rate2})} We first prove the first bound in 
equation (\ref{eq:rate2}).
For a large enough $n$ we have that $\mathcal{I}_{n} = \{\theta:\|\theta-\theta_0\| = O(n^{-1/2})\}\subset \widetilde{\Theta}$. Then by using 
Lemma \ref{lemma:BS} we have
 \begin{align*}
\sup_{\theta\in \mathcal{I}_{\delta}} 
|\mathcal{L}_{n}(\theta)-
Q_{n}(\theta-\theta_0)|    =O_{p}\left(n^{-3/2}\right).
\end{align*}   
Sufficient conditions for the second rate in (\ref{eq:rate2}) follow from 
Lemma \ref{lemma:rate} which gives 
  \begin{align*}
  \widehat{\theta}_{L} = \widehat{\theta}_{Q}+O_{p}(n^{-3/4}).
\end{align*} 
Altogether, we have shown that under 
Assumption \ref{assum:main} the rates in (\ref{eq:rate2}) hold.
\hfill $\Box$

\vspace{2mm} We briefly mention that 
\begin{align*}
Q_{n}(\delta) = 
\mathcal{L}_{n}(\theta_0)+
s^{\top}\delta + \frac{1}{2}\delta^{\top}W\delta = 
\mathcal{L}_{n}(\theta_0) +
\frac{1}{2}(\delta-\widetilde{\delta})^{\top}W
(\delta-\widetilde{\delta})-\frac{1}{2}\widetilde{\delta}^{\top}W\widetilde{\delta}
\end{align*}
where $\widetilde{\delta} = W^{-1}s$. The quadratic representation given above is commonly used in the literature on the log-likelihood ratio statistic with boundary constraints.

\section{Proof of results in Section \ref{sec:explicitexpression}}

To prove Theorem \ref{theorem:blocks} we require the following lemma, where
we obtain a decomposition for 
$\widehat{\delta}$ in the case that the boundary of the parameter space is at zero, i.e. 
$\xi = 0$ and therefore $\mathcal{C}_{0} = \mathcal{C} = \R_{+}^{b}\times \R^{p-b}$. We later 
generalize to 
any $\mathcal{C}_{\xi}$.

\begin{lemma}\label{lemma:blocks}
Suppose the conditions in Theorem \ref{theorem:blocks} hold. Recall $\widehat{\delta} = \argmin_{\delta \in
    \mathcal{C}} Q(\delta)$. For any $C \subseteq B$, recall $A_C = B \setminus C$ and $F_C = C \cup U$. When $\xi = 0$, $\delta^{(C)}=(0_{A_{C}}, \delta^{F_C})$, where
$\delta^{F_C} = \argmin_{\delta_{F_C} \in \R^{|F_C|}}Q(0_{A_C}, \delta_{F_C})
=-W_{F_C F_C}^{-1}s_{F_C}$. Then, 
\\
(i) There exists a unique set $C \in \mathcal{P}(B)$ such that $\widehat{\delta} = \delta^{(C)}$. Moreover,  the set $C$ satisfies $\delta_{C}^{F_C}> 0 
    \ \textrm{ and } \
    W_{A_C|F_C} \widetilde{\delta}_{A_C} \leq 0$. 
\\
(ii) Conversely, if there exists a set $C\in \mathcal{P}(B)$
  such that 
  \begin{align*}
     \delta_{C}^{F_C}> 0 
    \ \textrm{ and } \
    W_{A_C|F_C} \widetilde{\delta}_{A_C} \leq 0,
    \textrm{ then }
   \delta^{(C)} = \widehat{\delta}.
  \end{align*}
(iii) For $C \in \mathcal{P}(B)$, let 
\begin{align*}
R_{C} = 
\{a\in \mathbb{R}^{p}: (-W_{F_{C}F_{C}}^{-1}a_{F_{C}})_{C}>0,
  -W_{A_{C}|F_{C}}(W^{-1}a)_{A_{C}}\leq 0\}.
\end{align*}
The sets $\{R_C\}_{C \in \mathcal{P}(B)}$ form a disjoint partition of $\R^p$. Moreover, $\widehat{\delta} = \delta^{(C)}$ if and only if $s \in R_C$.
\end{lemma}
{\bf Proof.} We record that for any $\delta = (\delta_A, \delta_F)$, 
 \begin{align} \label{eqn:nabla_Q}
        \nabla Q(\delta) 
        = \left( 
        \begin{array}{c}
             \nabla_{A} Q(\delta)  \\
             \nabla_{F} Q(\delta) 
        \end{array}
   \right)     =
    \left( 
        \begin{array}{c}
            s_{A} + W_{AA} \delta_{A} + W_{AF} \delta_{F} \\
            s_{F} + W_{FA} \delta_{A} + W_{FF} \delta_{F}
        \end{array}
   \right).  
 \end{align}

(i) By definition, $\widehat{\delta}_B \in \R_+^b$. Let $C = \{j \in B: \widehat{\delta}_j > 0\} \subseteq B$ denote the (unique) set of {\it strictly} positive coordinates in $\widehat{\delta}_B$. Therefore, $\widehat{\delta}_j = 0$ for all $j \in A_C = B \setminus C$, and we can write $\widehat{\delta} = (0_{A_C}, \widehat{\delta}_{F_C})$. Since $\widehat{\delta}$ is the joint minimiser of $Q$, it follows from the coordinatewise optimality of a joint minimiser that $\widehat{\delta}_{F_C} = \argmin_{\delta_{F_C} \in \R_{+}^{|C|}\times \R^{|U|}} Q(0_{A_C}, \delta_{F_C})$. Further, since $\widehat{\delta}_{F_C}$ lies within the interior of $\R_{+}^{|C|}\times \R^{|U|}$, it follows from Lemma \ref{lemma:direction} that $\nabla_{F_C} Q(0_{A_C}, \widehat{\delta}_{F_C})=0$. Using the expression for the gradient from \eqref{eqn:nabla_Q}, it follows that $\widehat{\delta}_{F_C} = -(W_{F_CF_C})^{-1} s_{F_C}$, which by definition is $\delta^{F_C}$. Thus, by virtue of the minimiser being in the interior and the objective function being quadratic, the minimiser over $\R_{+}^{|C|}\times \R^{|U|}$ and over all of $\R^{|F_C|}$ agree. This shows $\widehat{\delta} = (0_{A_C}, \delta^{F_C}) = \delta^{(C)}$.


To complete the proof of this part, observe that $\delta_C^{F_C} = \widehat{\delta}_C > 0$ by assumption. Next, since $\delta^{(C)}=\argmin_{\delta\in \mathcal{C}}Q(\delta)$, Lemma \ref{lemma:direction} gives that $D_{\nu} Q(\delta^{(C)}) \ge 0$ for every feasible direction vector
$\nu$ at $\delta^{(C)}$ in $\mathcal{C}$. Substituting
 $\delta = \delta^{(C)}$ into \eqref{eqn:nabla_Q}, one has 
\begin{align}
     \nabla Q(\delta^{(C)}) =
\left(
\begin{array}{cc}
  \nabla_{A_C} Q(\delta^{(C)}) \\
  \nabla_{F_C} Q(\delta^{(C)}) \\
\end{array}
     \right) =
\left(
\begin{array}{c}
     s_{A_C} - W_{A_CF_C}(W_{F_CF_C})^{-1}s_{F_C}\\
  0 \\
\end{array}
     \right) = 
\left(
\begin{array}{c}
   -W_{A_C|F_C} \widetilde{\delta}_{A_C} \\
  0 \\
\end{array}
     \right), \label{eq:nabladelta}
   \end{align}
where the last line uses that 
 $s = -W\widetilde{\delta}$. Thus, for any $\nu \in \R^p$, 
   \begin{align*}
      D_{\nu} Q(\delta^{(C)}) = \nu^\top \nabla Q(\delta^{(C)})=
\nu_{A_C}^{\top}\nabla_{A_C} Q(\delta^{(C)})+
\nu_{F_C}^{\top}\nabla_{F_C} Q(\delta^{(C)}) = -\nu_{A_C}^{\top} W_{A_C|F_C} \widetilde{\delta}_{A_C}.
\end{align*}
Next, note that all the feasible directions in $\mathcal{C}$ at $\delta^{(C)}$ is
$S = \{ \nu =(\nu_{A_C}, \nu_{F_C}): \nu_{A_C}\in \R_{+}^{|A_C|} \}$. Thus, from the previous display, we have $-\nu_{A_C}^{\top} W_{A_C|F_C} \widetilde{\delta}_{A_C} \ge 0$ for every $\nu_{A_C} \in \R_+^{|A_C|}$. This  implies
that $W_{A_C|F_C} \widetilde{\delta}_{A_C} \leq 0$.
\\[2ex]
(ii) We show that if there exists a $C \in \mathcal{P}(B)$ such that
\begin{align*}
     \delta^{F_C}_{C} > 0 
    \ \textrm{ and } \
    W_{A_C|F_C} \widetilde{\delta}_{A_C} \leq 0,
  \end{align*}
  then $\delta^{(C)} = (0_{A_C}, \delta^{F_C}) = \widehat{\delta}$. To prove this, we treat
  $\delta^{(C)}$  as a candidate for $\widehat{\delta}$ and show that
  $D_{v} Q(\delta^{(C)})\geq 0$ for all feasible directions $v$ in $\mathcal{C}$ at $\delta^{(C)}$; recall the set of all such feasible directions $S = \{ \nu =(\nu_{A_C}, \nu_{F_C}): \nu_{A_C}\in \R_{+}^{|A_C|} \}$. The claim then follows from Lemma
  \ref{lemma:direction}. 

  The rest of the proof mirrors parts of the proof of (i). Recall that for any $\nu \in \R^p$, $D_{\nu} Q(\delta^{(C)}) = -\nu_{A_C}^{\top}W_{A_C|F_C}\widetilde{\delta}_{A_C}$. 
  Since $\delta^{F_C}_{C}>0$, we have  $\delta^{(C)}\in
  \mathcal{C}$ (thus it is a viable candidate). 
 Fix $\nu \in S$. Since by assumption $W_{A_C|F_C}\widetilde{\delta}_{A_C}\leq 0$
 and $\nu_A \in \R^{|A|}_{+}$, 
 then 
\begin{align*}
D_{\nu} Q(\delta^{(C)}) =
  -\nu_{A_C}^{\top}W_{A_C|F_C}\widetilde{\delta}_{A_C}\geq 0.
\end{align*}
This completes the proof of (ii). 
\\[1ex]
(iii) Parts (i) and (ii) together imply there exists a unique set $C \in \mathcal{P}(B)$ such that $\delta_C^{F_C}  > 0$ and $W_{A_C|F_C} \widetilde{\delta}_{A_C} \le 0$. For this set $C$, $\delta^{(C)} = \widehat{\delta}$. Moreover, since $\delta_C^{F_C} = (- W_{F_C F_C}^{-1} s_{F_C})_C$ and $W_{A_C|F_C} \widetilde{\delta}_{A_C} = - W_{A_C|F_C} (W^{-1} s)_{A_C}$, we have $s \in R_C$. It remains to show that the collection of sets $\{R_C\}_{C \in \mathcal{P}(B)}$ form a disjoint partition of $\R^p$. To that end, for each $a \in \R^p$, we let
  $\widehat{\delta}(a) = (\widehat{\delta}_{B}(a), \widehat{\delta}_{U}(a)) =
  \argmin_{\delta \in \mathcal{C}}Q(\delta; a,W)$
  with 
$\widehat{\delta}_{B}(a) \in
\R_{+}^{b}$ and $\widehat{\delta}_{U}(a) \in \R^{p-b}$. Our $\widehat{\delta}$ equals $\widehat{\delta}(s)$ in this notation. 

For each $C \subseteq B$, we define the set
\begin{align*}
\bar{R}_{C} = \{a \in \R^{p}: \widehat{\delta}_{A_C}(a)=0_{A_C} \textrm{ and
  }\widehat{\delta}_{C}(a) \in (0,\infty)^{|C|}\}.
\end{align*}
Since $\widehat{\delta}_{B}(a) \in \R^{b}_{+}$,
there must exist a unique set $C(a) =\{j\in B: \widehat{\delta}_j(a)>0 \} \subseteq B$ 
which determines the strictly positive coordinates in $\widehat{\delta}_B(a)$, while all remaining coordinates in  $A(a) = B \backslash C(a)$ are zero. 
Since every $a \in \R^P$ corresponds to exactly one set $C(a)$, it belongs to exactly one $\bar{R}_C$.
Thus, $\{ \bar{R}_{C} \}_{C \in \mathcal{P}(B)}$ is a disjoint partition
of $\R^{p}$.

We now show $R_C = \bar{R}_C$. Following the proof of part (i) (with $a$ in place of $s$), one obtains $\bar{R}_C \subseteq R_C$. Following the proof of part (ii) (again with $a$ in place of $s$), one obtains $R_C \subseteq \bar{R}_C$. This proves the claim.

We now use the lemma above to prove Theorem \ref{theorem:blocks}. 

\vspace{2mm}
\noindent \emph{Proof of Theorem \ref{theorem:blocks}}
 We first prove the result in the case $\mathcal{C}_{\xi}=\mathcal{C}$. 
 By using Lemma \ref{lemma:blocks} 
 $\{R_{C}\}_{C\in \mathcal{P}(B)}$ forms a
disjoint partition of $\R^{p}$, further Lemma \ref{lemma:blocks}(i,ii) implies 
    \begin{align}
    \label{eq:widehattilde}
    \widehat{\delta}
     &= \sum_{C\in \mathcal{P}(B)}
\left(
\begin{array}{c}
0_{A_C} \\
\delta^{F_C} \\
\end{array}
\right)
I\left[\delta_{C}^{F_{C}}> 0 \right]
        I\left[
        W_{A_{C}|F_{C}} \widetilde{\delta}_{A_{C}} \leq 0\right].
    \end{align}
    where $\delta^{F_C} = -W_{F_CF_C}^{-1}s_{F_C}$ and  
    $\widetilde{\delta} = -W^{-1}s$.
    
To prove the result in the general case that $\xi\neq 0$ we want to obtain an expression for 
\begin{align*}
\widehat{\delta}(s, W, \xi) :\, = \argmin_{\delta \in \mathcal{C}_\xi} Q(\delta;s, W)
\end{align*}
using 
the expression for $\widehat{\delta}(s, W,0)$ in (\ref{eq:widehattilde}), where recall that  $Q(\delta;s, W) = \mathcal{L}_n(\theta_0)+s^{\top}\delta+\frac{1}{2}\delta^{\top}W\delta$. To do this, we 
write $\eta = \delta - \xi$ so that 
\begin{align*}
Q(\delta;s, W) = Q(\eta + \xi; s, W) =
\frac{1}{2} \eta^\top W \eta + \eta^\top (W \xi + s) + Q(\xi; s, W), 
\end{align*}
where $\eta \in \mathcal{C}$. Then, using standard properties of switching min and location shift, we have 
\begin{align*}
\widehat{\delta}(s, W, \xi) = \xi + \widehat{\delta}(W \xi + s, W, 0), 
\end{align*}
where $\xi^{\top} = (\xi_{B}^{\top},0_{U})$ and 
\begin{align*}
W\xi = \left(
\begin{array}{c}
W_{AA}\xi_{A} + W_{AC}\xi_{C} \\
W_{CA}\xi_{A} + W_{CC}\xi_{C} \\
W_{UA}\xi_{A} + W_{UC}\xi_{C} \\
\end{array}
\right) = 
\left(
\begin{array}{c}
W_{AA}\xi_{A}+W_{AC}\xi_{C} \\
W_{FA}\xi_{A} + W_{FC}\xi_{C} \\
\end{array}
\right).
\end{align*}

We now use (\ref{eq:widehattilde}) to obtain an expression for 
$\widehat{\delta}(W \xi + s, W, 0)= \widehat{\delta}(s_{\xi}, W, 0)$
where $s_{\xi} = W \xi + s$. For each subset $F\subset [p]$ 
we define the unconstrained least squares estimator $\delta^{\xi,F}=-W_{FF}^{-1}(s_{\xi})_{F}$. It can be shown that
\begin{align*}
\delta^{\xi,F}&=
-W_{FF}^{-1}(s_{\xi})_{F} = -W_{FF}^{-1}(s_{F} + 
W_{FA}\xi_{A}+W_{FC}\xi_{C}) \\
&= -W_{FF}^{-1}(s_{F}+W_{FA}\xi_{A}) - 
\left(
\begin{array}{c}
\xi_{C} \\
0 \\
\end{array}
\right)
\end{align*}
where the last line is because $W_{FF}^{-1}W_{FC} = (I_{|C|},0_{U}^{\top})$. 
Similarly, we define $\widetilde{\delta}^{\xi} = -W^{-1}s_{\xi}$ and rewrite this as 
\begin{align*}
\widetilde{\delta}^{\xi} = -W^{-1}(s+W\xi) = -W^{-1}s-\xi.
\end{align*}
Substituting $\delta^{\xi,F}$ and $\widetilde{\delta}^{\xi}$ into 
(\ref{eq:widehattilde}) gives
 \begin{align*}
 \widehat{\delta}(W \xi + s, W, 0) &=
\sum_{C\in \mathcal{P}(B)}
\left(
\begin{array}{c}
0 \\
\delta^{\xi,F_C}_{C}\\
\delta^{\xi,F_C}_{U} \\
\end{array}
\right)
I\left[\delta_{C}^{\xi,F_{C}}> 0 \right]
        I\left[
        W_{A_{C}|F_{C}} \widetilde{\delta}^{\xi}_{A_{C}} \leq 0\right]\\ 
     &= \sum_{C\in \mathcal{P}(B)}
\left(
\begin{array}{c}
0 \\
-\xi_{C} + \delta^{F_C}_{C}\\
\delta^{F_C}_{U} \\
\end{array}
\right)
I\left[\delta_{C}^{F_{C}}> \xi_{C} \right]
        I\left[
        W_{A_{C}|F_{C}}(\widetilde{\delta}_{A_C}-\xi_{A_C}) \leq 0\right],
    \end{align*}
    where $\delta^{F_C} = -W_{F_CF_C}^{-1}(s_{F_C}+W_{F_CA_{C}}\xi_{A})$ and  
    $\widetilde{\delta} = -W^{-1}s$.

Altogether we have 
\begin{align*}
 \widehat{\delta}(s, W, \xi) &= \xi + \widehat{\delta}(W\xi+s, W, 0) \\
     &= \sum_{C\in \mathcal{P}(B)}
\left(
\begin{array}{c}
\xi_{A_C} \\
\delta^{F_C}_{C}\\
\delta^{F_C}_{U} \\
\end{array}
\right)
I\left[\delta_{C}^{F_{C}}> \xi_{C} \right]
        I\left[
        W_{A_{C}|F_{C}}(\widetilde{\delta}_{A_C}-\xi_{A_C}) \leq 0\right],
    \end{align*}
    since $\xi^{\top} = (\xi_{A}^{\top},\xi_{C}^{\top},0_{U})$.
    This gives the desired result.
   \hfill $\Box$

 \vspace{3mm}
\noindent\emph{Proof of Corollary \ref{cor:blocks}}   
  We will use the classical block inverse matrix identity
\begin{align}
 W^{-1} &=\left(
  \begin{array}{cc}
    W_{AA} & W_{AF} \\
    W_{FA} & W_{FF}
   \end{array} 
             \right)^{-1}=
 \left(
  \begin{array}{cc}
    W_{A|F}^{-1} & -W_{A|F}^{-1}W_{AF}W_{FF}^{-1} \\
   -W_{FF}^{-1}W_{FA}W_{A|F}^{-1} & W_{FF}^{-1}+W_{FF}^{-1}W_{FA}W_{A|F}^{-1}W_{AF}W_{FF}^{-1}\\
   \end{array} 
  \right) \nonumber\\
  &=\left(
  \begin{array}{cc}
    \Sigma_{AA} & \Sigma_{AF} \\
    \Sigma_{FA} & \Sigma_{FF}
   \end{array} 
             \right)\label{eq:schur}
\end{align}
where $W_{A|F} = W_{AA} -W_{AF}W_{FF}^{-1}W_{FA}  =
\Sigma_{AA}^{-1}$. Thus
\begin{align}
\label{eq:cor1}
  I\left[
        W_{A|F} (\widetilde{\delta}_{A} - \xi_A) \leq 0\right] =  I\left[
        \Sigma_{AA}^{-1}(\widetilde{\delta}_{A}-\xi_{A}) \leq 0\right]
\end{align}
Further, by using (\ref{eq:schur}), we observe that
$W_{F F}^{-1} W_{FA} = -\Sigma_{FA}\Sigma_{AA}^{-1}$, thus 
\begin{align}
\label{eq:cor2}  
  \delta^{F} =
  \widetilde{\delta}_{F}+W_{F F}^{-1} W_{F
  A}(\widetilde{\delta}_{A}-\xi_{A}) =  \widetilde{\delta}_{F}-\Sigma_{FA}\Sigma_{AA}^{-1}
  (\widetilde{\delta}_{A}-\xi_{A}).
\end{align}
Substituting (\ref{eq:cor1}) and (\ref{eq:cor2}) into
(\ref{eq:consequence}) proves the result. 
\hfill $\Box$

\section{Proof of results in Section \ref{sec:distribution}}

\noindent
\emph{Proof of Theorem \ref{theorem:nonstandard}} The result
immediately follows from (\ref{eq:rate2}) and Theorem \ref{theorem:blocks}.
\hfill $\Box$

\vspace{2mm}

\noindent \emph{Proof of Lemma \ref{lemma:LSrateG}}
See Lemma \ref{lemma:rateLS}. \hfill $\Box$

\vspace{2mm}

\begin{remark}[Connections between 
Definitions \ref{def:SN_Variants} and the classical skew normal]
\label{remark:supSN}
    Our one-sided skewed normal distribution $\widetilde{SN}$ includes
    the classical multivariate skewed normal distribution defined in \cite{azzalini1996multivariate}; see their equation~(2.3).
    Specifically, they define the joint density of multivariate skew normal as
    \begin{align*}
        2 \phi_k(x; \Omega) \Phi_1 (\alpha^\top x), \quad x\in \R^k,
    \end{align*}
    where $\Omega$ is a positive definite matrix, 
    $\phi_k(x; \Omega)$ is the density of the $k$ dimensional multivariate normal distribution with a correlation matrix $\Omega$, and
    $\Phi_1$ is the cumulative distribution function of the standard normal distribution.
    The matrix $\Omega$ and vector $\alpha \in \R^k$ are determined by their shape parameter $\lambda$ and dependence parameter $\Psi$ through their equations (2.4) - (2.6).
    This density is obtained from our definition of one-sided skewed normal $\widetilde{SN}$ by setting
    \begin{align*}
        d_1 = 1, \quad d_2 = k, \quad \zeta = 0, \quad E_2 = \Omega,\quad  B_{12} = \alpha^\top, \quad  E_1=1,
    \end{align*}
    since $\widetilde{\Phi}_1(-\alpha^\top x; 1) = \Phi_1(\alpha^\top x)$ and normalizing constant equals to $1/2$ by construction.
    Our definition extends the conditioning method in their Section~2.2 from a scalar Gaussian variable exceeding zero to a Gaussian vector exceeding arbitrary thresholds $\zeta$.
    
    Our one-sided skewed normal distribution is also a reparametrization of the unified skewed normal distribution introduced in \cite{arellano2006unification} with the location parameter set to zero (see their equation~(10)). 
    Specifically, their density in this case is
    \begin{align*}
        \phi_k(x; \Omega) 
        \frac{\Phi_m(\gamma+\Delta^\top \bar{\Omega}^{-1} \omega^{-1} x; \Gamma-\Delta^\top \bar{\Omega}^{-1}\Delta)}{\Phi_m(\gamma; \Gamma)}, 
        \quad x\in \R^k.
    \end{align*}
    where $\Omega = \omega \bar{\Omega} \omega$, $\omega$ is a $k \times k$ positive diagonal matrix, and $\Gamma$ and $\bar{\Omega}$ are $m\times m$ and $k \times k$ positive definite matrices, respectively. Also, $\gamma \in \R^m$ and $\Delta \in \R^{k \times m}$ with $\Gamma - \Delta^\top \bar{\Omega}^{-1} \Delta \succ 0$.
    This density coincides with our one-sided skewed normal under the choices
    \begin{eqnarray*}
        &d_1 = m, \quad d_2 = k, \quad \zeta = -\gamma, \\
        &E_2 =\Omega, \quad B_{12} = \Delta^\top \bar{\Omega}^{-1} \omega^{-1}, \quad  E_1 =\Gamma-\Delta^\top \bar{\Omega}^{-1}\Delta.
    \end{eqnarray*}
    Here $E_1+B_{12}E_{2}B_{12}^\top = \Gamma$, so our normalizing constant is $\Phi_m(\gamma; \Gamma)$.
    Also, our two-sided skewed normal $SN$ can be represented as the unified skewed normal distribution with zero location parameter.
    This follows by reversing the signs of the Gaussian components subject to upper bounds and transforming the covariance matrix accordingly. To be specific,
    \begin{align*}
        SN_{d_4,d_3}(\zeta_1, \zeta_2, B_{13}, B_{23}, E_{3}, E_4)
        =
        \widetilde{SN}_{d_4,d_3}\left(  \left(\begin{matrix}
            -\zeta_1 \\
            \zeta_2
        \end{matrix} \right), 
        \left(\begin{matrix}
            -B_{13} \\
            B_{23}
        \end{matrix} \right),
        E_3, JE_4J^\top\right),
    \end{align*}
    where $J = \operatorname{diag}(-I_{d_1}, I_{d_2})$.
\end{remark}

\noindent \emph{Proof of Lemma \ref{lemma:SN_TSN}}
To prove the result we use the following classical identity. If
\begin{align*}
\left(
\begin{matrix}
X_{1} \\
X_{2} \\
\end{matrix}
\right)\sim 
N\bigg(
\left( \begin{matrix}
    0 \\
    0
\end{matrix} \right),
\left( \begin{matrix}
\Sigma_{11} & \Sigma_{12} \\
\Sigma_{12}^{\top} & \Sigma_{22} \\
\end{matrix} \right)
\bigg),
\end{align*}
then 
\begin{align*}
X_{1}\mid X_{2} = x \sim N\left(\Sigma_{12}\Sigma_{22}^{-1}x,\Sigma_{1|2}\right). 
\end{align*}
Hence to prove (i) we use that 
$Y_{1}\mid Y_{2}=y_{2} = N\left(V_{12}V_{22}^{-1}y_{2},V_{1|2}\right)$, thus 
\begin{align*}
\operatorname{pr} (Y_{1}\leq \zeta_{1}\mid Y_{2}=y_{2}) = 
\Phi_{d_1}(\zeta_1 -V_{12}V_{22}^{-1}y_{2};V_{1|2}).
\end{align*}
Altogether, this gives 
\begin{align*}
f_{Y_{2}\mid Y_{1}<\zeta_{1},Y_{2}>\zeta_{2}}(y_{2}) 
&=
\frac{\operatorname{pr}(Y_{1}<\zeta_1\mid Y_{2}=y_2)
\phi_{d_2}(y_{2};V_{22})}{\operatorname{pr}(Y_{1}<\zeta_1,Y_{2}>\zeta_2)}I(y_{2}>\zeta_2) \\
&= \frac{\Phi_{d_2}(\zeta_1-V_{12}V_{22}^{-1}y_{2};V_{1|2})\phi_{d_2}(y_{2};V_{22})
I(y_{2}>\zeta_2)}{
\int \Phi_{d_2}(\zeta_1-V_{12}V_{22}^{-1}y_{2};V_{1|2})\phi_{d_2}(y_{2};V_{22})
I(y_{2}>\zeta_2)dy_2
},
\end{align*}
which by definition gives $Y_{2}\mid (Y_{1}<\zeta_1,Y_{2}>\zeta_2)\sim 
tSN_{d_1,d_2}(\zeta_1,\zeta_2,V_{12}V_{22}^{-1}, V_{22}, V_{1|2})$. 
In case $V_{12}=0$, $Y_1$ and $Y_2$ are independent. So $Y_{2}\mid (Y_{1}<\zeta_1,Y_{2}>\zeta_2)$ is equivalent to $Y_{2}\mid (Y_{2}>\zeta_2)$ and is
\begin{align*}
    f_{Y_{2}\mid (Y_{2}>\zeta_2)}(y_2) 
    &= \frac{\phi_{d_2}(y_2;V_{22})}{\operatorname{pr}(Y_2>\zeta_2)} I(y_2>\zeta_2),
\end{align*}
which by definition is $tN_{d_2}(\zeta_2, V_{22})$.

Now we prove (ii). Using a similar argument to the above, 
if
\begin{align*}
\left(
\begin{matrix}
X_{1} \\
X_{2} \\
X_{3}
\end{matrix}
\right)\sim 
N\bigg(
\left( \begin{matrix}
    0 \\
    0 \\
    0
\end{matrix} \right),
\left( \begin{matrix}
\Sigma_{11} & \Sigma_{12} & \Sigma_{13} \\
\Sigma_{12}^{\top} & \Sigma_{22} & \Sigma_{23} \\
\Sigma_{13}^{\top} & \Sigma_{23}^\top & \Sigma_{33} 
\end{matrix} \right)
\bigg),
\end{align*}
then
we have 
\begin{align*}
\left(
\begin{matrix}
X_{1} \\
X_{2} \\
\end{matrix}
\right )\mid  X_{3}=x_{3}\sim 
N\left(
\left(
\begin{matrix}
\Sigma_{13}\Sigma_{33}^{-1}x_{3} \\
\Sigma_{23}\Sigma_{33}^{-1}x_{3} \\
\end{matrix}
\right),
\Sigma_{(1,2)|3}
\right),
\end{align*}
where $\Sigma_{(1,2)|3} = \Sigma_{(1,2),(1,2)} - (\Sigma_{13}^\top, \Sigma_{23}^\top)^\top \Sigma_{33}^{-1} (\Sigma_{31}, \Sigma_{32})$ and $\Sigma_{(1,2),(1,2)}$ denote the $(1,2)$ sub-diagonal of the matrix $\Sigma$.
Thus 
 \begin{align*}
f_{Y_{3}\mid Y_{1}<\zeta_1,Y_{2}>\zeta_2}(y_{3})&= 
\frac{\operatorname{pr}(Y_{1}<\zeta_1,Y_{2}>\zeta_2\mid Y_{3}=y_3)\phi_{d_3}(y_{3};V_{33})}{\operatorname{pr}(Y_{1}<\zeta_1,Y_{2}>\zeta_2)} \\
&= \frac{\Phi_{d_1+d_2}^{*}( \zeta_1 - V_{13}V_{33}^{-1}y_3, \zeta_2 - V_{23}V_{33}^{-1}y_3; V_{(1,2)|3} ) \phi_{d_3}(y_3;V_{33})}{\int \Phi_{d_1+d_2}^{*}( \zeta_1 - V_{13}V_{33}^{-1}y_3, \zeta_2 - V_{23}V_{33}^{-1}y_3; V_{(1,2)|3} ) \phi_{d_3}(y_3;V_{33}) dy_3},
\end{align*}
which by definition gives the two sided skewed normal $SN(\zeta_1, \zeta_2, V_{13}V_{33}^{-1}, V_{23}V_{33}^{-1}, V_{33}, V_{(1,2)|3})$.

In case $V_{13}=0$, $Y_1$ and $Y_3$ are independent. Hence $Y_3\mid (Y_1<\zeta_1, Y_2>\zeta_2)$ is equivalent to $Y_3\mid (Y_2>\zeta_2)$ and is
\begin{align*}
    f_{Y_3 \mid (Y_2>\zeta_2)}(Y_3) 
    &= \frac{\operatorname{pr}(Y_2 > \zeta_2 \mid Y_3=y_3) \phi_{d_3}(y_3;V_{33})}{\operatorname{pr}(Y_2>\zeta_2)} \\
    &= \frac{\widetilde{\Phi}_{d_2}(\zeta_2-V_{23}V_{33}^{-1}y_3;V_{2|3})\phi_{d_3}(y_3;V_{33})}{\int \widetilde{\Phi}_{d_2}(\zeta_2-V_{23}V_{33}^{-1}y_3;V_{2|3})\phi_{d_3}(y_3;V_{33}) dy_3},
\end{align*}
which by definition is the one-sided skewed normal $\widetilde{SN}_{d_2,d_3}(\zeta_2, V_{23}V_{33}^{-1},V_{33}, V_{2|3})$.
\hfill $\Box$

\begin{theorem}\label{theorem:CLTpre}
    Suppose (\ref{eq:rate2}) and 
    Assumption \ref{assum:CLT} holds. 
    Then using the notation in Theorem ~\ref{theorem:blocks} for each 
    $C\in \mathcal{P}(B)$ we have
    \begin{align*}
        \sqrt{n} 
        \left(\begin{matrix}
            W_{A_C|F_C} \widetilde{\delta}_{A_C} \\
            \delta^{F_C}_C \\
            \delta^{F_C}_U
        \end{matrix} \right)
        \Dcon N(0, V_{C}),
    \end{align*}
    where $V_{C} =  T_{C} \Sigma T_{C}^\top$ and 
    \begin{align*}
        T_{C} &= 
        \left( \begin{matrix}
            W_{A_C|F_C} & 0 & 0\\
            W_{C|U}^{-1}W_{CA_C|U} & I_{|C|} & 0 \\
         W_{U|C}^{-1}W_{UA_C|C}   & 0 & I_{|U|}
        \end{matrix}\right),
    \end{align*}
with    $W_{CA_C|U} = W_{CA_C}-W_{CU}W_{UU}^{-1}W_{UA_C}$,
    $W_{UA_C|C}=W_{UA_C}-W_{UC}W_{CC}^{-1}W_{CA_C}$,
    $A_C = B\backslash C$ and $F_C = C\cup U$.  Further 
    \begin{itemize}
        \item[(i)] Asymptotically, the conditional distribution of $\sqrt{n}\delta^{F_C}_{C}$ given
    $\{ (\delta^{F_C}_C > 0), (W_{A_C |F_C}\widetilde{\delta}_{A_C} \leq 0) \}$
    is 
    \begin{align*}
        tSN_{|A_C|,|C|}(0,0,  V_{A_C C}V_{CC}^{-1}, V_{C C}, V_{A_{C}|C}),
    \end{align*}
where $V_{A_{C}|C} = V_{A_CA_C} - V_{A_C C}(V_{CC})^{-1}V_{CA_C}$.
        \item[(ii)] Asymptotically, 
        the conditional distribution of $\sqrt{n}\delta^{F_C}_{U}$ given
    $\{ (\delta^{F_C}_C > 0), (W_{A_C|F_{C}} \widetilde{\delta}_{A_C} \leq 0) \}$ is
    \begin{align*}
        SN_{b, p-b}\left(0, 0,
        V_{A_C U}V_{UU}^{-1}, V_{CU}V_{UU}^{-1}, V_{UU}, V_{B|U} \right),
    \end{align*}
    \end{itemize}
    where 
    \begin{align*}
       V_{B|U} = 
        \left( \begin{matrix}
            V_{A_C A_C} & V_{A_C C}  \\
            V_{C A_C} & V_{C C} 
        \end{matrix} \right) - 
        \left( \begin{matrix}
             V_{A_C U} \\
            V_{C U}
        \end{matrix} \right) 
        V_{UU}^{-1}
        \left( \begin{matrix}
             V_{U A_C} & V_{U C}
        \end{matrix} \right).
    \end{align*}
 \end{theorem}
\begin{proof}
    We first obtain the joint distribution of $(
      W_{A_{C}|F_{C}}\widetilde{\delta}_{A_C}, \delta^{F_C}_C,
      \delta^{F_C}_U)$. To simplify notation we drop the
      $C$ in $A_{C}$ and $F_{C}$.
      From Remark \ref{remark:SWinv}
we have $\delta^{F} =
  \widetilde{\delta}_{F}+(W_{F F})^{-1} W_{F
  A}\widetilde{\delta}_{A}$, thus by using this and the block matrix inversion 
      identity we have
\begin{align}
  \left(
  \begin{matrix}
    \delta^{F}_{C} \\
       \delta^{F}_{U} \\
   \end{matrix} 
  \right) &=
   \left(
  \begin{matrix}
    \widetilde{\delta}_{C} \\
       \widetilde{\delta}_{U} \\
   \end{matrix} 
  \right) +
  \left(
  \begin{matrix}
    W_{CC} & W_{CU} \\
    W_{UC} & W_{UU}
   \end{matrix} 
             \right)^{-1}
              \left(
  \begin{matrix}
    W_{CA} \\
    W_{UA} \\
   \end{matrix} 
  \right) \widetilde{\delta}_{A} \nonumber \\
&=   \left(
  \begin{matrix}
    \widetilde{\delta}_{C} \\
       \widetilde{\delta}_{U} \\
   \end{matrix} 
  \right) +
  \left(
  \begin{matrix}
    W_{C|U}^{-1} & -W_{C|U}^{-1}W_{CU}W_{UU}^{-1} \\
    -W_{U|C}^{-1}W_{UC}W_{CC}^{-1} & W_{U|C}^{-1}
   \end{matrix} 
             \right)
              \left(
  \begin{matrix}
    W_{CA} \\
    W_{UA} \\
   \end{matrix} 
  \right) \widetilde{\delta}_{A} \label{eq:deltaF_CU}
\end{align}
Using the above and the definition of $T_{C}$ it is easily seen that 
    \begin{align*}
        T_{C} \widetilde{\delta}
        &=         \left(\begin{matrix}
            W_{A|F} \widetilde{\delta}_{A} \\
            \delta^{F}_C \\
            \delta^{F}_U
        \end{matrix} \right).
    \end{align*}
    Since $\sqrt{n} \widetilde{\delta} \Dcon N(0, \Sigma)$, 
    \begin{align*}
        \sqrt{n} T \widetilde{\delta}
        \Dcon
        N(0, V_{C}) \textrm{ with } V = T_{C} \Sigma T_{C}^{\top}.
    \end{align*}
By using  Lemma~\ref{lemma:SN_TSN} we obtain the conditional
distributions in (i) and (ii).
\end{proof}

\vspace{1mm}
\noindent \emph{Proof of Theorem \ref{theorem:CLT}}
By using Theorems \ref{theorem:blocks} and \ref{theorem:CLTpre} we have
\begin{align*}
  \sqrt{n}\widehat{\delta}\Dcon 
    \sum_{C\in \mathcal{P}(B)}
    \left(
    \begin{matrix}
      m_{A_{C}}(0)\\
        tSN_{|A_C|,|C|}(0,0,  V_{A_C C}V_{CC}^{-1}, V_{C C}, V_{A_{C}|C})\\
     SN_{b, p-b}\left(0, 0,
        V_{A_C U}V_{UU}^{-1}, V_{CU}V_{UU}^{-1}, V_{UU}, V_{B|U} \right)
      \end{matrix}
    \right)
       \operatorname{pr} \big( \delta^{F_C}_C > 0, W_{A_{C}|F_{C}} \widetilde{\delta}_{A_C} \leq 0 \big).
\end{align*}
The above together with (\ref{eq:Lest}) and that $\xi_{B}=0$ gives the result.
\hfill $\Box$

\begin{lemma}\label{cor:WSigmaInve}
  Suppose (\ref{eq:rate2}) and 
    Assumption \ref{assum:CLT} holds, 
    where $\Sigma = W^{-1}$. Then for each $C \in \mathcal{P}(B)$ we have
    \begin{align}
      \label{eq:Waa}
        \sqrt{n} 
        \left(\begin{matrix}
            W_{A_C|F_C} \widetilde{\delta}_{A_C} \\
            \delta^{F_C}_C \\
            \delta^{F_C}_U\\
        \end{matrix} \right)
        \Dcon N(0, V_C),
    \end{align}
where 
    \begin{align*}
      V_{C} &=
\left(
  \begin{matrix}
    \Sigma_{A_C A_C}^{-1} & 0 & 0\\
 0 &   \Sigma_{C|A_C} & \Sigma_{CU|A_C} \\
 0 &  \Sigma_{CU|A_C}^{-1}  & \Sigma_{U|A_C} \\
   \end{matrix} 
             \right),
    \end{align*}
    $\Sigma_{C|A_C} = \Sigma_{CC} -
    \Sigma_{CA_C}\Sigma_{A_CA_C}^{-1}\Sigma_{A_CC}$,
    $\Sigma_{U|A_C} = \Sigma_{UU} -
    \Sigma_{UA_C}\Sigma_{A_C A_C}^{-1}\Sigma_{A_CU}$,
    and $\Sigma_{CU|A_C} = \Sigma_{CU}-\Sigma_{CA_C}\Sigma_{A_CA_C}^{-1}\Sigma_{A_CU}$.
    Further,
    \begin{itemize}
        \item[(i)] The conditional distribution of $\sqrt{n}\delta^{F_C}_{C}$ given
    $\{ (\delta^{F_C}_C > 0), (W_{A_C|F_{C}} \widetilde{\delta}_{A_C} \leq 0) \}$ 
    is asymptotically equivalent to the 
    conditional distribution of $\sqrt{n} \delta^{F_C}_{C}$ given
    $(\delta^{F_C}_{C} > 0)$, which   is a truncated normal 
        \begin{align*}
        tN_{|C|}(0, \Sigma_{C|A_{C}})
        \end{align*}
                \item[(ii)] The conditional distribution of $\sqrt{n} \delta^{F_C}_{U}$ given
    $\{ (\delta^{F_C}_C > 0), (W_{A_C|F_C} \widetilde{\delta}_{A_C} \leq 0) \}$ 
    asymptotically equivalent  to the 
    conditional distribution of $\sqrt{n} \delta^{F_C}_{U}$ given
    $(\delta^{F_C}_C > 0)$ and  
    is
    \begin{align*}
    \widetilde{SN}_{|C|, p-b}\big(0, \Sigma_{CU|A_C} \Sigma_{U|A_C}^{-1},
    \Sigma_{C|U,A_{C}} \big).
\end{align*}
where
$\Sigma_{C|U,A_{C}}=\Sigma_{C|A_C}- \Sigma_{CU|A_C}  \Sigma_{U|A_C}^{-1} \Sigma_{UC|A_C}$.
    \end{itemize}
\end{lemma}  
\begin{proof}
  To simplify notation we let $A := A_C= B\backslash C$ and
  $F :=F_C= C\cup U$.      
  By using Remark \ref{remark:SWinv} and $\Sigma=W^{-1}$, we have
  $\delta^{F}=\widetilde{\delta}_F+W_{FF}^{-1}W_{FA}\widetilde{\delta}_A$, that is,
\begin{align*}
  \left(
  \begin{matrix}
    \delta^{F}_{C} \\
       \delta^{F}_{U} \\
   \end{matrix} 
  \right) &=
   \left(
  \begin{matrix}
    \widetilde{\delta}_{C} \\
       \widetilde{\delta}_{U} \\
   \end{matrix} 
  \right) -  
  \left(
  \begin{matrix}
    \Sigma_{CA} \\
    \Sigma_{UA} \\
   \end{matrix} 
             \right)\Sigma_{AA}^{-1}\widetilde{\delta}_{A}.
\end{align*}  
Thus $\delta^{F}$ and $\widetilde{\delta}_{A}$
are asymptotically uncorrelated, which gives
\begin{align*}
V_{C} = \left(
  \begin{matrix}
    \Sigma_{A A}^{-1} & 0 & 0\\
 0 &   \Sigma_{C|A} & \Sigma_{CU} - \Sigma_{CA}\Sigma_{AA}^{-1}\Sigma_{AU} \\
 0 &  \Sigma_{UC} - \Sigma_{UA}\Sigma_{AA}^{-1}\Sigma_{AC}  & \Sigma_{U|A} \\
   \end{matrix} 
             \right),
\end{align*}
and proves (\ref{eq:Waa}).

To prove (i) and (ii) we note that by (\ref{eq:Waa})
$\widetilde{\delta}_{A_C}$ and $\delta^{F_C}$ are asymptotically independent. 
Thus, using $\sqrt{n}\delta^{F_C}_{C}\Dcon N(0,
\Sigma_{C|A})$, the conditional distribution of
$\sqrt{n}\delta^{F_C}_{C}$ given
    $( \delta^{F_C}_C > 0)$ is the truncated normal
    \begin{align*}
  \frac{\phi_{|C|}(x; \Sigma_{C|A})}{\operatorname{pr}(\delta^{F_C}_C>0)} I(x>0)
    \end{align*}
This proves (i). 
To prove (ii), note that 
\begin{align*}
    \sqrt{n} 
        \left(\begin{matrix}
            \delta^{F}_C \\
            \delta^{F}_U\\
        \end{matrix} \right)
        \Dcon N_{|F|} \left(
        \left( \begin{matrix}
            0 \\ 0
        \end{matrix} \right)
        , \left( \begin{matrix}
            \Sigma_{C|A} & \Sigma_{CU|A} \\
            \Sigma_{CU|A}^\top   & \Sigma_{U|A}
        \end{matrix}\right) \right),
\end{align*}
where $\Sigma_{CU|A_C} = \Sigma_{CU}-\Sigma_{CA_C}\Sigma_{A_CA_C}^{-1}\Sigma_{A_CU}$.
Then the conditional distribution of $\sqrt{n}\delta^{F}_U$ given $(\delta^{F}_C > 0)$ is 
\begin{align*}
    \frac{\operatorname{pr}(\delta^{F}_C > 0 \mid \sqrt{n} \delta^{F}_U= x)}{\operatorname{pr}(\delta^{F}_C > 0) }\phi_{p-b}(x;\Sigma_{U|A}),
\end{align*}
where $\sqrt{n} \delta^{F}_C \mid \sqrt{n} \delta^{F}_U= x \Dcon N_{|C|}(\Sigma_{CU|A}\Sigma_{U|A}^{-1}x,
\Sigma_{C|U,A}) $ and
$\Sigma_{C|U,A}=\Sigma_{C|A}- \Sigma_{CU|A}  \Sigma_{U|A}^{-1} \Sigma_{UC|A}$.
From Definition~\ref{def:SN_Variants}, it is
\begin{align*}
   \widetilde{SN}_{|C|, p-b}\big(0, \Sigma_{CU|A} \Sigma_{U|A}^{-1},
    \Sigma_{U|A}, \Sigma_{C|U,A} \big)
\end{align*}
where we use in the last line that
$\Sigma_{C|U,A}=\Sigma_{C|A}- \Sigma_{CU|A}  \Sigma_{U|A}^{-1} \Sigma_{UC|A}$.
Alternatively, this can be proved directly from Lemma~\ref{lemma:SN_TSN} (ii).
\end{proof}

\vspace{1mm}
\noindent
\emph{Proof of Corollary \ref{theorem:CLT3}} 
By using Corollary~\ref{cor:blocks} and Lemma~\ref{cor:WSigmaInve} we have
\begin{align*}
\sqrt{n}\widehat{\delta} \Dcon 
    \sum_{C\in \mathcal{P}(B)}
    \left(
    \begin{matrix}
      m_{A_{C}}(0)\\
        tN_{|C|}(0, \Sigma_{C|A_C}) \\
        \widetilde{SN}_{|C|, p-b}\big(0, \Sigma_{CU|A_C} \Sigma_{U|A_C}^{-1},
    \Sigma_{U|A_C}, \Sigma_{C|U,A_{C}} \big) \\
      \end{matrix}
    \right)       \operatorname{pr}\big( \delta^{F_C}_C > 0\big)
    \operatorname{pr}\big( \Sigma_{A_C A_C}^{-1}\widetilde{\delta}_{A_C} \leq 0 \big),  \end{align*}
The above together with (\ref{eq:Lest}) and that $\xi_{B}=0$ gives the result.
\hfill $\Box$

\vspace{3mm}

\begin{lemma}\label{lemma:finite_approx}
Suppose (\ref{eq:rate2}) and 
Assumption \ref{assum:CLT} hold. 
Further, 
suppose  
$\xi_{B} = -\theta_{0,B}\neq 0$. 
We define the disjoint sets 
$B_{1}$ and $B_{2}$ such that $B=B_{1}\cup B_{2}$, 
$\xi_{B_{1}}=0$ and 
$\xi_{j} <0$ for all $j\in B_{2}$, hence  $\xi = (0_{B_{1}},\xi_{B_{2}},0_{U})$.
Let 
$U_{1} = U\cup B_{2}$. Then Theorem \ref{theorem:CLT}
holds with $B_{1}$ replacing $B$ and $U_{1}$ replacing $U$. 
However, for small samples 
and $\xi_{B}$ close to the boundary 
\begin{align}
\label{eq:finite}
   \widehat{\delta} &\approxD 
   \sum_{\stackrel{C_{1}\in \mathcal{P}(B_{1}),
C_{2}\in \mathcal{P}(B_{2})}{C =C_{1}\times C_{2}}} 
f_{C}(\widetilde{\delta}) \operatorname{pr}\big( \delta^{F_C}_{C} > \xi_{C},W_{A_{C}|F_{C}} (\widetilde{\delta}_{A_C}-\xi_{A_C}) \leq 0 \big)
\end{align} 
where 
\begin{align*}
f_{C}(\widetilde{\delta}) = 
 \left(
    \begin{matrix}
      m_{A_{C}}(\xi_{A_C})\\
        tSN_{|A_C|,|C|}(\zeta_{A_C},\zeta_{C}, V_{A_C C}V_{CC}^{-1}, V_{C C}/n,
        V_{A_C|C}/n)  \\
     SN_{b, p-b}\left(\zeta_{A_C},\zeta_{C},
     V_{A_CU}V_{UU}^{-1}, V_{CU}V_{UU}^{-1}, V_{UU}/n, V_{B|U}/n \right)  
      \end{matrix}
  \right) -
  \left(
  \begin{matrix} 
  0 \\
  \zeta_C-\xi_C\\
  \zeta_U\\
  \end{matrix}
  \right)
\end{align*}
where 
$\zeta_{A_{C}} = (T_{C}\xi)_{A_C},\quad 
\zeta_{C} = (T_{C}\xi)_{C}, \quad
\zeta_U = (T_C \xi)_U$, 
and $m_{A_C}(\xi_{A})$ is the point mass centered at $\xi_{A}$.
\end{lemma}
\begin{proof}
We recall that $B = B_1 \cup B_2$, where $\xi_{B_1} = 0$. Also 
we denote $C \in \mathcal{P}(B)$ as $C=C_1 \times C_2$ where $C_1 \in \mathcal{P}(B_1)$ and $C_2 \in \mathcal{P}(B_2)$, and
the event $R_C = \{ \delta^{F_C}_C > \xi_C, W_{A_C|F_C}(\widetilde{\delta}_{A_C}-\xi_{A_C})\leq 0\}$.
We apply 
Theorem~\ref{theorem:blocks} to have 
\begin{align} \label{eq:delta_hat_E2}
    \widehat{\delta} 
    &= \sum_{C \in \mathcal{P}(B)} \left( \begin{matrix}
        \xi_{A_C} \\
        \delta^{F_C}
    \end{matrix} \right) I_{R_C}
    = \sum_{\stackrel{C_{1}\in \mathcal{P}(B_{1}),
C_{2}\in \mathcal{P}(B_{2})}{C =C_{1}\times C_{2}}} 
\left( \begin{matrix}
        \xi_{A_C} \\
        \delta^{F_C}
    \end{matrix} \right) I_{R_C},
\end{align}
where $\delta^{F_C} = \argmin_{\delta_{F_C} \in \R^{|F_C|}}Q(\xi_{A_C}, \delta_{F_C})$.

In order to approximate the finite sample distribution of $\widehat{\delta}$, 
we use the transformation matrix $T_C$ defined in Theorem~\ref{theorem:CLT}. 
Under $\widetilde{\delta} \approxD N(0, \Sigma/n)$,
we have
\begin{align*}
    \Delta^{(C)} 
    := T_C \widetilde{\delta}
    = \left( \begin{matrix}
        \Delta_{A_C} \\
        \Delta_{C}\\
        \Delta_{U}
    \end{matrix} \right) \approxD N(0, V_C/n)
\end{align*}
where $V_C = T_C \Sigma T_C^\top $ and by the definition of $T_C$
\begin{align*}
    \Delta_{A_C} = W_{A_C|F_C} \widetilde{\delta}_{A_C},
    \quad
    \Delta_C = W_{C|U}^{-1}W_{CA_C|U}\widetilde{\delta}_{A_C}+\widetilde{\delta}_C,
    \quad
    \Delta_U = W_{U|C}^{-1}W_{UA_C|C}\widetilde{\delta}_{A_C}+\widetilde{\delta}_U.
\end{align*}
From Remark \ref{remark:SWinv} and the decomposition as in \eqref{eq:deltaF_CU}, 
we have
\begin{align*}
    \delta^{F_C}_C = \widetilde{\delta}_C + W_{C|U}^{-1}W_{CA_C|U} (\widetilde{\delta}_{A_C}-\xi_{A_C})
    \quad \textrm{and} \quad
    \delta^{F_C}_U = \widetilde{\delta}_U + W_{U|C}^{-1}W_{UA_C|C} (\widetilde{\delta}_{A_C}-\xi_{A_C}),
\end{align*}
hence we can express
$\Delta_C = \delta^{F_C}_C + (\zeta_C -\xi_C)$ and $\Delta_U = \delta^{F_C}_U + \zeta_U$
(where $\zeta=T_C\xi$).
Moreover, the event $R_C$ can be written as
$R_C = \{\Delta_C > \zeta_C, \Delta_{A_C} < \zeta_{A_C} \}$.
Hence \eqref{eq:delta_hat_E2} can be expressed as
\begin{align*}
        \widehat{\delta} = \sum_{\stackrel{C_{1}\in \mathcal{P}(B_{1}),
C_{2}\in \mathcal{P}(B_{2})}{C =C_{1}\times C_{2}}}  \left( \begin{matrix}
        \xi_{A_C} \\
        \Delta_C - (\zeta_C-\xi_C) \\
        \Delta_U - \zeta_U
    \end{matrix}\right) I[\Delta_C > \zeta_C, \Delta_{A_C} < \zeta_{A_C}].
\end{align*}
Therefore, under $\Delta^{(C)} \approxD N(0, V_C/n)$, we are interested in two quantities:
(a) the conditional distribution of $(\Delta_C -(\zeta_C- \xi_C) )$ given $R_C$
and (b) the conditional distribution of $\Delta_U - \zeta_U$ given $R_C$ to approximate the finite sample distribution of $\widehat{\delta}$.
To derive (a), Lemma~\ref{lemma:SN_TSN}(i) gives
\begin{align*}
    \Delta_C \mid (\Delta_C > \zeta_C, \Delta_{A_C} < \zeta_{A_C}) 
    \approxD
    tSN_{|A_C|,|C|}(\zeta_{A_C}, \zeta_C, V_{A_CC}V_{CC}^{-1}, V_{CC}/n, V_{A_C|C}/n).
\end{align*}
Since $(\zeta_C- \xi_C)$ is deterministic, we can express (a) into
\begin{align*}
    &(\Delta_C -(\zeta_C- \xi_C))  \mid   (\Delta_C > \zeta_C, \Delta_{A_C} < \zeta_{A_C}) 
    \approxD 
    tSN_{|A_C|,|C|}(\zeta_{A_C}, \zeta_C, V_{A_CC}V_{CC}^{-1}, V_{CC}/n, V_{A_C|C}/n)
    -(\zeta_C- \xi_C).
\end{align*}
Using the same argument and Lemma~\ref{lemma:SN_TSN}(ii), we can derive (b)
\begin{align*}
    (\Delta_U - \zeta_U) \mid (\Delta_C > \zeta_C, \Delta_{A_C} < \zeta_{A_C})
    \approxD
    SN_{b, p-b}\left(\zeta_{A_C},\zeta_{C},
     V_{A_CU}V_{UU}^{-1}, V_{CU}V_{UU}^{-1}, V_{UU}/n, V_{B|U}/n \right) 
     - \zeta_U.
\end{align*}
Combining these yields
\begin{align*}
f_{C}(\widetilde{\delta}) = 
 \left(
    \begin{matrix}
      m_{A_{C}}(\xi_{A_C})\\
        tSN_{|A_C|,|C|}(\zeta_{A_C},\zeta_{C}, V_{A_C C}V_{CC}^{-1}, V_{C C}/n,
        V_{A_C|C}/n)  \\
     SN_{b, p-b}\left(\zeta_{A_C},\zeta_{C},
     V_{A_CU}V_{UU}^{-1}, V_{CU}V_{UU}^{-1}, V_{UU}/n, V_{B|U}/n \right)  
      \end{matrix}
  \right) -
  \left(
  \begin{matrix} 
  0 \\
  \zeta_C-\xi_C\\
  \zeta_U\\
  \end{matrix}
  \right),
\end{align*}
which completes the proof.
\end{proof}

Finally, we show how the above lemma leads to 
the 
Theorem \ref{theorem:CLT} with $C_{1}$ replacing $C$.
If $C\in \mathcal{P}(B_{1})\times (\mathcal{P}(B_{2})\backslash B_{2})$, then $\xi_{A_C}\neq 0$ (hence it has to contain some strictly negative entries). 
Since
$W_{A_C|F_C}$ is positive definite and 
$\widetilde{\delta}_{A_C}\Pcon 0$ we can apply Lemma \ref{lemma:convergence} to give
\begin{align}
\label{eq:probzero}
\operatorname{pr}(\delta^{F_C}_{C} > \xi_{C},W_{A_C|F_C}(\widetilde{\delta}_{A_C}-\xi_{A_C})\leq 0)\leq 
\operatorname{pr}(W_{A_C|F_C}(\widetilde{\delta}_{A_C}-\xi_{A_C})\leq 0)\rightarrow 0.
\end{align}
 Thus 
the terms in (\ref{eq:finite}) with
$C\in \mathcal{P}(B_{1})\times (\mathcal{P}(B_{2})\backslash B_{2})$ converge to zero,
which leads to the large sample approximation given in 
Theorem \ref{theorem:CLT} with $C_{1}$ replacing $C$.

\section{Proof of the results in Section \ref{sec:LLRT}}\label{sup:LLRT}

We first prove the quadratic approximation
of LR given in (\ref{eq:LRclassicalquad}).

Using $\mathcal{L}_{n}$ and $Q_{n}$ we define the minimisers
\begin{align}
\label{eq:L_est_H1}
    \widehat{\theta}_{1,L} = \argmin_{\theta \in \mathcal{H}_{1}} \mathcal{L}_n(\theta)
\quad\textrm{ and }\quad
    \widehat{\theta}_{0,L} = \argmin_{\theta \in
  \mathcal{H}_{0}}\mathcal{L}_{n}(\theta).
\end{align}
and
\begin{align}
\label{eq:Q_est_H1}
\widehat{\delta}_{1,Q} = \argmin_{\delta \in \mathcal{H}_{1}-\theta_0}
  Q_n(\delta) \quad\textrm{ and }\quad  
\widehat{\delta}_{0,Q} = \argmin_{\delta \in \mathcal{H}_{0}-\theta_0}Q_n(\delta),
\end{align}
where $\mathcal{H}_{j}-\theta_0$ denotes the parameter space 
$\mathcal{H}_{j}$ shifted by $\theta_0$.

\begin{lemma}\label{lemma:QuadLLRT}
  Suppose the null is true with population parameter $\theta_0$, 
  Assumption \ref{assum:CLT}(a)  and
  condition (\ref{eq:rate2}) holds for 
   $\widehat{\theta}_{j,L}$
  and 
  $\widehat{\theta}_{j,Q}$ in the case of
  $j=0$ and $1$
  (where $\widehat{\theta}_{j,Q}= \widehat{\delta}_{j,Q}+\theta_0$). 
Then
\begin{align*}
     LR   = 2n[Q_n(\widehat{\delta}_{0,Q}) - Q_n(\widehat{\delta}_{1,Q})] + o_p(1).
    \end{align*}
 where $Q_{n}(\delta) = \mathcal{L}_{n}(\theta_0) + s^{\top}\delta+\delta^{\top}W\delta/2$.   
\end{lemma}
\begin{proof}
In the proof we will
let $\widehat{\delta}_{j,L} = 
\widehat{\theta}_{j,L}-\theta_{0}$. Since
$\widehat{\delta}_{j,Q}$ minimises 
$Q_{n}(\cdot)$ over $\mathcal{C}_{j}$,  then
\begin{align}
 0\leq   Q_{n}(\widehat{\delta}_{j,L})-Q_{n}(\widehat{\delta}_{j,Q})
    &= Q_{n}(\widehat{\delta}_{j,L}) - \mathcal{L}_n(\widehat{\theta}_{j,L}) 
    + \underbrace{\mathcal{L}_n(\widehat{\theta}_{j,L}) - \mathcal{L}_n(\widehat{\theta}_{j,Q})}_{\leq 0}
    + \mathcal{L}_n(\widehat{\theta}_{j,Q}) -Q_{n}(\widehat{\delta}_{j,Q}) \nonumber\\
 &\leq Q_{n}(\widehat{\delta}_{j,L}) - \mathcal{L}_n(\widehat{\theta}_{j,L}) +
    \mathcal{L}_n(\widehat{\theta}_{j,Q}) -Q_{n}(\widehat{\delta}_{j,Q})  \nonumber\\
    &\leq |Q_{n}(\widehat{\delta}_{j,L}) - \mathcal{L}_n(\widehat{\theta}_{j,L}) | + |\mathcal{L}_n(\widehat{\theta}_{j,Q}) -Q_{n}(\widehat{\delta}_{j,Q})| \nonumber\\
    & = o_p(n^{-1}) \label{eq:rate5}
\end{align}

Now we return to the proof.
Under the stated assumptions of the lemma 
(specifically (\ref{eq:rate2}) and the null holding with parameter $\theta_0$) we have 
\begin{align}
\label{eq:RATE}
\sup_{\theta\in \mathcal{I}_{n}}|\mathcal{L}_{n}(\theta)- Q_{n}(\delta) | = o_{p}(n^{-1}) \textrm{ and }
\|\widehat{\theta}_{j,L}-\widehat{\theta}_{j,Q}\|= o_{p}(n^{-1/2}),
\end{align}
where $\mathcal{I}_{n} = \{\theta:\|\theta-\theta_0\| = O(n^{-1/2})\}$ and
$\delta=\theta-\theta_0$.
Further, by using Lemma \ref{lemma:LSrateG} we have 
$\|\widehat{\theta}_{j,Q}-\theta_0\|_2=
O_p(n^{-1/2})$. 
Thus by the triangle inequality we have
\begin{align*}
\|\widehat{\theta}_{j,L}-\theta_0\|_2
\leq \|\widehat{\theta}_{j,Q}-\theta_0\|_2+\|\widehat{\theta}_{j,Q}-\widehat{\theta}_{j,L}\|_2 = O_{p}(n^{-1/2}).
\end{align*}
The above allows us to  apply (\ref{eq:RATE}) to give 
\begin{align}
\mathcal{L}_{n}(\widehat{\theta}_{j,L}) =
Q_{n}(\widehat{\delta}_{j,L}) + o_{p}(n^{-1}) =
Q_{n}(\widehat{\delta}_{j,Q}) + o_{p}(n^{-1})
\end{align}
where the last line in the above follows from 
(\ref{eq:rate5}).

As the above holds for $j=0$ and $j=1$ we have 
\begin{align*}
LR &= 2n\left(\mathcal{L}_{n}(\widehat{\theta}_{0,L})-\mathcal{L}_{n}(\widehat{\theta}_{1,L})\right)
=2n\left(Q_{n}(\widehat{\delta}_{0,L})- Q_{n}(\widehat{\delta}_{1,L})\right) + o_{p}(1) \\
&=2n\left(Q_{n}(\widehat{\delta}_{0,Q})-Q_{n}(\widehat{\delta}_{1,Q})\right) + o_{p}(1)
\end{align*}
thus proving the desired result. 
\end{proof}

The lemma states that if the null is true and the estimators of the null parameters satisfy (\ref{eq:rate2}) under both the 
null and alternative then we can
replace $LR$ with its quadratic approximation.

We will often move between $Q_n(\delta) = \mathcal{L}_n(\theta_0)+s^{\top}\delta + \frac{1}{2}\delta^{\top}W\delta$ and 
\begin{align*}
Q_n(\delta) = \mathcal{L}_n(\theta_0)+\frac{1}{2}(\delta-\widetilde{\delta})^{\top}
W(\delta-\widetilde{\delta})-\frac{1}{2}\widetilde{\delta}^{\top}W\widetilde{\delta}.
\end{align*}
For example, using this representation we have
  \begin{align*}
LR &= 2n(Q_{n}(\widehat{\delta}_{0,Q})-Q_{n}(\widehat{\delta}_{1,Q})) + o_{p}(1) \\
&=n(\widehat{\delta}_{0,L}-\widetilde{\delta})^{\top}W
(\widehat{\delta}_{0,L}-\widetilde{\delta})-
n(\widehat{\delta}_{1,L}-\widetilde{\delta})^{\top}W
(\widehat{\delta}_{1,L}-\widetilde{\delta})
+o_{p}(1),
\end{align*}
which connects the classical quadratic expansion with Lemma \ref{lemma:QuadLLRT}.

\subsection{Some preliminary results}

In the following lemma we apply Theorem \ref{theorem:blocks} to obtain
an explicit expression for the quadratic loss $Q_{n}(\cdot)$. We show
that it can be written as the sum of
quadratic forms over disjoint sets. It is
this results that will facilitate the asymptotic analysis of $LR$.

\begin{lemma}\label{lemma:Quasi-Like}
  Let $\widehat{\delta} = 
  \arg\min_{\delta\in \mathcal{C}_{\varepsilon}}Q_{n}(\delta)$, then
\begin{align}
\label{eq:Qexpansion}
&  Q_{n}(\widehat{\delta}) - \mathcal{L}_n(\theta_0) \nonumber\\ 
=&  \sum_{C\in \mathcal{P}(B)}
 \left(
 s_{A_C}^\top \xi_{A_C} + \frac{1}{2} \xi_{A_C}^\top W_{A_CA_C} 
 \xi_{A_C} - \frac{1}{2} (s_{F_C} + W_{F_CA_C} \xi_A)^\top W_{F_CF_C}^{-1} (s_{F_C} + W_{F_CA_C} \xi_{A_C})
 \right) \nonumber\\
&\qquad \times I\left[\delta^{F_C}_{C} > \xi_{C} \right]
        I\left[
                             W_{A_{C}|F_{C}}(\widetilde{\delta}_{A_{C}}-\xi_{A_C}) \leq 0\right] \nonumber\\
=& \frac{1}{2} \sum_{C \in \mathcal{P}(B)} (\xi_{A_C}-\widetilde{\delta}_{A_C})^\top W_{A_C|F_C} (\xi_{A_C}-\widetilde{\delta}_{A_C})
    I[\delta^{F_C}_C>\xi_C] I[W_{A_C|F_C} (\widetilde{\delta}_{A_C}-\xi_{A_C}) \leq 0] - 
    \frac{1}{2}\widetilde{\delta}^{\top}W\widetilde{\delta}
\end{align}
where $\widetilde{\delta} = -W^{-1}s$ (the latter representation can also be useful). 

In the special case that $\xi=0$ and 
$\widehat{\delta} = 
  \argmin_{\delta\in \mathcal{C}} Q_{n}(\delta)$
the above reduces to 
\begin{align*}
  Q_{n}(\widehat{\delta}) &= \mathcal{L}_n(\theta_0)
  -\frac{1}{2}\sum_{C\in \mathcal{P}(B)}s^{\top}_{F_C}W_{F_CF_C}^{-1}s_{F_C}I\left[\delta^{F_C}_{C} > 0 \right] I\left[
                             W_{A_{C}|F_{C}} \widetilde{\delta}_{A_{C}} \leq 0\right] \\
 &= \mathcal{L}_n(\theta_0)
  -\frac{1}{2}\sum_{C\in \mathcal{P}(B)}(\delta^{F_C})^{\top} W_{F_CF_C}\delta^{F_C}I\left[\delta^{F_C}_{C} > 0 \right] I\left[
                             W_{A_{C}|F_{C}} \widetilde{\delta}_{A_{C}} \leq 0\right]                            
\end{align*}
\end{lemma}  

\begin{proof}
From Theorem~\ref{theorem:blocks}  we have 
\begin{align*}
    \widehat{\delta}
     &= \sum_{C\in \mathcal{P}(B)}
\underbrace{\left(
\begin{array}{c}
\xi_{A_C} \\
-W_{F_CF_C}^{-1}(s_{F_C}+W_{F_CA_C}\xi_{A_C})
\end{array}
     \right)}_{=\delta^{(C)}}
I\left[\delta^{F_C}_{C} > \xi_{C} \right]
        I\left[
        W_{A_{C}|F_{C}}(\widetilde{\delta}_{A_{C}}-\xi_{A_C}) \leq 0\right],
\end{align*}
where the events 
\begin{align*}
R_{C} = (\delta^{F_C}_{C} > \xi_{C},W_{A_{C}|F_{C}} (\widetilde{\delta}_{A_{C}}-\xi_{A_C})
\leq 0)
\end{align*}
form a disjoint partition of $\R^{p}$. 
Since 
\begin{align*}
  Q_n(\widehat{\delta}) &= \mathcal{L}_n(\theta_0) + s^{\top}\widehat{\delta} +
                          \frac{1}{2}\widehat{\delta}^{\top}W \widehat{\delta}.
\end{align*}
Substituting $\widehat{\delta}$ into the above and using that 
$R_{C}$ are disjoint sets we have
\begin{align*}
& s^{\top}\widehat{\delta}+
   \frac{1}{2}\widehat{\delta}^{\top}W \widehat{\delta}
  =\sum_{C\in \mathcal{P}(B)}
 \left(
 s^{\top}\delta^{(C)}+\frac{1}{2}(\delta^{(C)})^{\top}W\delta^{(C)}
 \right)
I\left[\delta^{F_C}_{C} > \xi_{C} \right]
        I\left[
                             W_{A_{C}|F_{C}}(\widetilde{\delta}_{A_{C}}-\xi_{A_C}) \leq 0\right]
\end{align*}
Now we simplify the term inside the summand (omitting the dependencies of $A$ and $F$ on $C$ to simplify notation) to give 
\begin{align*}
s^{\top}\delta^{(C)}+\frac{1}{2}(\delta^{(C)})^{\top}W\delta^{(C)}=
s_A^\top \xi_A + \frac{1}{2} \xi_A^\top W_{AA} \xi_A + (s_F + W_{FA} \xi_A)^\top \delta_F + \frac{1}{2} \delta_F^\top W_{FF} \delta_F.
\end{align*}
Since $\delta_{F} = -W_{FF}^{-1}(s_{F}+W_{FA}\xi_{A})$ we let 
$u  = s_{F}+W_{FA}\xi_{A}$. Using this notation the last two terms in the above expression can be written as 
\begin{align*}
u^\top \delta_F + \frac{1}{2} \delta_F^\top W_{FF} \delta_F &= u^\top (-W_{FF}^{-1} u) + \frac{1}{2} (-W_{FF}^{-1} u)^\top W_{FF} (-W_{FF}^{-1} u) \\
&= -u^\top W_{FF}^{-1} u + \frac{1}{2} u^\top W_{FF}^{-1} u \\
&= -\frac{1}{2} u^\top W_{FF}^{-1} u.
\end{align*}
Substituting $u = s_F + W_{FA} \xi_A$ back into the expression gives
\begin{align*}
s^\top \delta + \frac{1}{2} \delta^\top W \delta = s_A^\top \xi_A + \frac{1}{2} \xi_A^\top W_{AA} \xi_A - \frac{1}{2} (s_F + W_{FA} \xi_A)^\top W_{FF}^{-1} (s_F + W_{FA} \xi_A).
\end{align*}
Finally, placing the above back into the summand gives 
\begin{align*}
  Q_n(\widehat{\delta}) &= \mathcal{L}_n(\theta_0) + 
  \sum_{C\in \mathcal{P}(B)}
 \left(
 s_{A_C}^\top \xi_{A_C} + \frac{1}{2} \xi_{A_C}^\top W_{A_CA_C} 
 \xi_{A_C} - \frac{1}{2} (s_{F_C} + W_{F_CA_C} \xi_A)^\top W_{F_CF_C}^{-1} (s_{F_C} + W_{F_CA_C} \xi_{A_C})
 \right) \\
&\qquad \times I\left[\delta^{F_C}_{C} > \xi_{C} \right]
        I\left[
                             W_{A_{C}|F_{C}}(\widetilde{\delta}_{A_{C}}-\xi_{A_C}) \leq 0\right].
\end{align*}
To prove the second part of (\ref{eq:Qexpansion})
we recall that
$Q_n(\widehat{\delta}) = \mathcal{L}_n(\theta_0)+\frac{1}{2}(\widehat{\delta}-\widetilde{\delta})^{\top}
W(\widehat{\delta}-\widetilde{\delta}) - \frac{1}{2}\widetilde{\delta}^{\top}W\widetilde{\delta}$.
Thus we obtain an expression for the first term in this quadratic form. For this 
we use (\ref{eq:deltaproject}):
\begin{align*}
\delta^{F} =
\widetilde{\delta}_{F}+(W_{FF})^{-1} W_{FA}\,(\widetilde{\delta}_{A} -
\xi_{A}).
\end{align*}
Thus by using the above and Theorem \ref{theorem:blocks} we have 
\begin{align*}
    \widehat{\delta} - \widetilde{\delta}
    = \sum_{C \in \mathcal{P}(B)} \left( \begin{matrix}
        \xi_{A_C} - \widetilde{\delta}_{A_C} \\
        -W_{F_C F_C}^{-1}W_{F_C A_C}(\xi_{A_C} - \widetilde{\delta}_{A_C})
    \end{matrix}\right)
    I[\delta^{F_C}_C>\xi_C] I[W_{A_C|F_C} (\widetilde{\delta}_{A_C}-\xi_{A_C}) \leq 0]. 
\end{align*}
Plugging the above into 
$\|\widehat{\delta}-\widetilde{\delta}\|_{W}$ gives 
\begin{align*}
  \|\widehat{\delta}-\widetilde{\delta}\|_{W}^{2}  
    =  \sum_{C \in \mathcal{P}(B)} (\xi_{A_C}-\widetilde{\delta}_{A_C})^\top W_{A_C|F_C} (\xi_{A_C}-\widetilde{\delta}_{A_C})
    I[\delta^{F_C}_C>\xi_C] I[W_{A_C|F_C} (\widetilde{\delta}_{A_C}-\xi_{A_C}) \leq 0]. 
\end{align*}
This proves the second part of (\ref{eq:Qexpansion}).

In the special case that $\xi=0$ the above reduces to 
\begin{align*}
  Q_n(\widehat{\delta}) &= \mathcal{L}_n(\theta_0)
  -\frac{1}{2}\sum_{C\in \mathcal{P}(B)}s^{\top}_{F_C}W_{F_CF_C}^{-1}s_{F_C}I\left[\delta^{F_C}_{C} > 0 \right] I\left[
                             W_{A_{C}|F_{C}} \widetilde{\delta}_{A_{C}} \leq 0\right] \\
&= \mathcal{L}_n(\theta_0)
  -\frac{1}{2}\sum_{C\in \mathcal{P}(B)}(\delta^{F_C})^{\top}W_{F_CF_C}\delta^{F_C}I\left[\delta^{F_C}_{C} > 0 \right] I\left[
                             W_{A_{C}|F_{C}} \widetilde{\delta}_{A_{C}} \leq 0\right] 
\end{align*}
where in the last line we use that 
$\delta^{F} = -W_{FF}^{-1}s_{F}$, 
which gives the desired result. 
\end{proof}

In the next lemma we obtain a simple identity for the difference between two quadratic forms. 
\begin{lemma}\label{lemma:quadratic}
Suppose $\delta = (\delta_1^{\top},\delta_2^{\top})^{\top}$
($\delta_{1}$ and $\delta_{2}$ are $q_1$ and $q_2$-dimensional vectors),
$s = (s_1^{\top}, s_2^{\top})^{\top}$ and 
\begin{align*}
W = 
\left(
\begin{array}{cc}
W_{11} & W_{12} \\
W_{12}^{\top} & W_{22} \\
\end{array}
\right)
\end{align*}
where $W$ is a strictly positive definite, symmetric matrix.
Let $\widetilde{\delta} = (\widetilde{\delta}_1^{\top},\widetilde{\delta}_2^{\top})^{\top}
= -W^{-1}s$ and $\delta^{(2)} = (0,\delta^{2})
= (0, -W_{22}^{-1}s_2)$.
Then 
\begin{align*}
\frac{1}{2}\left(\widetilde{\delta}^{\top}W\widetilde{\delta} - 
(\delta^{2})^{\top}W_{22}\delta^2\right)
=\frac{1}{2}\widetilde{\delta}^{\top}_1
W_{1|2}\widetilde{\delta}_1
\end{align*}
where $W_{1|2} = (W_{11}-W_{12}W_{22}^{-1}W_{12}^{\top})$.
\end{lemma}
\begin{proof}
Expanding $W$ in terms of block matrices gives 
\begin{align}
\label{eq:L1}
\frac{1}{2}\widetilde{\delta}^{\top}W
\widetilde{\delta} = \frac{1}{2}(\widetilde{\delta}_1^{\top}W_{11}
\widetilde{\delta}_1+2
\widetilde{\delta}^{\top}_{1}W_{12}
\widetilde{\delta}_{2}+\widetilde{\delta}^{\top}_{2}W_{22}
\widetilde{\delta}_2).
\end{align}
Next, we consider 
\begin{align*}
-\frac{1}{2} (\delta^{2})^{\top}W_{22}
\delta^{2},
\end{align*}
where we recall that $\delta^{2} = -W_{22}^{-1}s_{2}$. We now rewrite
$\delta^{2}$ in terms of $\widetilde{\delta}$ (to align it with 
$\widetilde{\delta}^{\top}W
\widetilde{\delta}$). By using $s = -W \widetilde{\delta}$. This gives 
$\delta^2=\widetilde{\delta}_2 + W_{22}^{-1}W_{21}\widetilde{\delta}_{1}$ and 
\begin{align}
\label{eq:L2}
-\frac{1}{2} (\delta^{2})^{\top}W_{22}
\delta^{2}=
-\frac{1}{2}\left(
\widetilde{\delta}_{2}^{\top}W_{22}^{\top}
\widetilde{\delta}_{2} + 2
\widetilde{\delta}_{1}^{\top}W_{12}\widetilde{\delta}_{2} + \widetilde{\delta}_{1}^{\top}W_{12}W_{22}^{-1}W_{12}^{\top}
\widetilde{\delta}_{1}\right)
\end{align}
Taking differences between (\ref{eq:L1}) and (\ref{eq:L2})
and using that 
$W_{1|2} = W_{11}-W_{12}^{\top}W_{22}W_{21}$
gives 
\begin{align*}
\frac{1}{2}\widetilde{\delta}^{\top}W
\widetilde{\delta}  - \frac{1}{2} (\delta^{2})^{\top}W_{22}
\delta^{2} = 
(\widetilde{\delta}_{1})^{\top}W_{1|2}
\widetilde{\delta}_{1}
\end{align*}
and the desired result. 
\end{proof}

Finally, in the following lemma we show that 
$Q_n(\delta_{\xi})$ monotonically decreasing in 
$\xi$ where $\delta_{\xi} = \argmin_{\delta\in \mathcal{C}_{\xi}}Q_n(\delta)$. 
\begin{lemma}\label{lemma:monotonic}
Suppose $W$ is a positive definite matrix and 
$\xi\leq 0$. If $\xi_{1}<\xi_{2}\leq 0$, then
\begin{align*}
\min_{\delta \in C_{\xi_{1}}}Q_n(\delta) \leq 
\min_{\delta \in C_{\xi_{2}}}Q_n(\delta).
\end{align*}
\end{lemma}
\begin{proof} The proof is based on the fact that
\begin{align*}
Q_n(\delta) = \mathcal{L}_n(\theta_0)+ \frac{1}{2}(\delta-\widetilde{\delta})^{\top}W
(\delta-\widetilde{\delta}) -\frac{1}{2}\widetilde{\delta}^{\top}W\widetilde{\delta}
\end{align*}
where $\widetilde{\delta} = -W^{-1}s$ and 
$\frac{1}{2}\widetilde{\delta}^{\top}W\widetilde{\delta}$ is simply an intercept. 
Let 
\begin{align*}
\|\delta-\widetilde{\delta}\|_{W}^{2} = 
(\delta-\widetilde{\delta})^{\top}W
(\delta-\widetilde{\delta}),
\end{align*}
thus 
\begin{align*}
Q_n(\delta) = \mathcal{L}_n(\theta_0)+ \frac{1}{2}\|\delta-\widetilde{\delta}\|_{W}^{2}-
\frac{1}{2}\widetilde{\delta}^{\top}W\widetilde{\delta}.
\end{align*}
If $\xi_{1}<\xi_{2}$ then $\mathcal{C}_{\xi_2}\subset \mathcal{C}_{\xi_{1}}$. Consequently
\begin{align*}
\min_{\delta\in \mathcal{C}_{\xi_1}}\|\delta-\widetilde{\delta}\|_{W}\leq
\min_{\delta\in \mathcal{C}_{\xi_2}}\|\delta-\widetilde{\delta}\|_{W},
\end{align*}
which gives the result. 
\end{proof}

\subsection{Asymptotic approximation of LR}\label{sup:hyp}

Before getting into details of the the proof we first describe the main strategy that is common to most of the tests we consider. 
Suppose $J_{C}\subseteq [p]$,
$G\subseteq [p]$ and let
$\delta_{K} = - W_{K}^{-1}s_{K}$. We 
assume that 
$J_{C}\subseteq G$ for all $C\in \mathcal{P}(B)$ or $G\subseteq J_{C}$ for all $C\in \mathcal{P}(B)$.
The comment thread in all of the approximations that we derive below is that under certain types of 
null and alternative hypothesis one can apply
Lemma \ref{lemma:QuadLLRT}, to 
obtain the approximation
\begin{align*}
LR = 
\textrm{sgn}(H_0)
\left(n\sum_{C\in \mathcal{P}(B)}
\delta_{J_C}^{\top}W_{J_C}\delta_{J_C}I(s\in R_{C}) - n\delta_{G}^{\top}W_{G}\delta_{G}\right)+o_{p}(1)
\end{align*}
where $\textrm{sgn}(H_0)\in \{-1,1\}$
(depending on the form of the
null and alternative)
A key observation is that 
since $\{R_{C}\}$ is a partition of $\R^{p}$, then the above can be rewritten as 
\begin{align}
\label{eq:LRargument}
LR = \textrm{sgn}(H_0) n\sum_{C\in \mathcal{P}(B)}
\left(\delta_{J_C}^{\top}W_{J_C}\delta_{J_C} - \delta_{G}^{\top}W_{G}\delta_{G}\right)
I(s\in R_{C}) 
+o_{p}(1).
\end{align}
Here we can apply
Lemma \ref{lemma:quadratic}
\begin{align*}
\textrm{sgn}(H_0)[
\delta_{J_C}^{\top}W_{J_C}\delta_{J_C} - \delta_{G}^{\top}W_{G}\delta_{G}] = 
\left\{
\begin{array}{cc}
(\delta_{J_C})_{J_{C}-G}^{\top}
W_{(J_C-G)|G}(\delta_{J_C})_{J_C-G}
& G \subset J_{C} \\
(\delta_{G})_{G-J_{C}}^{\top}W_{(G-J_C)|J_C}(\delta_{G})_{G-J_C} & J_{C} \subset G \\
\end{array}
\right.
\end{align*}
We show below that this
reduction applies for 
Hypothesis 1-3 but 
not to Hypothesis 4.


\vspace{2mm}
\noindent {\bf Hypothesis 1}
\vspace{2mm}

\noindent \emph{Proof of Theorem \ref{theorem:LLRexpansion}}
To prove the result we will use the decomposition  $W_{FF}$ in terms of block matrices
\begin{align}
\label{eq:WFschur}  
W_{FF} = 
\left(
\begin{matrix}
W_{CC} & W_{CU} \\
W_{CU}^{\top} & W_{UU} \\
\end{matrix}
\right) \textrm{ and let }W_{C|U} = W_{CC} - W_{CU}W_{UU}^{-1}W_{UC}.
\end{align}
Returning to the proof, by using Lemma \ref{lemma:Quasi-Like} we have 
\begin{align}
  \label{eq:Qalt}
  Q_{n}(\widehat{\delta}) &= \mathcal{L}_{n}(\theta_0)-\frac{1}{2}\sum_{C\in \mathcal{P}(B)}
 (\delta^{F_C})^{\top}W_{F F}\delta^{F_C}
  I\left[(\delta^{F_{C}})_{C} > 0\right]
        I\left[W_{A_{C}|F_{C}}\widetilde{\delta}_{A_{C}} \leq 0\right].
                                \end{align}   
On the other hand, since $\{R_{C}\}_{C\in \mathcal{P}(C)}$ forms a
disjoint  partition of $\mathbb{R}^{p}$ we have 
\begin{align}
\label{eq:Qnull}
  Q_{n}(0_B,\delta^{U}) =&\mathcal{L}_{n}(\theta_0)-\frac{1}{2}
  (\delta^{U})^{\top}W_{UU}\delta^{U}\nonumber\\
  =& \mathcal{L}_{n}(\theta_0)
   -\frac{1}{2}\sum_{C \in \mathcal{P}(B)}
 (\delta^{U})^{\top}W_{UU}\delta^{U}
  I\left[(\delta^{F_{C}})_{C} > 0\right]
        I\left[W_{A_{C}|F_{C}}\widetilde{\delta}_{A_{C}} \leq 0\right].
\end{align}
Taking the difference between (\ref{eq:Qalt}) and 
(\ref{eq:Qnull}) and using  Lemma \ref{lemma:QuadLLRT} we have
\begin{align}
\label{eq:QAA}
  & LR_{1}
   =
 n \sum_{C\in \mathcal{P}(B)}
     \left\{(\delta^{F_{C}})^{\top}W_{F_CF_C}
     \delta^{F_C}-
      (\delta^{U})^{\top}W_{UU}\delta^{U}
     \right\}
  I\left[\delta^{F_{C}}_{C} > 0\right]
        I\left[W_{A_{C}|F_{C}}\widetilde{\delta}_{A_{C}} \leq 0\right]+o_{p}(1).
\end{align}
For all $C\in \mathcal{P}(B)\backslash \emptyset$ we have
\begin{align*}
 \left\{(\delta^{F_C})^{\top}W_{F_CF_C}\delta^{F_C}-
      (\delta^{U})^{\top}W_{UU}\delta^{U}
     \right\}  = 
(\delta^{F_{C}})_{C}^{\top}W_{C|U}
                  (\delta^{F_{C}})_{C}.
\end{align*}
where the display on the right hand side of the above follows from 
Lemma \ref{lemma:quadratic} and (\ref{eq:WFschur}).
Substituting the above into (\ref{eq:QAA}) gives the result. 
\hfill $\Box$

\vspace{2mm}

\noindent \emph{Proof of Corollary \ref{cor:pivotal}}
Using the Schur complement we have
    $W_{C|U} = \Sigma_{C|A_C}^{-1}$.
    Hence 
from Theorem~\ref{theorem:LLRexpansion}, the $LR$ in terms of $\Sigma$ is
\begin{align*}
    LR_1 = n \sum_{C \in \mathcal{P}(B)\backslash \emptyset} 
    (\delta^{F_C}_C)^\top \Sigma_{C|A_C}^{-1} \delta^{F_C}_C
    I[\delta^{F_C}_C > 0 ] I[ \Sigma_{A_C A_C}^{-1}\widetilde{\delta}_{A_C} \leq 0 ]
    + 0 \cdot I[\Sigma_{BB}^{-1} \widetilde{\delta}_B  \leq 0  ]
    +o_p(1)
\end{align*}

To derive the asymptotic distribution of $LR_1$, 
we recall that under Assumption \ref{assum:CLT}
we have 
$\sqrt{n} \widetilde{\delta} \Dcon N(0, \Sigma)$, thus 
  $\sqrt{n} \delta^{F_C} \Dcon N(0, \Sigma_{F_C|A_C})$, 
    implying $\sqrt{n} \delta^{F_C}_{C} \Dcon N_{|C|}(0, \Sigma_{C|A_C})$ by $(\Sigma_{F_C|A_C})_{CC} = \Sigma_{C|A_C}$. Thus 
\begin{align*}
Z = \sqrt{n} \Sigma_{C|A_C}^{-1/2} \delta^{F_C}_C \Dcon N_{|C|}(0, I_{|C|}).
\end{align*}
Using the above notation and that  
 $\Sigma_{C|A_C}^{1/2}$ is positive definite we can 
apply Lemma~\ref{lemma:trun_chi} to  give
\begin{align*}
&n(\delta^{F_C}_C)^\top \Sigma_{C|A_C}^{-1} \delta^{F_C}_C
    \mid (\delta^{F_C}_C > 0,\Sigma_{A_C A_C}^{-1}\widetilde{\delta}_{A_C} \leq 0) \\
 = &n(\delta^{F_C}_C)^\top \Sigma_{C|A_C}^{-1} \delta^{F_C}_C
    \mid (\delta^{F_C}_C > 0) \\   
=&    
    Z^\top Z \mid (\Sigma_{C|A_C}^{1/2} Z>0)
    \Dcon \chi^2_{|C|},
\end{align*}
where the second last line in the display uses that 
$\delta^{F_C}_C$ and $\widetilde{\delta}_{A_C}$ are asymptotically independent (see Corollary \ref{theorem:CLT3}).
This proves 
\begin{align*}
    LR_1 \Dcon  
    \sum_{C \in \mathcal{P}(B)\backslash \emptyset} 
   \chi^{2}_{|C|}
    \operatorname{pr}[\delta^{F_C}_C > 0 ] 
    \operatorname{pr}[ \Sigma_{A_C A_C}^{-1}\widetilde{\delta}_{A_C} \leq 0 ]
    + 0 \cdot \operatorname{pr}[\Sigma_{BB}^{-1} \widetilde{\delta}_B  \leq 0  ],
\end{align*}
and thus the desired result. 
\hfill $\Box$

Below we give the pseudo code for the 
conditional test described in 
Section \ref{sec:simtest}.

\begin{algorithm}[H]
\caption{Simulation of  conditional distribution of
$LR\mid\widehat{C}$ under $H_{0}:\theta_{B}=0$}
\label{alg:conditional_LR}


\KwIn{$\Sigma$, $W$, sample size $n$, index sets $B$ and $U$,
observed set $\widehat{C}$, and number of accepted draws $M$.}

\BlankLine
Set $C=\widehat{C}$\;

\eIf{$C=\emptyset$}{
    Set $LR^{(1)}=\cdots=LR^{(M)}=0$\;
}{
    Set $A_C=B\setminus C$ and $F_C=C\cup U$\;
    Compute $W_{A_C|F_C}$ and $W_{C|U}$\;
    Set $m=0$\;

    \While{$m<M$}{
        Draw $\widetilde{\delta}\sim N(0,\Sigma/n)$\;
        Compute $\delta^{F_C}
        =\widetilde{\delta}_{F_C}
        +W_{F_CF_C}^{-1}W_{F_CA_C}\widetilde{\delta}_{A_C}$\;

        \If{$\delta_C^{F_C}>0$ \textbf{and}
        $W_{A_C|F_C}\widetilde{\delta}_{A_C}\leq0$}{
            Set $m\leftarrow m+1$\;
            Set $LR^{(m)}
            =n(\delta_C^{F_C})^\top W_{C|U}\delta_C^{F_C}$\;
        }
    }
}

\KwOut{$\{LR^{(1)},\ldots,LR^{(M)}\}$.}
\end{algorithm}

\vspace{2mm}
\noindent {\bf Hypothesis 2}

\vspace{2mm}
\noindent \emph{Proof of Lemma \ref{cor:LLRexpansion}}
           We denote $C_1=C \cup U_1$ and decompose $W_{FF}$
    \begin{align*}
        W_{FF} = \left( \begin{matrix}
            W_{C_1 C_1} & W_{C_1 U_2} \\
            W_{U_2 C_1} & W_{U_2U_2}
        \end{matrix} \right)
        \textrm{ and let }W_{C_1|U_2} = W_{C_1C_1} - W_{C_1U_2}W_{U_2U_2}^{-1}W_{U_2C_1}.
    \end{align*}
 The proof is similar to the proof of Theorem~\ref{theorem:LLRexpansion}. The difference is that 
 the case $C_{1}=\emptyset$ does not arise. Whereas
 in the set-up of  
 Theorem 3.1 this does happen. We will show below, this means there cannot be a positive probability that $LR=0$. 
    
By using Lemma \ref{lemma:QuadLLRT}, under the null hypothesis we have
\begin{align*}
LR = 2n\left[\min_{\theta \in \R^{|U_{2}|}}\mathcal{L}_{n}(0_{B},0_{U_{1}},\theta)  -  \min_{\theta \in \R_{+}^{b}\times\R^{p-b}}\mathcal{L}_{n}(\theta)\right] 
  =   2n[Q_n(0_B,0_{U_1},\delta^{U_2}) - Q_n(\widehat{\delta})] + o_p(1), 
\end{align*}
where $\delta^{U_2}= -W_{U_2U_2}^{-1}s_{U_2}$ 
and
$\widehat{\delta} = \argmin_{\delta\in \mathcal{C}}Q_{n}(\delta)$. Hence by using     
 Lemma~\ref{lemma:Quasi-Like} we have 
    \begin{align*}
& 2n[Q_n(0_B,0_{U_1},\delta^{U_2}) - Q_n(\widehat{\delta})] \\  =&   n \sum_{C \in \mathcal{P}(B)} (\delta^{F_C})^\top W_{F_C F_C} \delta^{F_C} 
        I\left[\delta^{F_{C}}_{C} > 0\right]
        I\left[W_{A_{C}|F_{C}}\widetilde{\delta}_{A_{C}} \leq 0 \right] - n(\delta^{U_2})^\top W_{U_2U_2}\delta^{U_2}\\
=& n \sum_{C \in \mathcal{P}(B)}
\left\{
(\delta^{F_C})^\top W_{F_C F_C} \delta^{F_C} -
(\delta^{U_2})^\top W_{U_2U_2}\delta^{U_2}
\right\}
        I\left[\delta^{F_{C}}_{C} > 0\right]
        I\left[W_{A_{C}|F_{C}}\widetilde{\delta}_{A_{C}} \leq 0 \right].
    \end{align*}
    Altogether this gives
    \begin{align*}
        LR 
        = n \sum_{C\in \mathcal{P}(B)} \left\{ (\delta^{F_C})^\top W_{FF} \delta^{F_C} - (\delta^{U_2})^\top W_{U_2U_2} \delta^{U_2} \right\} I\left[\delta^{F_{C}}_{C} > 0\right]
        I\left[W_{A_{C}|F_{C}}\widetilde{\delta}_{A_{C}} \leq 0 \right] + o_p(1).
    \end{align*}
    Now we focus on simplifying the
    the terms within the above summand. 
  By  using Lemma~\ref{lemma:quadratic} we have the 
  simplification 
    \begin{align*}
        \left\{ (\delta^{F_C})^\top W_{FF} \delta^{F_C} - (\delta^{U_2})^\top W_{U_2U_2} \delta^{U_2} \right\} 
        = (\delta^{F_C}_{C_1})^\top W_{C_1|U_2} \delta^{F_C}_{C_1}.
    \end{align*}
    Hence
    \begin{align*}
         LR 
        = n \sum_{\stackrel{C\in \mathcal{P}(B)}{C_{1}=C\cup U_1}} 
        (\delta^{F_C}_{C_1})^\top W_{C_1|U_2} \delta^{F_C}_{C_1} I\left[\delta^{F_{C}}_{C} > 0\right]
        I\left[W_{A_{C}|F_{C}}\widetilde{\delta}_{A_{C}} \leq 0 \right] + o_p(1).
    \end{align*}

    In the special case where $\Sigma=W^{-1}$, then using the Schur complement we have
    $W_{C_1|U_2} = \Sigma_{C_1|A_C}^{-1}$.
Thus 
\begin{align}
\label{eq:LRagain}
         LR 
        = n \sum_{\stackrel{C\in \mathcal{P}(B)}{C_{1}=C\cup U_1}} 
        (\delta^{F_C}_{C_1})^\top \Sigma_{C_1|A_C}^{-1}\delta^{F_C}_{C_1} I\left[\delta^{F_{C}}_{C} > 0\right]
        I\left[\Sigma_{A_{C} A_{C}}^{-1} \widetilde{\delta}_{A_{C}} \leq 0 \right] + o_p(1),
    \end{align}
where we recall from Corollary \ref{theorem:CLT3} that 
$\delta^{F_{C}}$ is asymptotically
independent of 
$\widetilde{\delta}_{A_{C}}$. 
    We recall that under Assumption \ref{assum:CLT}
we have 
$\sqrt{n} \widetilde{\delta} \Dcon N(0, \Sigma)$, 
    $\sqrt{n} \delta^{F_C} \Dcon N(0, \Sigma_{F_C|A_C})$, 
    implying $\sqrt{n} \delta^{F_C}_{C_1} \Dcon N_{|C_1|}(0, \Sigma_{C_1|A_C})$ by $(\Sigma_{F_C|A_C})_{C_1 C_1} = \Sigma_{C_1|A_C}$. Thus we will show that
  \begin{align}
  \label{eq:conditional1}
        n (\delta^{F_C}_{C_1})^\top \Sigma_{C_1|A_C}^{-1} \delta^{F_C}_{C_1}
        \mid (\delta^{F_C}_C>0) \Dcon \chi^2_{|C_1|},
    \end{align}
     In contrast to Corollary~\ref{cor:pivotal}, 
    only a subset $C \subset C_1$ of the coordinates in the quadratic form is constrained to be positive.
    Let $\delta^{F_C}_{C_1}=(X_{1},X_{2})$ and 
    with $X_{1} = \delta^{F_C}_{C}$ and 
    $\Sigma = \Sigma_{C_1|A_C}$, then by applying 
 Lemma~\ref{lemma:trun_chi2} proves (\ref{eq:conditional1}). 
 Combining this with (\ref{eq:LRagain}) proves the final result. 
\hfill $\Box$

\vspace{2mm}

\vspace{2mm}
\noindent {\bf Hypothesis 3}
\vspace{2mm}

\noindent \emph{Proof of Theorem \ref{theorem:typeB}}
We first evaluate an approximation of $LR_3$ under the null hypothesis. Since there does not exist a specific null we  consider an arbitrary $\theta\in \mathcal{C}$ (where $\theta_{B}\in \R_{+}^{b}$ and $\theta_{U}\in \R^{p-b}$). Expanding $LR_{3}$ 
above $\theta\in \mathcal{C}$ we have 
\begin{align*}
LR_{3} \mid \theta_B
=&2n\left[ \min_{\theta\in \mathcal{C}}\mathcal{L}_{n}(\theta) - \min_{\theta\in \R^{p}}\mathcal{L}_{n}(\theta)\right]\nonumber\\
=&2n\left[\min_{\delta\in \mathcal{C}_{\xi}}\left(\delta^{\top}s+\frac{1}{2}\delta^{\top}W\delta \right)+\frac{1}{2}\widetilde{\delta}^{\top}W\widetilde{\delta}\right] + o_{p}(1) \nonumber\\
=&2n \left[Q_{n}(\widehat{\delta}_{\xi}) -\mathcal{L}_n(\theta_0)+ \frac{1}{2}\widetilde{\delta}^{\top}W\widetilde{\delta}\right]+ o_{p}(1) 
\end{align*}
where 
\begin{align*}
\widehat{\delta}_{\xi}
= \sum_{C\in \mathcal{P}(B)}
\begin{pmatrix}
\xi_{A_{C}} \\
\delta^{F_{C}}
\end{pmatrix}
I\!\left[\big(\delta^{F_C}\big)_{C} > \xi_{C} \right]
\,
I\!\left[
W_{A_{C}|F_{C}} \,(\widetilde{\delta}_{A_{C}}-\xi_{A_C}) \leq 0\right]
\end{align*}
Thus by using Lemma \ref{lemma:Quasi-Like}
we obtain (\ref{eq:normexpansion}).

Since $\mathcal{C}_{0}\subseteq\mathcal{C}_{\xi}$ (for $\xi_B\leq 0$), then for any $\theta_{B}\in 
\R_{+}^{b}$ we have  
$Q(\widehat{\delta}_{\xi})\leq Q(\widehat{\delta}_{0})$.
This shows that $(LR_3\mid \theta_B) \leq (LR_3\mid \theta_B=0)+o_p(1)$ and
the asymptotic expression for 
$(LR_{3}\mid \theta_{B}=0)$ is easily deduced from 
(\ref{eq:normexpansion}) by setting ($\xi=0$).

 In the case that $\Sigma = W^{-1}$,  then 
    $W_{A_C|F_C} = \Sigma_{A_C A_C}^{-1}$ (this follows from the block inverse matrix identity).
    Since $\sqrt{n}\widetilde{\delta} \Dcon N(0, \Sigma)$, it follows that $\sqrt{n} \widetilde{\delta}_{A_C} \Dcon N(0, \Sigma_{A_CA_C})$. Finally, by following the same proof as 
    in Lemma \ref{cor:LLRexpansion} we have
\begin{align*}
n(\widetilde{\delta}_{A_C})^\top \Sigma_{A_C A_C}^{-1} \widetilde{\delta}_{A_C} \mid (\delta^{F_{C}}_{C} > 0,
        \Sigma_{A_C A_C}^{-1} \widetilde{\delta}_{A_{C}} <0)\Dcon 
        \chi^{2}_{|A_C|}.
\end{align*}    
This gives the remainder of the result. 
\hfill $\Box$

\vspace{2mm}
\noindent {\bf Hypothesis 4}
The hypothesis is 
\begin{align*}
    H_0: \theta_{B_{1}} =  0 
    \quad \textrm{vs.} \quad
    H_1: \theta_{B_1} \geq 0 \textrm{ and } \theta_{B_1} \neq 0
\end{align*}
where $B_{1}\subset B$ (one can also include a subset of $\theta_{U}$ in the above too).
Let $B_{2} = B\backslash B_{1}$, such that 
$\theta_{B_{2}}$ are the boundary constrained nuisance parameters, $b_{1}=|B_{1}|$ and
$b_2=|B_{2}|$. Let $\xi_{B_2} = -\theta_{B_2} = 
(\xi_{2,1},\ldots,\xi_{2,b_2})$.
Using Lemma \ref{lemma:QuadLLRT} the LR under the null is 
\begin{align*}
LR_4\mid \theta_{B_2} = 2n\left[
\inf_{\delta\in \mathcal{C}_{0,\xi_{B_2}}}Q_{n}(0_{B_1},\delta) - \inf_{\delta\in \mathcal{C}_{1,\xi_{B_2}}}Q_{n}(\delta)
\right] + o_{p}(1),
\end{align*}
where
$\mathcal{C}_{0,\xi_{B_2}} = 
\prod_{j=1}^{b_2} [\xi_{2,j}, \infty)\otimes \R^{p-b}$ 
and 
$\mathcal{C}_{1,\xi_{B_2}} = 
\R_{+}^{b_1}\otimes \mathcal{C}_{0,\xi_{B_2}}$
We now obtain an expression for the above and show that 
it is not possible to obtain the simplifications that we obtained for the previous hypothesis. 

We define
\begin{align*}
W = 
\left( 
\begin{matrix}
W_{B_{1}B_{1}} & W_{B_1B_2} & W_{B_{1}U} \\
W_{B_{2}B_{1}} & W_{B_2B_2} & W_{B_{2}U} \\
W_{UB_{1}} & W_{UB_2} & W_{UU} 
\end{matrix}
\right),
\end{align*}
and
$s = (s_{B_{1}},s_{B_{2}},s_{U})$,
$\widetilde{\delta} = -W^{-1}s$,
$\xi = (0_{B_1},-\theta_{B_2},0_{U})$, $A_{C} = B\backslash C$, $F_{C} = C\cup U$.
We also define the submatrices 
\begin{align*}
\mathcal{W} = 
\left( 
\begin{matrix}
 W_{B_2B_2} & W_{B_{2}U} \\
 W_{UB_2} & W_{UU} 
\end{matrix}
\right),
\end{align*}
and $\mathcal{A}_{C} = B_{2}\backslash C$.

Using the above notation and Lemma \ref{lemma:Quasi-Like}, we have 
\begin{align*}
\min_{\delta \in C_{1,\xi_{B_{2}}}}Q_{n}(\delta) &= \mathcal{L}_{n}(\theta_0)+
\frac{1}{2} P_{1}(\theta_{B_2}) - \frac{1}{2}\widetilde{\delta}^{\top}W\widetilde{\delta} \\
\textrm{ and }
\min_{\delta \in C_{0,\xi_{B_{2}}}}Q_{n}(\delta) &= \mathcal{L}_{n}(\theta_0)+
\frac{1}{2} P_{0}(\theta_{B_2}) - \frac{1}{2}\widetilde{\delta}^{\top}W\widetilde{\delta}
\end{align*}
where 
\begin{align*}
P_{1}(\theta_{B_2}) = 
 \sum_{C \in \mathcal{P}(B_{1}\cup B_{2})} (\xi_{A_C}-\widetilde{\delta}_{A_C})^\top W_{A_C|F_C} (\xi_{A_C}-\widetilde{\delta}_{A_C})
    I[\delta^{F_C}_C>\xi_C] I[W_{A_C|F_C} (\widetilde{\delta}_{A_C}-\xi_{A_C}) \leq 0]
\end{align*}
and 
\begin{align*}
P_{0}(\theta_{B_2}) = 
 \sum_{C_2 \in \mathcal{P}(B_{2})} (\xi_{\mathcal{A}_{C_2}}-\widetilde{\delta}_{\mathcal{A}_{C_2}})^\top \mathcal{W}_{\mathcal{A}_{C_2}|F_{C_2}} (\xi_{\mathcal{A}_{C_2}}-\widetilde{\delta}_{\mathcal{A}_{C_2}})
    I[\delta^{F_{C_2}}_{C_2}>\xi_{C_2}] I[\mathcal{W}_{\mathcal{A}_{C_2}|F_{C_2}} (\widetilde{\delta}_{\mathcal{A}_{C_2}}-\xi_{\mathcal{A}_{C_2}}) \leq 0].
\end{align*}
The other expression in Lemma \ref{lemma:Quasi-Like} can also be used to express $P_{0}(\theta_{B_2})$ and $P_{1}(\theta_{B_2})$ in terms of $\widehat{\delta}$.
Hence 
\begin{align}
\label{eq:LRB}
LR_4 \mid \theta_{B_2} = n \left(P_{0}(\theta_{B_2})-P_{1}(\theta_{B_2}) \right)
+o_{p}(1).
\end{align}
The important observation is that the partition sets 
under the null is not a coarser version of the partition set under the alternative. Precisely, for $C_2 \in \mathcal{P}(B_{2})$ we define
\begin{align*}
R_{0,C_2} &= \{a\in \R^{p};\delta^{F_{C_2}}_{C_2}>\xi_{C_2},\mathcal{W}_{\mathcal{A}_{C_2}|F_{C_2}}(\widetilde{\delta}_{\mathcal{A}_{C_2}}-\xi_{\mathcal{A}_{C_2}}) \leq 0\} 
\end{align*}
and for $C\in \mathcal{P}(B_{1}\times B_{2})$ we define
\begin{align*}
R_{1,C} &= \{a\in \R^{p};\delta^{F_C}_C>\xi_C,W_{A_C|F_C}(\widetilde{\delta}_{A_C}-\xi_{A_C}) \leq 0\},
\end{align*}
where 
$\{R_{1,C};C\in \mathcal{P}(B)\}\not\leq
\{R_{0,C_2};C_{2}\in \mathcal{P}(B_{2})\}$. This condition is necessary to cancel our the nuisance parameters in the $P_{0}$ and $P_{1}$ 
(see Hypothesis 1-3). However, it does not hold for Hypothesis 4.

We cannot make claims such as  
$R_{1,C_{1}\times C_{2}}\subseteq R_{0,C_{2}}$.
Consequently, unlike the previous hypothesis, the nuisance parameters cannot be eliminated by differencing 
$P_{0}(\theta_{B_2})-P_{1}(\theta_{B_2})$ and the distribution of 
$LR$ depends on $\theta_{B_2}$. However, for a given $\theta_{B_2}$
it is possible to use (\ref{eq:LRB}) to simulate from the distribution. 

\vspace{1.5mm}
\noindent {\bf The \cite{p:gui-mc-25} approach}
We now place the method proposed in \cite{p:gui-mc-25}
within our framework. They use the fact that
$P_{0}(\theta_{B_2})$ and $P_{1}(\theta_{B_2})$ are both monotonically decreasing as $\theta_{B_2}$ moves away from the origin (see Lemma \ref{lemma:monotonic}). Using this intriguing observation they 
propose to use pointwise confidence interval for each entry in 
$\theta_{B_2}$ (based on $\widetilde{\theta}_{B_2,L}$)
as a means of obtaining an upper bound for 
$LR$ under the null. For each $1\leq j \leq b_2$ let 
$[\widetilde{\theta}_{j,\alpha/2}, \widetilde{\theta}_{j,1-\alpha/2}]$ denote a $(1-\alpha)$-level confidence interval for $\theta_{B,j}$. 
Then they define 
\begin{align*}
LR(\widetilde{\theta}_{B,\alpha/2}, \widetilde{\theta}_{B,1-\alpha/2})
= P_{0}(\widetilde{\theta}_{B_2,\alpha/2}) - P_{1}(\widetilde{\theta}_{B_2,1-\alpha/2}).
\end{align*}
They show that the critical values of the above statistic
keep the level of the test and is also powerful against all alternatives. It is interesting to note that 
$P_{0}(\widetilde{\theta}_{B_2,\alpha/2}) - P_{1}(\widetilde{\theta}_{B_2,1-\alpha/2})$ can be simulated using the formulas above.

\section{Examples}\label{sup:examples}

\subsection{Examples (a short summary)}

\begin{example}
We illustrate the case $p=1$ and $B=\{1\}$ considered in Example~\ref{example:p1}.
Under $\sqrt{n} \widetilde{\delta} \Dcon N(0,\sigma^{2})$,
the asymptotic distribution of $\widehat{\delta}$ in this case is immediate,
but we here derive it using Theorem~\ref{theorem:CLT}.
Since $U=\emptyset$, 
$T_{C}$ is a scalar.
Specifically, for the case (i) $C=\emptyset$, then $T_{\emptyset} = W$ and for the case (ii) $C=\{1\}$, then $T_{\{1\}}=1$.
Using Theorem~\ref{theorem:CLT},
\begin{align*}
    \sqrt{n} \widehat{\delta}
    \Dcon m_{\{1\}}(0) \Phi\left(0,\frac{\sigma^{2}}{n}\right)
    + tN_1(0, \sigma^{2})\left[1-\Phi\left(0,\frac{\sigma^{2}}{n}\right)\right].
\end{align*}
\end{example}

\begin{example}\label{example:p2-dist}
We return to the example $p=2$ and $B=\{\1,2\}$ considered in Example \ref{example:p2}. 
Given $W$ and $\Sigma$, 
we use 
Theorem \ref{theorem:CLT} to obtain the asymptotic distribution of $\widehat{\delta}$ in this case.
Since $U=\emptyset$, the block corresponding to $U$ disappears, and 
$T_{C}$ is a $2\times 2$-dimensional matrix.
We first obtain $V_{C}$ for $C\in \{\{1\},\{2\},\{1,2\},\emptyset\}$
\begin{itemize}
    \item[(i)] If $C=\{1\}$, then $A_C=\{2\}$ and $\sqrt{n} T_{\{1\}} \widetilde{\delta} \Dcon N(0, V^{\{1\}})$ where
    \begin{align*}
        V^{\{1\}} 
        &= T_{\{1\}} \Sigma T_{\{1\}}^\top
        = \left( \begin{matrix}
            \Sigma_{11} + 2\Sigma_{12}W_{12}W_{11}^{-1}+ \Sigma_{22}W_{12}^{2} W_{11}^{-2} &
            (\Sigma_{12}+W_{11}^{-1}W_{12}\Sigma_{22}) W_{2|1} \\
            (\Sigma_{12}+W_{11}^{-1}W_{12}\Sigma_{22}) W_{2|1}
            & \Sigma_{22}W_{2|1}^2
        \end{matrix} \right),        
    \end{align*}
    with $T_{\{1\}} = \left( \begin{matrix}
                1 & W_{11}^{-1}W_{12}  \\
                0 & W_{2|1}
            \end{matrix}\right).$
    \item[(ii)] If $C=\{2\}$, then $A_C=\{1\}$  and $\sqrt{n} T_{\{2\}} \widetilde{\delta} \Dcon N(0, V^{\{2\}})$ where
    \begin{align*}
        V^{\{2\}}
        &= T_{\{2\}} \Sigma T_{\{2\}}^\top
        = \left( \begin{matrix}
            \Sigma_{11}W_{1|2}^2 & (\Sigma_{12} + W_{22}^{-1}W_{21}\Sigma_{11})W_{1|2} \\
            (\Sigma_{12} + W_{22}^{-1}W_{21}\Sigma_{11})W_{1|2} 
            & \Sigma_{22} + 2\Sigma_{12}W_{22}^{-1}W_{21} + \Sigma_{11}W_{22}^{-2}W_{21}^{2}
        \end{matrix} \right)
    \end{align*}
    with $T_{ \{2\} } = \left( \begin{matrix}
                W_{1|2} & 0 \\
                W_{22}^{-1}W_{21} & 1
            \end{matrix}\right)$.
\end{itemize}
Using Theorem \ref{theorem:CLT} the asymptotic distribution of 
$\widehat{\delta}$ is
\begin{align*}
    \sqrt{n}\widehat{\delta} 
    &\Dcon tN_2(0,\Sigma) \operatorname{pr}(\widetilde{\delta}_1>0, \widetilde{\delta}_2>0) 
    + \left( \begin{matrix}
        tSN_{1,1}(0, 0, V^{\{1\}}_{21}/ V^{\{1\}}_{11}, V^{\{1\}}_{11}, V^{\{1\}}_{2|1}) \\
        m_{\{2\}}(0)
    \end{matrix}\right) \operatorname{pr}(\delta^1>0, \widetilde{\delta}_2<0) \\
    &+ \left( \begin{matrix}
         m_{\{1\}}(0) \\
         tSN_{1,1}(0, 0, V^{\{2\}}_{12}/ V^{\{2\}}_{22}, V^{\{2\}}_{22}, V^{\{2\}}_{1|2})
    \end{matrix}\right) \operatorname{pr}(\delta^2>0, \widetilde{\delta}_1<0)
    + \left(
    \begin{matrix}
    m_{\{1\}}(0) \\
    m_{\{2\}}(0) \\
    \end{matrix}
    \right)
    \operatorname{pr}(\delta^1<0, \delta^2<0).
\end{align*}
Note that, when $C=\{1,2\}$, then $A_C=\emptyset$ and
$tSN_{0,2}(0,0,0,V_{CC},0)=tN(0,\Sigma)$.

\end{example}

\begin{example}[The case that $p=2$, $B=\{1,2\}$ and $\Sigma = W^{-1}$]
By using Corollary \ref{theorem:CLT3} the asymptotic distribution in Example \ref{example:p2-dist} reduces
\begin{align*}
    \sqrt{n} \widehat{\delta}
    &\Dcon
    tN_2(0, \Sigma) \operatorname{pr}(\widetilde{\delta}_1>0, \widetilde{\delta}_2>0)
    + \left( \begin{matrix}
        tN_1(0, \Sigma_{1|2}) \\
        m_{\{2\}}(0)
    \end{matrix}\right) \operatorname{pr}(\delta^1>0, \widetilde{\delta}_2<0) \\
    & + \left( \begin{matrix}
        m_{\{1\}}(0)\\
        tN_1(0, \Sigma_{2|1}) 
    \end{matrix}\right) \operatorname{pr}(\delta^2>0, \widetilde{\delta}_1<0) 
    + \left(
    \begin{matrix}
    m_{\{1\}}(0) \\
    m_{\{2\}}(0) \\
    \end{matrix}
    \right)\operatorname{pr}(\delta^1<0,\delta^2<0).
\end{align*}
\end{example}

\subsection{Case $p=2$ (in detail)} \label{sec:example_p2}
Consider the case that $p=2$, $B=\{1,2\}, \mathcal{P}(B)=\{\emptyset, \{1\}, \{2\}, \{1,2\}\}$ and $U=\emptyset$.
Suppose the population parameter lies on the boundary (thus $\xi=0$).
Recall from \eqref{eq:casep2},
\begin{align*}
    \widehat{\delta} 
    &= \left(\begin{matrix}
        \widetilde{\delta}_1 \\
        \widetilde{\delta}_2
    \end{matrix}\right) I(\widetilde{\delta}_1>0, \widetilde{\delta}_2>0)
    + \left(\begin{matrix}
        \delta^1 \\
        0
    \end{matrix}\right) I(\delta^1>0, \widetilde{\delta}_2<0) \\
    & \quad + \left(\begin{matrix}
        0 \\
        \delta^2
    \end{matrix}\right) I(\widetilde{\delta}_1<0, \delta^2>0)
    + \left(\begin{matrix}
        0 \\
        0
    \end{matrix}\right) I(\delta^1<0, \delta^2<0),
\end{align*}
where $\widetilde{\delta}= \argmin_{\delta \in \R^2} Q(\delta) = -W^{-1}s$, 
$\delta^1 = \argmin_{\delta \in \R} Q(\delta,0) = \widetilde{\delta}_1 + W_{11}^{-1}W_{12} \widetilde{\delta}_2$, and 
$\delta^2 = \argmin_{\delta \in \R} Q(0,\delta) = \widetilde{\delta}_2 + W_{22}^{-1}W_{21} \widetilde{\delta}_1$.

To study the asymptotic distribution of $\widehat{\delta}$, we consider the conditional distribution associated with each term in the above.
To this end, we use a matrix $T_C$ defined in Theorem~\ref{theorem:CLT}.
Since $U=\emptyset$, the block corresponding to $U$ disappears, and $T_C$ reduces to a $2 \times 2$ matrix.
Specifically,
\begin{itemize}
    \item[(i)] If $C=\{1\}$, then $A_C=\{2\}$ and $\sqrt{n} T_{\{1\}} \widetilde{\delta} \Dcon N(0, V^{\{1\}})$ where
    \begin{align*}
        V^{\{1\}} 
        &= T_{\{1\}} \Sigma T_{\{1\}}^\top
        = \left( \begin{matrix}
            \Sigma_{11} + 2\Sigma_{12}W_{12}W_{11}^{-1}+ \Sigma_{22}W_{12}^{2} W_{11}^{-2} &
            (\Sigma_{12}+W_{11}^{-1}W_{12}\Sigma_{22}) W_{2|1} \\
            (\Sigma_{12}+W_{11}^{-1}W_{12}\Sigma_{22}) W_{2|1}
            & \Sigma_{22}W_{2|1}^2
        \end{matrix} \right) \\
        &=
        \left(
        \begin{matrix}
                v_{11}^{\{1\}} & v_{12}^{\{1\}} \\
         v_{21}^{\{1\}} & v_{22}^{\{1\}} \\
\end{matrix}
\right)
    \end{align*}
    and $T_{\{1\}} = \left( \begin{matrix}
                1 & W_{11}^{-1}W_{12}  \\
                0 & W_{2|1}
            \end{matrix}\right).$
    \item[(ii)] If $C=\{2\}$, then $A_C=\{1\}$  and $\sqrt{n} T_{\{2\}} \widetilde{\delta} \Dcon N(0, V^{\{2\}})$ where
    \begin{align*}
        V^{\{2\}}
        &= T_{\{2\}} \Sigma T_{\{2\}}^\top
        = \left( \begin{matrix}
            \Sigma_{11}W_{1|2}^2 & (\Sigma_{12} + W_{22}^{-1}W_{21}\Sigma_{11})W_{1|2} \\
            (\Sigma_{12} + W_{22}^{-1}W_{21}\Sigma_{11})W_{1|2} 
            & \Sigma_{22} + 2\Sigma_{12}W_{22}^{-1}W_{21} + \Sigma_{11}W_{22}^{-2}W_{21}^{2}
        \end{matrix} \right)\\
        &=
        \left(
        \begin{matrix}
        v_{11}^{\{2\}} & v_{12}^{\{2\}} \\
         v_{21}^{\{2\}} & v_{22}^{\{2\}} \\
         \end{matrix}
        \right)
    \end{align*}
    where $T_{ \{2\} } = \left( \begin{matrix}
                W_{1|2} & 0 \\
                W_{22}^{-1}W_{21} & 1
            \end{matrix}\right)$.
\end{itemize}
Using Lemma~\ref{lemma:SN_TSN} gives the conditional distributions for $C=\{1\}$ and $C=\{2\}$.
As a result, the asymptotic distribution is
\begin{align*}
    \sqrt{n}\widehat{\delta} 
    &\Dcon tN_2(0,\Sigma) \operatorname{pr}(\widetilde{\delta}_1>0, \widetilde{\delta}_2>0) 
    + \left( \begin{matrix}
        tSN_{1,1}(0, 0, v^{\{1\}}_{21}/ v^{\{1\}}_{11}, v^{\{1\}}_{11}, v^{\{1\}}_{2|1}) \\
        m_{\{2\}}(0)
    \end{matrix}\right) \operatorname{pr}(\delta^1>0, \widetilde{\delta}_2<0) \\
    &+ \left( \begin{matrix}
         m_{\{1\}}(0) \\
         tSN_{1,1}(0, 0, v^{\{2\}}_{12}/ v^{\{2\}}_{22}, v^{\{2\}}_{22}, v^{\{2\}}_{1|2})
    \end{matrix}\right) \operatorname{pr}(\delta^2>0, \widetilde{\delta}_1<0)
    + \left( \begin{matrix}
        m_{\{1\}}(0) \\
        m_{\{2\}}(0)
    \end{matrix} \right) \operatorname{pr}(\delta^1<0, \delta^2<0).
\end{align*}
Note that, when $C=\{1,2\}$, then $A_C=\emptyset$, thus   
$tSN_{0,2}(0,0,0,V_{CC},0)=tN_{2}(0,\Sigma)$.

\vspace{2mm}
\noindent {\bf Special Case $\Sigma = W^{-1}$}
Applying Corollary \ref{theorem:CLT3}
\begin{align*}
    \sqrt{n} \widehat{\delta}
    &\Dcon
    tN_2(0, \Sigma) \operatorname{pr}(\widetilde{\delta}_1>0, \widetilde{\delta}_2>0)
    + \left( \begin{matrix}
        tN_1(0, \Sigma_{1|2}) \\
        m_{\{2\}}(0)
    \end{matrix}\right) \operatorname{pr}(\delta^1>0, \widetilde{\delta}_2<0) \\
    & + \left( \begin{matrix}
        m_{\{1\}}(0)\\
        tN_1(0, \Sigma_{2|1}) 
    \end{matrix}\right) \operatorname{pr}(\delta^2>0, \widetilde{\delta}_1<0) 
    + \left( \begin{matrix}
        m_{\{1\}}(0) \\
        m_{\{2\}}(0)
    \end{matrix} \right) \operatorname{pr}(\delta^1<0,\delta^2<0).
\end{align*}
Indeed, when $\Sigma =W^{-1}$, we have $v_{21}^{\{1\}}=v_{12}^{\{2\}}=0$. Thus the $tSN_{1,1}$ reduces to ordinary univariate truncated normal distributions.

\vspace{2mm}
\noindent {\bf The Quasi LR test for the case $p=2$}
We consider the one-sided hypothesis 
\begin{align*}
    H_0: \theta_B = 0 \quad \textrm{vs.} \quad
    H_1: \theta_B \geq 0 \textrm{ and }\theta_{B}\neq 0,
\end{align*}
in the case that $B=\{1,2\}$ and $U=\emptyset$.
From Theorem~\ref{theorem:LLRexpansion},
\begin{align*}
    LR 
    &= 
    n \bigg( \widetilde{\delta}^\top W \widetilde{\delta} \cdot I(\widetilde{\delta}_1>0, \widetilde{\delta}_2>0)
    + W_{11} (\delta^1)^2 \cdot I(\delta^1>0, \widetilde{\delta}_2<0) 
    + W_{22} (\delta^2)^2 \cdot I(\delta^2>0, \widetilde{\delta}_1<0)
    \bigg) \\
    &+ 0 \cdot I(\delta^1<0, \delta^2<0) + o_p(1).
\end{align*}
In the special case where $\Sigma = W^{-1}$, 
applying 
Corollary \ref{cor:pivotal}
gives the well known result that 
$LR$ is a mixture of 
$\chi^{2}$-distributions
\begin{align*}
    LR
    &\Dcon \chi_2^2 \cdot \operatorname{pr}(\widetilde{\delta}_1>0, \widetilde{\delta}_2>0)
    + \chi_1^2 \cdot \big( \operatorname{pr}(\delta^1>0, \widetilde{\delta}_2<0)
    +  \operatorname{pr}(\delta^2>0, \widetilde{\delta}_1<0) \big)
    + \chi_0^2 \cdot \operatorname{pr}(\delta^1<0, \delta^2<0).
\end{align*}

\vspace{2mm}
\noindent {\bf Impact of boundary on finite sample distribution (for the case $p=2$)}\label{sec:Boundary}
We consider the case where $p=2, B=\{1\}, U=\{2\}$ and allow for the case that population parameter does not lie on the boundary (that is $\xi_1 < 0$).
From Theorem~\ref{theorem:blocks},
\begin{align} \label{eq:sim_p2_xi}
    \widehat{\delta}
    = \underbrace{\left( \begin{matrix}
        \widetilde{\delta}_1 \\
        \widetilde{\delta}_2
    \end{matrix} \right) I[\widetilde{\delta}_1 > \xi_1]}_{C=\{1\}}
    + \underbrace{\left( \begin{matrix}
        \xi_1 \\
        \widetilde{\delta}_2 + W_{22}^{-1}W_{21}(\widetilde{\delta}_1-\xi_1)
    \end{matrix} \right) I[\widetilde{\delta}_1 < \xi_1]}_{C=\emptyset}.
\end{align}
From Equation~\eqref{eq:finite_approx},
\begin{align*}
    \widehat{\theta}-\theta_0
    &\approxD 
    \underbrace{\left( \begin{matrix}
        tN_1(\xi_1, \Sigma_{11}/n) \\
        \widetilde{SN}_{1,1}(\xi_1, \Sigma_{12}\Sigma_{22}^{-1},\Sigma_{22}/n, \Sigma_{1|2}/n)
    \end{matrix} \right)
    \operatorname{pr}[\widetilde{\delta}_1 > \xi_1]}_{C=\{1\}} \\
    &+ \underbrace{\left[
    \left( \begin{matrix}
        m_{\{1\}}(\xi_1) \\
        \widetilde{SN}_{1,1}(-\xi_1, -V_{12}V_{22}^{-1},V_{22}/n, V_{1|2}/n)
    \end{matrix} \right)
    - \left( \begin{matrix}
        0 \\
        W_{22}^{-1}W_{21}\xi_1
    \end{matrix} \right)\right]
    \operatorname{pr}[\widetilde{\delta}_1 < \xi_1]}_{C=\emptyset},
\end{align*}
where 
$$V = \left( \begin{matrix}
    \Sigma_{11} & \Sigma_{12} + W_{22}^{-1}W_{21}\Sigma_{11} \\
    \Sigma_{12} + W_{22}^{-1}W_{21}\Sigma_{11} &
    \Sigma_{22} + 2 W_{22}^{-1}W_{21}\Sigma_{12} +( W_{22}^{-1}W_{21})^2 \Sigma_{11}
\end{matrix}\right).$$

We focus on the unconstrained 
estimator $\widehat{\theta}_{2}$.
The conditional distribution of
$\widehat{\theta}_{2} \mid C=\emptyset$
 is one-sided skewed normal distribution with shift $-W_{22}^{-1}W_{21}\xi_1$. 
 Whereas 
the conditional distribution of $\widehat{\theta}_{2} \mid C=\{1\}$ is a one-sided skewed normal distribution (no shift).
We consider three distinct examples below. 

\begin{figure}
    \centering
    \begin{subfigure}[t]{0.48\linewidth}
        \centering
        \includegraphics[width=\linewidth]{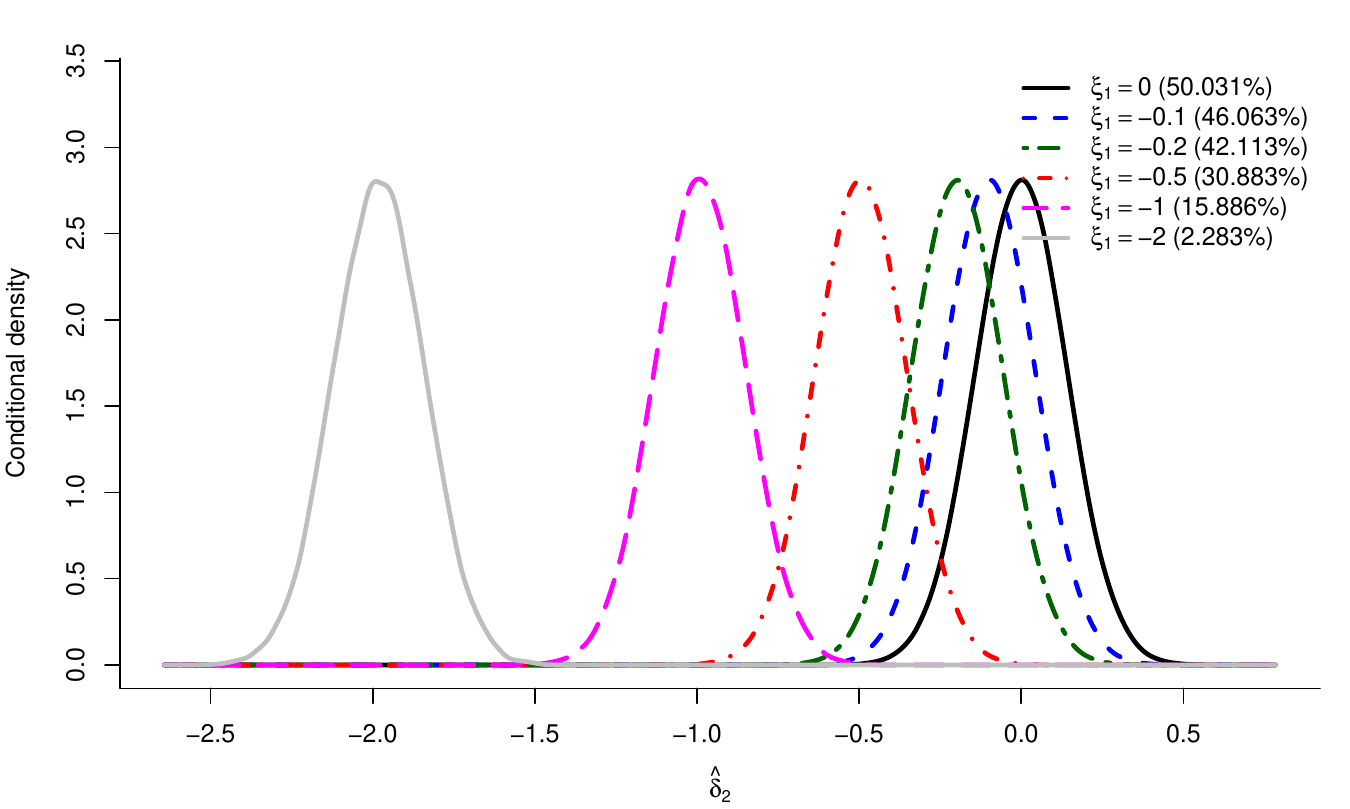}
        \caption{$\widehat{\delta}_2 \mid C=\emptyset$}
        \label{fig:p2_xi_cond_C0}
    \end{subfigure}\hfill%
    \begin{subfigure}[t]{0.48\linewidth}
        \centering
        \includegraphics[width=\linewidth]{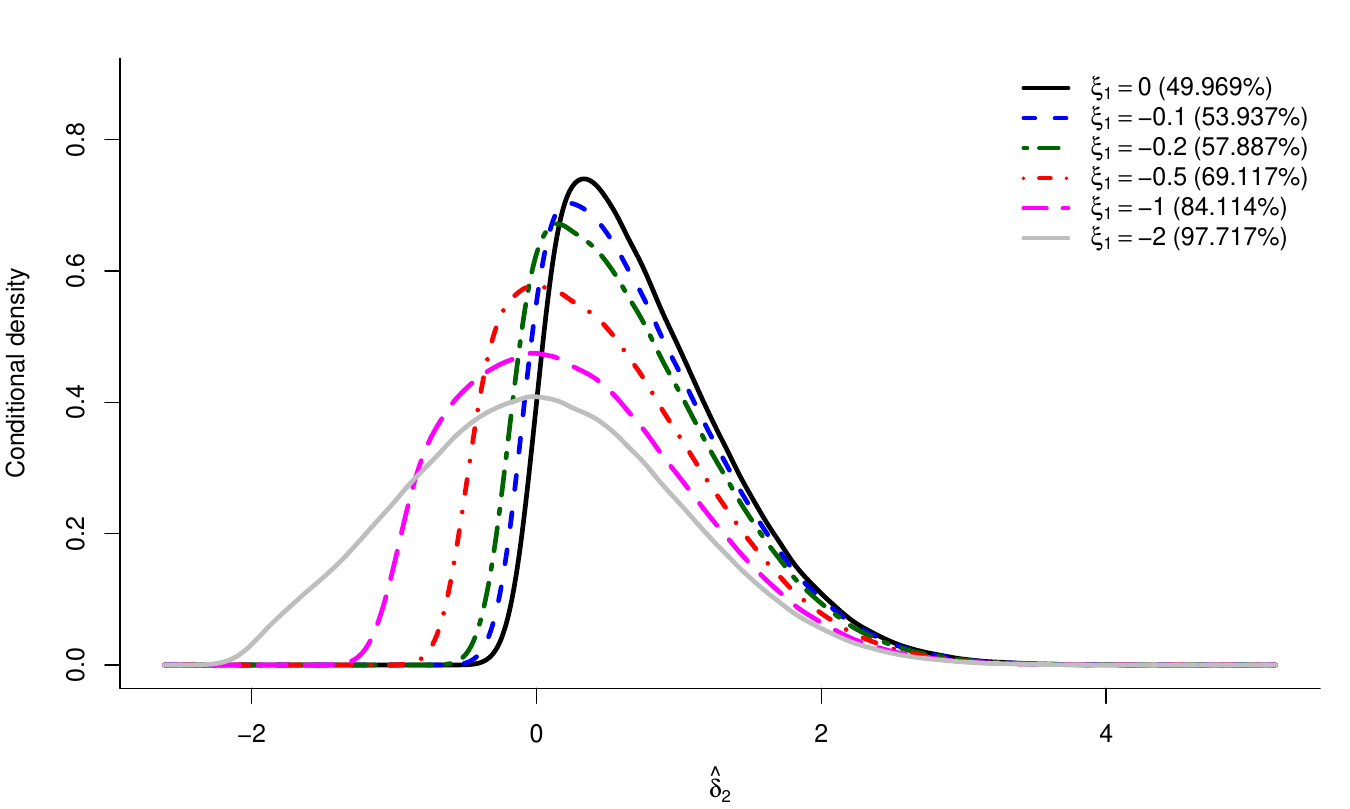}
        \caption{$\widehat{\delta}_2 \mid C=\{1\}$}
        \label{fig:p2_xi_cond_C1}
    \end{subfigure}

    \par\medskip

    \begin{subfigure}[t]{0.48\linewidth}
        \centering
        \includegraphics[width=\linewidth]{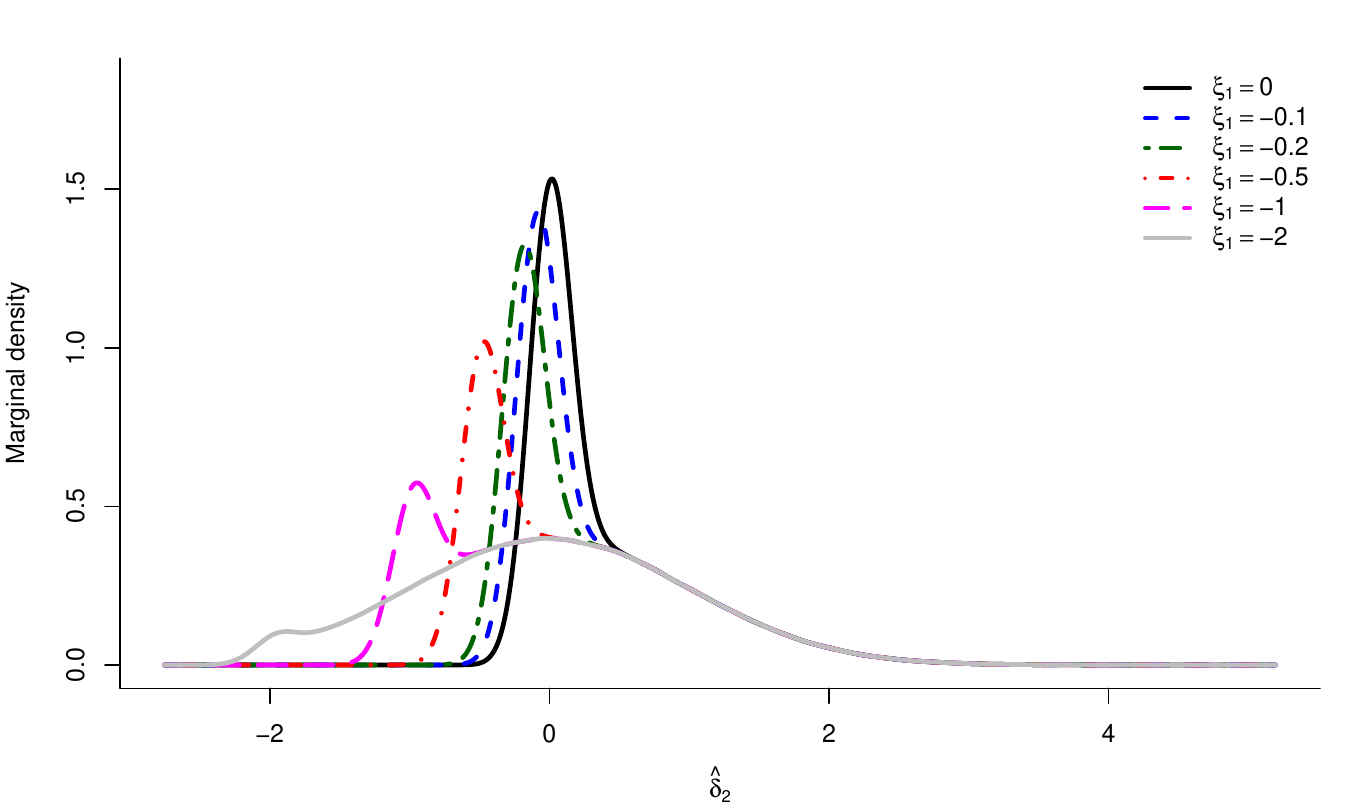}
        \caption{Marginal}
        \label{fig:p2_xi_marginal}
    \end{subfigure}
    \caption{
        Example 1. Density estimates of $\widehat{\delta}_2$ for $\xi_1\in\{0,-0.1,-0.2,-0.5,-1, -2\}$:
        (a) conditional on $C=\emptyset$,
        (b) conditional on $C=\{1\}$, and
        (c) marginal.
    }
    \label{fig:p2_xi_cond}
\end{figure}

\vspace{2mm}
\noindent {\bf Example 1}
To explore the boundary effect on the finite sample distribution,
we set
\begin{align}
\label{eq:example1}
    W = \left( \begin{matrix}
        1 & -0.99 \\
        -0.99 & 1
    \end{matrix}\right),
    \qquad
    \Sigma = \left( \begin{matrix}
        1 & 0.99 \\
        0.99 & 1
    \end{matrix} \right).
\end{align}
For each $\xi_1 \in \{0, -0.1, -0.2, -0.5, -1, -2\}$,
we generate $M=2000000$ of $\widetilde{\delta} \sim N(0, \Sigma)$ and compute $\widehat{\delta}$ using the representation in \eqref{eq:sim_p2_xi}.
Figure~\ref{fig:p2_xi_cond} shows the estimated conditional
and marginal densities of $\widehat{\delta}_2$, and Table~\ref{tab:p2_xi_skewness} reports their skewness as measured by 
\begin{align*}
   S_{3}= \frac{1/M \sum_{i=1}^M(x_i-\bar{x})^3}{(1/M  \sum_{i=1}^M (x_i-\bar{x})^2  )^{3/2}}.
\end{align*}
We observe that the conditional density of $\widehat{\delta}_{2}\mid C=\emptyset$ looks close to normal with almost zero skewness. The reason is that for the given $W$ and $\Sigma$ in (\ref{eq:example1}), we can show that
$V_{12} = \Sigma_{12}+W_{22}^{-1}W_{21}\Sigma_{11} = 0$.
This implies the limiting distribution for $C = \emptyset$ is
\\*
$\widetilde{SN}_{1,1}(-\xi_1, -V_{12}V_{22}^{-1},V_{22}/n, V_{1|2}/n)-W_{22}^{-1}W_{21}\xi_{1} = N(0, 0.0199)+0.99\xi_{1}$, whose skewness is zero.
In Figure~\ref{fig:p2_xi_cond_C0}, the conditional densities for $C=\emptyset$ 
closely follow these shifted normal densities.
On the other hand, 
the conditional densities for $C=\{1\}$ show clear right skewness and the skewness  decreases as $\xi_1$ becomes more negative.
The marginal densities also exhibits right skewness and bimodality
which decreases as $\xi_1$ approaches to $-2$.
Each marginal density is a mixture of the two conditional densities weighted by $\operatorname{pr}(C=\emptyset)$ and $\operatorname{pr}(C=\{1\})$, respectively.
In Figure~\ref{fig:p2_xi_cond}, the approximate probabilities for each mixture is reported in the legend.

\begin{table}[htbp]
    \centering
    \caption{Example 1. Skewness of the conditional and marginal distributions of $\widehat{\delta}_2$.}
    \label{tab:p2_xi_skewness}
    \begin{tabular}{rrrr}
        \toprule
        $\xi_1$ & $C=\emptyset$ & $C=\{1\}$ & Marginal \\
        \midrule
         0.0 & -0.0007 & 0.9187 &  1.5056 \\
        -0.1 & -0.0008 & 0.8875 &  1.3846 \\
        -0.2 & -0.0007 & 0.8551 &  1.2698 \\
        -0.5 & -0.0013 & 0.7514 &  0.9613 \\
        -1.0 & -0.0045 & 0.5645 &  0.5597 \\
        -2.0 & -0.0100 & 0.2139 &  0.1209 \\
        \bottomrule
    \end{tabular}
\end{table}

\begin{figure}[h!]
    \centering
    \begin{subfigure}[t]{0.48\linewidth}
        \centering
        \includegraphics[width=\linewidth]{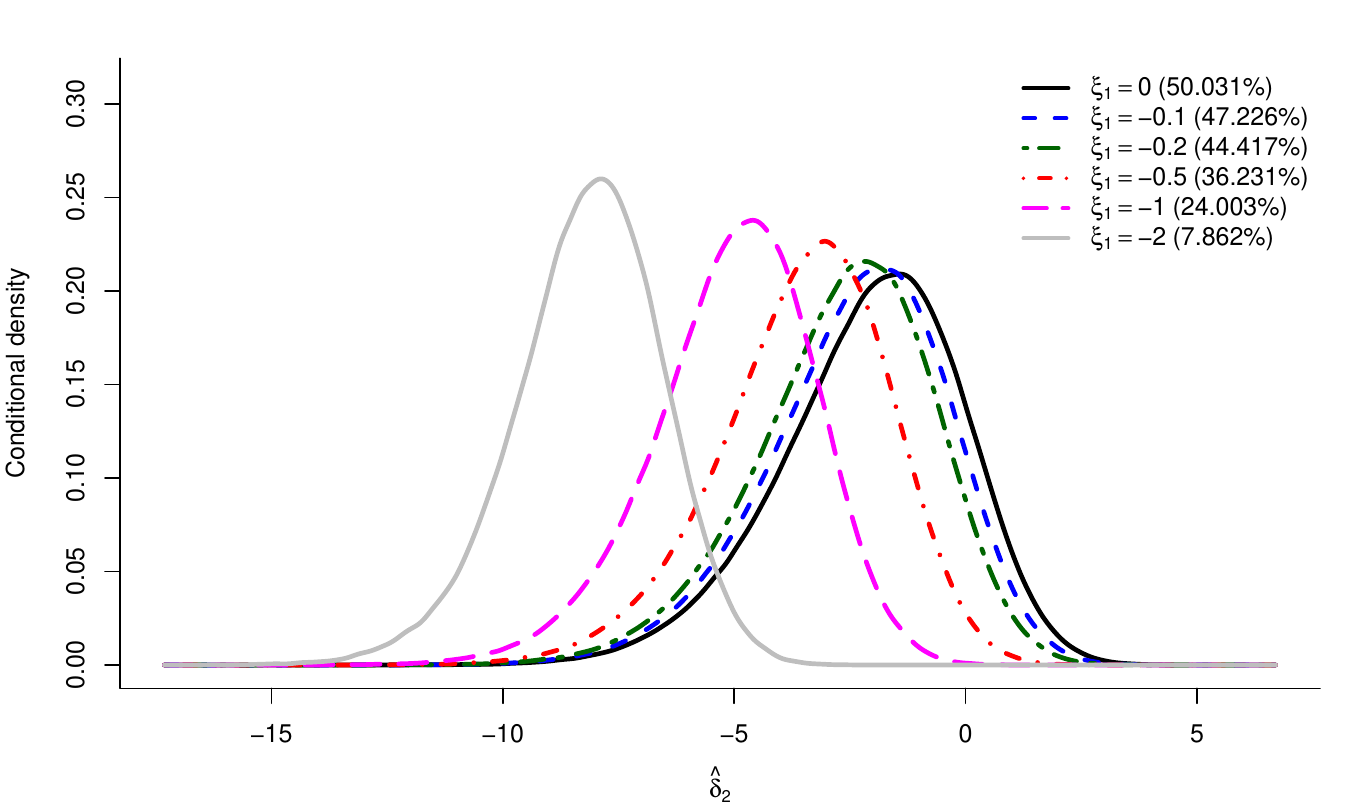}
        \caption{$\widehat{\delta}_2 \mid C=\emptyset$}
        \label{fig:p2_xi_cond_C0_example2}
    \end{subfigure}\hfill%
    \begin{subfigure}[t]{0.48\linewidth}
        \centering
        \includegraphics[width=\linewidth]{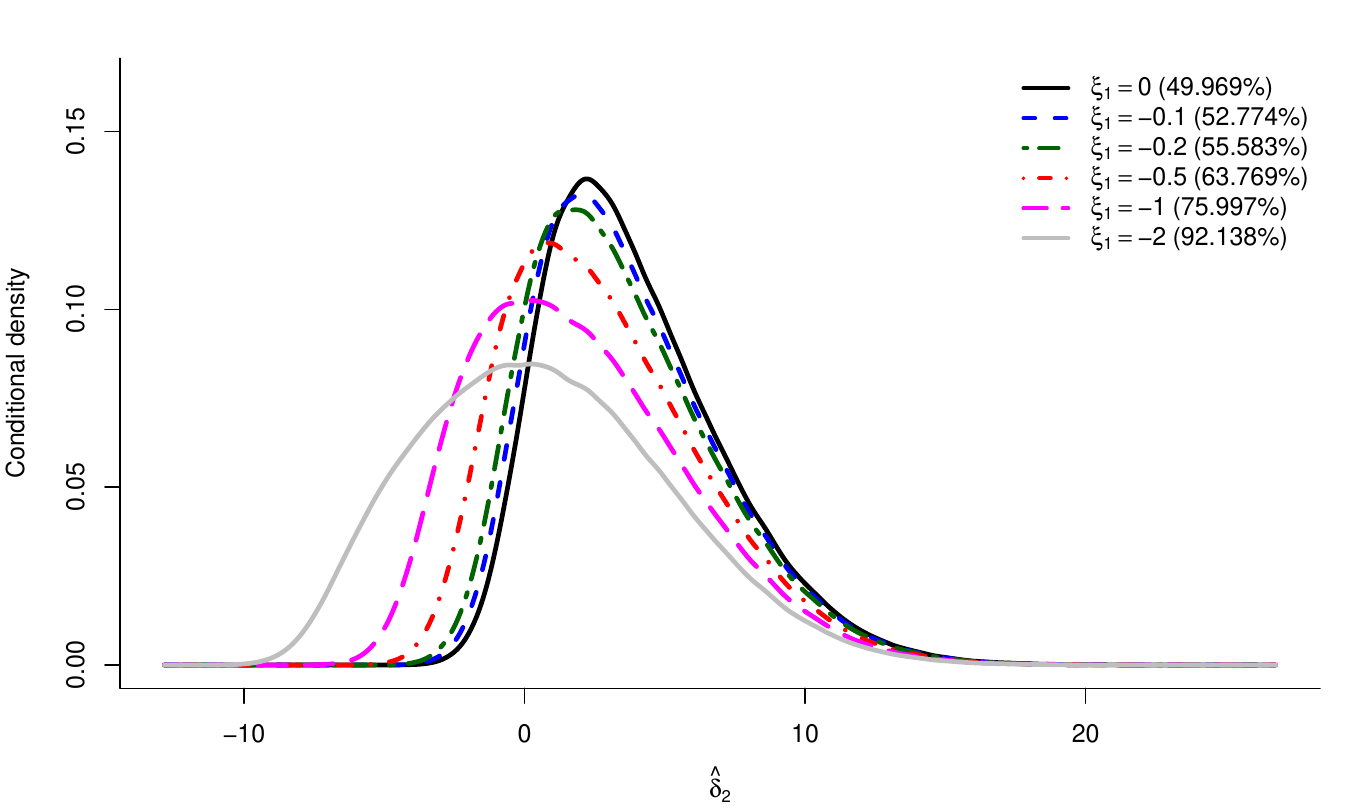}
        \caption{$\widehat{\delta}_2\mid C=\{1\}$}
        \label{fig:p2_xi_cond_C1_example2}
    \end{subfigure}

    \par\medskip

    \begin{subfigure}[t]{0.48\linewidth}
        \centering
        \includegraphics[width=\linewidth]{figures/sim_p2_xi_marginal_v3.pdf}
        \caption{Marginal}
        \label{fig:p2_xi_marginal_example2}
    \end{subfigure}
    \caption{Example 2.
        Density estimates of $\widehat{\delta}_2$ for $\xi_1\in\{0,-0.1,-0.2,-0.5,-1, -2\}$:
        (a) conditional on $C=\emptyset$,
        (b) conditional on $C=\{1\}$, and
        (c) marginal.
    }
    \label{fig:p2_xi_cond_example2}
\end{figure}

\vspace{2mm}
\noindent {\bf Example 2}
For another example,
we consider
\begin{align*}
    W = \left( \begin{matrix}
        10 & -5 \\
        -5 & 3
    \end{matrix}\right),
    \qquad
    \Sigma = \left( \begin{matrix}
        2 & 7 \\
        7 & 26
    \end{matrix} \right),
\end{align*}
and follow the same procedure as in the previous example for each $\xi_1 \in \{0, -0.1, -0.2, -0.5, -1, -2\}$.
Figure~\ref{fig:p2_xi_cond_example2} presents the estimated conditional and marginal densities of $\widehat{\delta}_2$ and Table~\ref{tab:p2_xi_skewness_example2} reports the skewness computed using the same formula as in the previous example.
Unlike in the previous example, the conditional densities of $\widehat{\delta_2}\mid C=\emptyset$ are skewed, with skewness decreasing as $\xi_1$ decreases. 
This pattern is clear in both Figure~\ref{fig:p2_xi_cond_C0_example2} and Table~\ref{tab:p2_xi_skewness_example2}.
In this example, $V_{12}$ is no longer zero, so the limiting distribution is one-sided skewed normal.
In addition, the conditional densities of $\widehat{\delta_2}\mid C=\{1\}$ show the right skewness and the skewness get weakened as $\xi_1$ gets closer to $-2$.
The marginal densities also show the clear right skewness, but bimodality is not obvious in this example.
Recall that each marginal densities is a mixture of the two conditional densities weighted by $\operatorname{pr}(C=\emptyset)$ and $\operatorname{pr}(C=\{1\})$, respectively. 
The legend in Figure~\ref{fig:p2_xi_cond_example2} shows the approximated probabilities for each mixture.

\begin{table}[htbp]
    \centering
    \caption{Example 2. 
    Skewness of the conditional and marginal distributions of $\widehat{\delta}_2$.}
    \label{tab:p2_xi_skewness_example2}
    \begin{tabular}{rrrr}
        \toprule
        $\xi_1$ & $C=\emptyset$ & $C=\{1\}$ & Marginal \\
        \midrule
         0.0 & -0.4839 & 0.7881 &  0.6190 \\
        -0.1 & -0.4850 & 0.7734 &  0.6016 \\
        -0.2 & -0.4844 & 0.7581 &  0.5823 \\
        -0.5 & -0.4810 & 0.7072 &  0.5157 \\
        -1.0 & -0.4581 & 0.6078 &  0.3909 \\
        -2.0 & -0.3810 & 0.3767 &  0.1699 \\
        \bottomrule
    \end{tabular}
\end{table}

\begin{figure}[h!]
    \centering
    \begin{subfigure}[t]{0.48\linewidth}
        \centering
        \includegraphics[width=\linewidth]{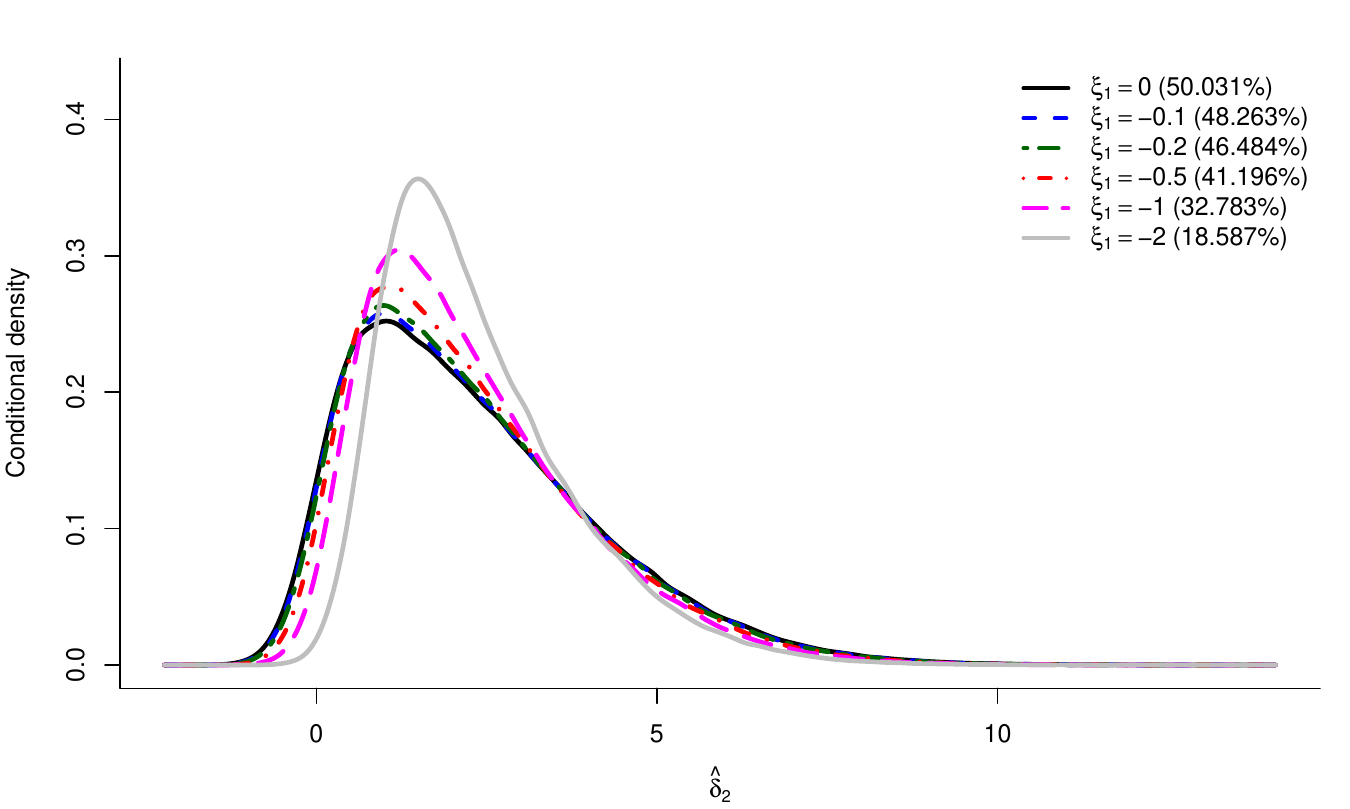}
        \caption{$\widehat{\delta}_2\mid C=\emptyset$}
        \label{fig:p2_xi_cond_C0_example3}
    \end{subfigure}\hfill%
    \begin{subfigure}[t]{0.48\linewidth}
        \centering
        \includegraphics[width=\linewidth]{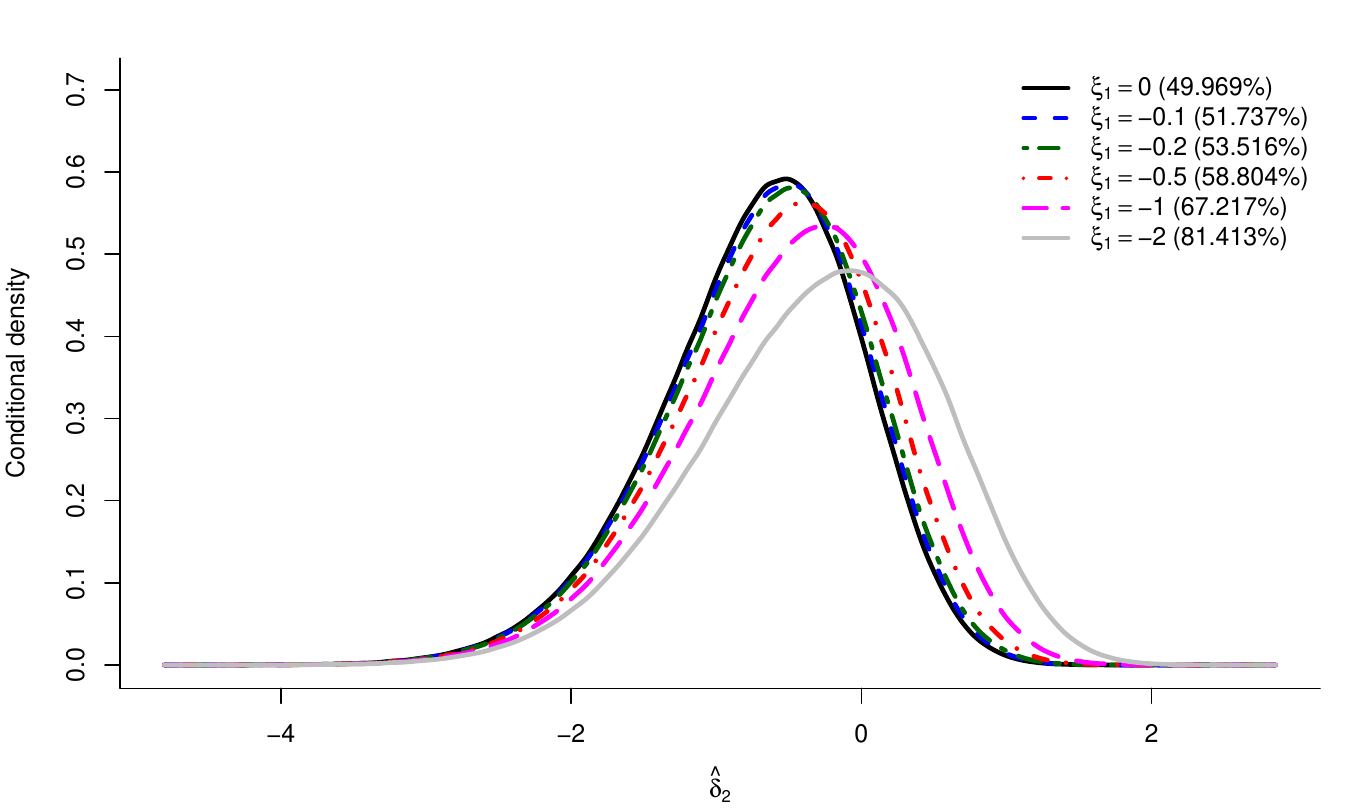}
        \caption{$\widehat{\delta}_2\mid C=\{1\}$}
        \label{fig:p2_xi_cond_C1_example3}
    \end{subfigure}

    \par\medskip

    \begin{subfigure}[t]{0.48\linewidth}
        \centering
        \includegraphics[width=\linewidth]{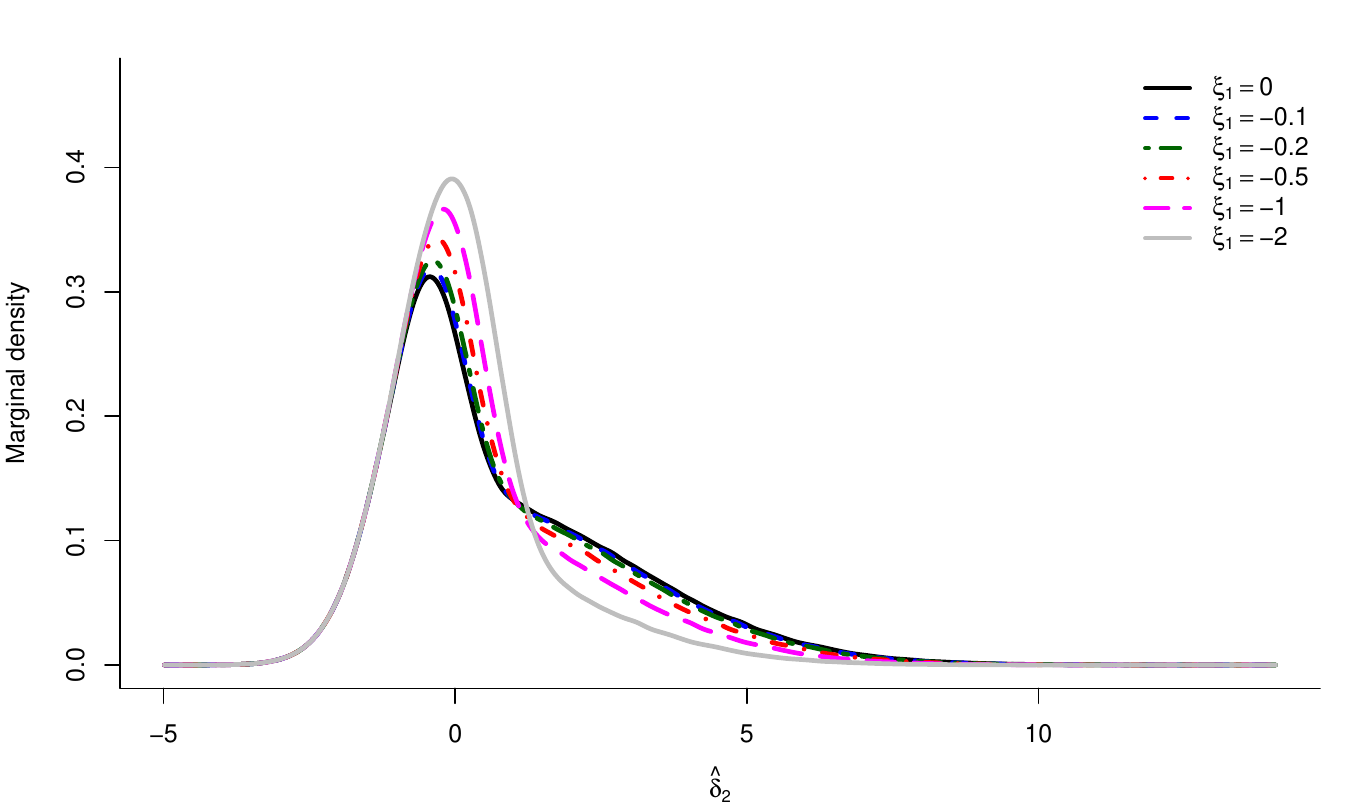}
        \caption{Marginal}
        \label{fig:p2_xi_marginal_example3}
    \end{subfigure}
    \caption{Example 3.
        Density estimates of $\widehat{\delta}_2$ for $\xi_1\in\{0,-0.1,-0.2,-0.5,-1, -2\}$:
        (a) conditional on $C=\emptyset$,
        (b) conditional on $C=\{1\}$, and
        (c) marginal.
    }
    \label{fig:p2_xi_cond_example3}
\end{figure}

\vspace{2mm}
\noindent {\bf Example 3}
For another example,
we consider
\begin{align*}
    W = \left( \begin{matrix}
        1 & -0.9 \\
        -0.9 & 1
    \end{matrix}\right),
    \qquad
    \Sigma = \left( \begin{matrix}
        5 & -2 \\
        -2 & 1
    \end{matrix} \right),
\end{align*}
and follow the same procedure as in the previous example for each $\xi_1 \in \{0, -0.1, -0.2, -0.5, -1, -2\}$.
Figure~\ref{fig:p2_xi_cond_example3} presents the estimated conditional and marginal densities of $\widehat{\delta}_2$ and Table~\ref{tab:p2_xi_skewness_example3} reports the skewness computed using the same formula as in the previous example.
In contrast to the previous examples, the skewness of the marginal distribution of $\widehat{\delta}_2$ increases as $\xi_1$ decreases from zero to $-2$.
However, the skewness eventually decreases as $\xi_1 $ decreases further, falling to $0.6099$ at $\xi_1=-4$ and $0.2495$ at $\xi_1=-5$.
This decrease is consistent with the distribution of $\widehat{\delta}_2$.
Specifically, as $\xi_1$ decreases, the condition $\widetilde{\delta}_1 > \xi_1$ in \eqref{eq:sim_p2_xi} becomes dominant, implying $\widehat{\delta}_2$ is closer to $\widehat{\delta}_2$ whose asymptotic distribution is $N(0,1)$ in this example.

\begin{table}[htbp]
    \centering
    \caption{Example 3. Skewness of the conditional and marginal distributions of $\widehat{\delta}_2$.}
    \label{tab:p2_xi_skewness_example3}
    \begin{tabular}{rrrr}
        \toprule
        $\xi_1$ & $C=\emptyset$ & $C=\{1\}$ & Marginal \\
        \midrule
         0.0 & 0.9049 & -0.4561 &  1.0426 \\
        -0.1 & 0.9178 & -0.4557 &  1.0716 \\
        -0.2 & 0.9304 & -0.4550 &  1.1000 \\
        -0.5 & 0.9665 & -0.4506 &  1.1808 \\
        -1.0 & 1.0189 & -0.4358 &  1.2917 \\
        -2.0 & 1.0981 & -0.3777 &  1.3437 \\
        \bottomrule
    \end{tabular}
\end{table}


\subsection{Case $p=3$}

\noindent {\bf Notation}
Before presenting the example,
we introduce notation used frequently in this section.
Suppose we are given a symmetric positive definite matrix
\begin{align*}
    W = \left( \begin{matrix}
        W_{11} & W_{12} & W_{13} \\ 
        W_{21} & W_{22} & W_{23} \\ 
        W_{31} & W_{32} & W_{33} \\ 
    \end{matrix}\right) \in \R^{3 \times 3}.
\end{align*}
Then we denote the following:
\begin{align*}
   & W_{(1,2), 3} = \left( \begin{matrix}
        W_{13} \\
        W_{23}
    \end{matrix} \right), 
    \quad
    W_{1, (2,3)} = \left( \begin{matrix}
        W_{12} & W_{13} 
    \end{matrix} \right),
    \quad
    W_{(1,2), (1,2)} = \left( \begin{matrix}
        W_{11} & W_{12} \\
        W_{21} & W_{22}
    \end{matrix}\right) \succ 0
\end{align*}
Using this notation, we further denote
\begin{align*}
     W_{1|2} &= W_{11}-W_{12}W_{22}^{-1}W_{21} > 0 , \\
    W_{12|3} &= W_{12} - W_{13} W_{33}^{-1} W_{32} \in \R, \\
     W_{1|(2,3)} &= W_{11} - W_{1,(2,3)} W_{(2,3),(2,3)}^{-1} W_{1,(2,3)}^\top >0,
    \\
    W_{(1,2)|3} &= W_{(1,2),(1,2)} - W_{(1,2),3}W_{33}^{-1}W_{(1,2),3}^\top \succ 0,
\end{align*}
with analogous notation for other index sets.
The same notation is used for $V^C$ and $\Sigma$.

\vspace{2mm}
\noindent {\bf Example 1. $B=\{1,2,3\}$ and $U=\emptyset$}
\vspace{2mm}
Consider an example with 
$p=3$, $B=\{1,2,3\}, \mathcal{P}(B)=\{\emptyset, \{1\}, \{2\}, \{3\} ,\{1,2\}, \{1,3\}, \{2,3\}, \{1,2,3\}\}$ and $U=\emptyset$.
Suppose the population parameter lies on the boundary (thus $\xi=0$).
From Theorem~\ref{theorem:blocks},
\begin{align*}
    \widehat{\delta}
    &=
    \left( \begin{matrix}
        \widetilde{\delta}_1 \\
        \widetilde{\delta}_2 \\
        \widetilde{\delta}_3
    \end{matrix} \right) I[\widetilde{\delta}_1 > 0, \widetilde{\delta}_2>0, \widetilde{\delta}_3>0]
    +
    \left( \begin{matrix}
        \delta^{12}_1 \\
        \delta^{12}_2 \\
        0
    \end{matrix} \right) I[\delta^{12}_1>0, \delta^{12}_2>0, \widetilde{\delta}_3<0] \\
    &+
    \left( \begin{matrix}
        \delta^{13}_1 \\
        0 \\
        \delta^{13}_3
    \end{matrix} \right) I[\delta^{13}_1>0, \delta^{13}_3>0, \widetilde{\delta}_2<0]
    +
    \left( \begin{matrix}
        0 \\
        \delta^{23}_2 \\
        \delta^{23}_3
    \end{matrix} \right) I[\delta^{23}_2>0, \delta^{23}_3>0, \widetilde{\delta}_1<0] \\
    &+
    \left( \begin{matrix}
        \delta^1 \\
        0 \\
        0
    \end{matrix} \right) I\bigg[ \delta^1>0, W_{(2,3)|1} \left( \begin{matrix}
        \widetilde{\delta}_2 \\
        \widetilde{\delta}_3
    \end{matrix}\right) < 0 \bigg]
    +
    \left( \begin{matrix}
        0 \\
        \delta^2 \\
        0
    \end{matrix} \right) I\bigg[\delta^2>0, W_{(1,3)|2} \left( \begin{matrix}
        \widetilde{\delta}_1 \\
        \widetilde{\delta}_3
    \end{matrix}\right) < 0 \bigg] \\
    &+
    \left( \begin{matrix}
        0 \\
        0 \\
        \delta^3
    \end{matrix} \right) I\bigg[ \delta^3>0, W_{(1,2)|3} \left( \begin{matrix}
        \widetilde{\delta}_1 \\
        \widetilde{\delta}_2
    \end{matrix}\right) < 0 \bigg]
    + \left( \begin{matrix}
        0 \\
        0 \\
        0
    \end{matrix} \right) I[\delta^1<0, \delta^2<0, \delta^3<0],
\end{align*}
where we define
\begin{align*}
    &\delta^{12} 
    = \left( \begin{matrix}
        \delta^{12}_1 \\
        \delta^{12}_2
    \end{matrix} \right)
    = \left( \begin{matrix}
        \widetilde{\delta}_1 \\
        \widetilde{\delta}_2
    \end{matrix} \right) + W_{(1,2),(1,2)}^{-1}W_{(1,2),3}\widetilde{\delta}_3,
    \qquad
    \delta^{13} 
    = \left( \begin{matrix}
        \delta^{13}_1 \\
        \delta^{13}_3
    \end{matrix} \right)
    = \left( \begin{matrix}
        \widetilde{\delta}_1 \\
        \widetilde{\delta}_3
    \end{matrix} \right) + W_{(1,3),(1,3)}^{-1}W_{(1,3),2}\widetilde{\delta}_2 \\
    &\delta^{23} 
    = \left( \begin{matrix}
        \delta^{23}_2 \\
        \delta^{23}_3
    \end{matrix} \right)
    = \left( \begin{matrix}
        \widetilde{\delta}_2 \\
        \widetilde{\delta}_3
    \end{matrix} \right) + W_{(2,3),(2,3)}^{-1}W_{(2,3),1}\widetilde{\delta}_1,
    \qquad
    \delta^1 = \widetilde{\delta}_1 + W_{11}^{-1}(W_{12}\widetilde{\delta}_2+W_{13}\widetilde{\delta}_3), \\
    &\delta^2 = \widetilde{\delta}_2 + W_{22}^{-1}(W_{21}\widetilde{\delta}_1+W_{23}\widetilde{\delta}_3),
    \qquad
    \delta^3 = \widetilde{\delta}_3 + W_{33}^{-1}(W_{31}\widetilde{\delta}_1+W_{32}\widetilde{\delta}_2).
\end{align*}

To study the asymptotic distribution of $\widehat{\delta}$, 
we consider the conditional distribution associated with each term in the above.
To this end, we use a matrix $T_C$ defined in Theorem~\ref{theorem:CLT}.
Since $U=\emptyset$, the block corresponding to $U$ disappears, and $T_C$ reduces to a $3 \times 3$ matrix.
\begin{itemize}
    \item[(i)] If $C=\{1,2\}$, then $A_C=\{3\}$ and $\sqrt{n}T_{\{1,2\}} \widetilde{\delta} \Dcon N(0, V^{\{1,2\}})$ where
    \begin{align*}
        V^{\{1,2\}} 
        &= T_{\{1,2\}}\Sigma T_{\{1,2\}}^\top
        =
    \left(
        \begin{matrix}
        v_{11}^{\{1,2\}} & v_{12}^{\{1,2\}} & v_{13}^{\{1,2\}} \\
         v_{21}^{\{1,2\}} & v_{22}^{\{1,2\}} & v_{23}^{\{1,2\}} \\
         v_{31}^{\{1,2\}} & v_{32}^{\{1,2\}} & v_{33}^{\{1,2\}} 
         \end{matrix}
        \right)
        \quad \textrm{where }
        T_{\{1,2\}} = \left( \begin{matrix}
        I_{2} & W_{(1,2),(1,2)}^{-1}W_{(1,2),3} \\
        0 & W_{3|(1,2)}
    \end{matrix} \right)
    \end{align*}
    \item[(ii)] If $C=\{1,3\}$, then $A_C=\{2\}$ and $\sqrt{n}T_{\{1,3\}} \widetilde{\delta} \Dcon N(0, V^{\{1,3\}})$ where
    \begin{align*}
        V^{\{1,3\}} = T_{\{1,3\}}\Sigma T_{\{1,3\}}^\top
        =
    \left(
        \begin{matrix}
        v_{11}^{\{1,3\}} & v_{12}^{\{1,3\}} & v_{13}^{\{1,3\}} \\
         v_{21}^{\{1,3\}} & v_{22}^{\{1,3\}} & v_{23}^{\{1,3\}} \\
         v_{31}^{\{1,3\}} & v_{32}^{\{1,3\}} & v_{33}^{\{1,3\}} 
         \end{matrix}
        \right)
        \quad \textrm{where }
        T_{\{1,3\}} = \left( \begin{matrix}
        1 & W_{1|3}^{-1} W_{12|3} & 0 \\
        0 & W_{2|(1,3)} & 0 \\
        0 & W_{3|1}^{-1} W_{32|1} & 1
    \end{matrix}\right),
    \end{align*}
    where we use $W_{(1,3),(1,3)}^{-1}W_{(1,3),2} = \left( \begin{matrix}
        W_{1|3}^{-1} W_{12|3} \\
        W_{3|1}^{-1} W_{32|1}
    \end{matrix}\right)$.
    \item[(iii)] If $C=\{2,3\}$, then $A_C=\{1\}$ and $\sqrt{n}T_{\{2,3\}} \widetilde{\delta} \Dcon N(0, V^{\{2,3\}})$ where
    \begin{align*}
        V^{\{2,3\}} = T_{\{2,3\}}\Sigma T_{\{2,3\}}^\top
        =
    \left(
        \begin{matrix}
        v_{11}^{\{2,3\}} & v_{12}^{\{2,3\}} & v_{13}^{\{2,3\}} \\
         v_{21}^{\{2,3\}} & v_{22}^{\{2,3\}} & v_{23}^{\{2,3\}} \\
         v_{31}^{\{2,3\}} & v_{32}^{\{2,3\}} & v_{33}^{\{2,3\}} 
         \end{matrix}
        \right)
        \quad \textrm{where }
        T_{\{2,3\}} = 
    \left( \begin{matrix}
        W_{1|(2,3)} & 0 \\
        W_{(2,3),(2,3)}^{-1}W_{(2,3),1} & I_2
    \end{matrix} \right)
    \end{align*}
    \item[(iv)] If $C=\{1\}$, then $A_C=\{2,3\}$ and $\sqrt{n}T_{\{1\}} \widetilde{\delta} \Dcon N(0, V^{\{1\}})$ where
    \begin{align*}
        V^{\{1\}} = T_{\{1\}}\Sigma T_{\{1\}}^\top
         =
    \left(
        \begin{matrix}
        v_{11}^{\{1\}} & v_{12}^{\{1\}} & v_{13}^{\{1\}} \\
         v_{21}^{\{1\}} & v_{22}^{\{1\}} & v_{23}^{\{1\}} \\
         v_{31}^{\{1\}} & v_{32}^{\{1\}} & v_{33}^{\{1\}} 
         \end{matrix}
        \right)
        \quad \textrm{where }
        T_{\{1\}} = \left(
        \begin{matrix}
            1 & W_{11}^{-1}W_{1,(2,3)} \\
            0 & W_{(2,3)|1}
        \end{matrix}
        \right)
    \end{align*}
    \item[(v)] If $C=\{2\}$, then $A_C=\{1,3\}$ and $\sqrt{n}T_{\{2\}} \widetilde{\delta} \Dcon N(0, V^{\{2\}})$ where
    \begin{align*}
        V^{\{2\}} = T_{\{2\}}\Sigma T_{\{2\}}^\top
        =
    \left(
        \begin{matrix}
        v_{11}^{\{2\}} & v_{12}^{\{2\}} & v_{13}^{\{2\}} \\
         v_{21}^{\{2\}} & v_{22}^{\{2\}} & v_{23}^{\{2\}} \\
         v_{31}^{\{2\}} & v_{32}^{\{2\}} & v_{33}^{\{2\}} 
         \end{matrix}
        \right)
        \quad \textrm{where }
        T_{\{2\}} = \left(
        \begin{matrix}
            W_{1|2} & 0 & W_{13|2} \\
            W_{22}^{-1}W_{21} & 1 & W_{22}^{-1}W_{23} \\
            W_{31|2} & 0 & W_{3|2}
        \end{matrix}
        \right)
    \end{align*}
    \item[(vi)] If $C=\{3\}$, then $A_C=\{1,2\}$ and $\sqrt{n}T_{\{3\}} \widetilde{\delta} \Dcon N(0, V^{\{3\}})$ where
    \begin{align*}
        V^{\{3\}} = T_{\{3\}}\Sigma T_{\{3\}}^\top 
        =
    \left(
        \begin{matrix}
        v_{11}^{\{3\}} & v_{12}^{\{3\}} & v_{13}^{\{3\}} \\
         v_{21}^{\{3\}} & v_{22}^{\{3\}} & v_{23}^{\{3\}} \\
         v_{31}^{\{3\}} & v_{32}^{\{3\}} & v_{33}^{\{3\}} 
         \end{matrix}
        \right)
        \quad \textrm{where }
        T_{\{3\}} =  \left(
        \begin{matrix}
            W_{(1,2)|3} & 0 \\
            W_{33}^{-1}W_{3,(1,2)} & 1
        \end{matrix}
        \right).
    \end{align*}
\end{itemize}
Using Lemma~\ref{lemma:SN_TSN} gives the conditional distributions for $C \in \{ \{1\}, \{2\}, \{3\}, \{1,2\}, \{1,3\}, \{2,3\} \}$.
As a result, the asymptotic distribution is
\begin{align*}
    \sqrt{n} \widehat{\delta}
    &\Dcon
    tN_3(0, \Sigma) \operatorname{pr}[\widetilde{\delta}_1 > 0, \widetilde{\delta}_2>0, \widetilde{\delta}_3>0] \\
    &+
    \left( \begin{matrix}
        tSN_{1,2}(0, 0, v_{3,(1,2)}^{\{1,2\}} (v_{(1,2),(1,2)}^{\{1,2\}})^{-1}
        ,
        v_{(1,2),(1,2)}^{\{1,2\}}, v_{3|(1,2)}^{\{1,2\}}) \\
        m_{\{3\}}(0)
    \end{matrix} \right) \operatorname{pr}[\delta^{12}_1>0, \delta^{12}_2>0, \widetilde{\delta}_3<0] \\
    &+
    \left( \begin{matrix}
        [tSN_{1,2}(0, 0, v_{2,(1,3)}^{\{1,3\}} (v_{(1,3),(1,3)}^{\{1,3\}})^{-1},
        v_{(1,3),(1,3)}^{\{1,3\}}, v_{2|(1,3)}^{\{1,3\}})]_1 \\
        m_{\{2\}}(0) \\
        [tSN_{1,2}(0, 0, v_{2,(1,3)}^{\{1,3\}} (v_{(1,3),(1,3)}^{\{1,3\}})^{-1},
        v_{(1,3),(1,3)}^{\{1,3\}}, v_{2|(1,3)}^{\{1,3\}})]_2
    \end{matrix} \right) \operatorname{pr}[\delta^{13}_1>0, \delta^{13}_3>0, \widetilde{\delta}_2<0] \\
    &+
    \left( \begin{matrix}
    m_{\{1\}}(0) \\
        tSN_{1,2}(0, 0, v_{1,(2,3)}^{\{2,3\}} (v_{(2,3),(2,3)}^{\{2,3\}})^{-1},
        v_{(2,3),(2,3)}^{\{2,3\}}, v_{1|(2,3)}^{\{2,3\}})
    \end{matrix} \right) \operatorname{pr}[\delta^{23}_2>0, \delta^{23}_3>0, \widetilde{\delta}_1<0] \\
    &+
     \left( \begin{matrix}
        tSN_{2,1}(0, 0, v_{(2,3),1}^{\{1\}}/ v_{11}^{\{1\}},
        v_{11}^{\{1\}}, v_{(2,3)|1}^{\{1\}}) \\
        m_{\{2\}}(0)\\
        m_{\{3\}}(0)
    \end{matrix} \right)
    \operatorname{pr}\bigg[ \delta^1>0, W_{(2,3)|1} \left( \begin{matrix}
        \widetilde{\delta}_2 \\
        \widetilde{\delta}_3
    \end{matrix}\right) < 0 \bigg] \\
    &+
    \left( \begin{matrix}
        m_{\{1\}}(0) \\
        tSN_{2,1}(0, 0, v_{(1,3),2}^{\{2\}}/ v_{22}^{\{2\}},
        v_{22}^{\{2\}}, v_{(1,3)|2}^{\{2\}}) \\
        m_{\{3\}}(0)
    \end{matrix} \right)
    \operatorname{pr} \bigg[ \delta^2>0, W_{(1,3)|2} \left( \begin{matrix}
        \widetilde{\delta}_1 \\
        \widetilde{\delta}_3
    \end{matrix}\right) < 0 \bigg] \\
    &+
    \left( \begin{matrix}
        m_{\{1\}}(0) \\
        m_{\{2\}}(0) \\
        tSN_{2,1}(0, 0, v_{(1,2),3}^{\{3\}}/ v_{33}^{\{3\}},
        v_{33}^{\{3\}}, v_{(1,2)|3}^{\{3\}}) \\
    \end{matrix} \right)
    \operatorname{pr} \bigg[ \delta^3>0, W_{(1,2)|3} \left( \begin{matrix}
        \widetilde{\delta}_1 \\
        \widetilde{\delta}_2
    \end{matrix}\right) < 0 \bigg] \\
    &+
    \left(\begin{matrix}
        m_{\{1\}}(0) \\
        m_{\{2\}}(0)\\
        m_{\{3\}}(0)
    \end{matrix} \right)
    \operatorname{pr}[\delta^1<0, \delta^2<0, \delta^3<0],
\end{align*}
which also comes directly from Theorem~\ref{theorem:CLT}.
Here, we use $[\cdot]_j$ to denote its $j$th component for a multivariate distribution.

\vspace{2mm}
\noindent {\bf Special case $\Sigma=W^{-1}$}
Applying Corollary \ref{theorem:CLT3} gives
\begin{align*}
    \sqrt{n} \widehat{\delta}
    &\Dcon
    tN_3(0, \Sigma) \operatorname{pr}[\widetilde{\delta}_1 > 0, \widetilde{\delta}_2>0, \widetilde{\delta}_3>0] 
    +
    \left( \begin{matrix}
        tN_2(0, \Sigma_{(1,2)|3}) \\
        m_{\{3\}}(0)
    \end{matrix} \right) \operatorname{pr}[\delta^{12}_1>0, \delta^{12}_2>0, \widetilde{\delta}_3<0] \\
    &+
    \left( \begin{matrix}
        [tN_2(0, \Sigma_{(1,3)|2})]_1 \\
        m_{\{2\}}(0) \\
        [tN_2(0, \Sigma_{(1,3)|2})]_2
    \end{matrix} \right) \operatorname{pr}[\delta^{13}_1>0, \delta^{13}_3>0, \widetilde{\delta}_2<0] 
    +
    \left( \begin{matrix}
    m_{\{1\}}(0) \\
        tN_2(0, \Sigma_{(2,3)|1})
    \end{matrix} \right) \operatorname{pr}[\delta^{23}_2>0, \delta^{23}_3>0, \widetilde{\delta}_1<0] \\
    &+
     \left( \begin{matrix}
        tN_1(0, \Sigma_{1|(2,3)}) \\
        m_{\{2\}}(0)\\
        m_{\{3\}}(0)
    \end{matrix} \right)
    \operatorname{pr}\bigg[ \delta^1>0, \Sigma_{(2,3),(2,3)}^{-1} \left( \begin{matrix}
        \widetilde{\delta}_2 \\
        \widetilde{\delta}_3
    \end{matrix}\right) < 0 \bigg]  \\
    &+
    \left( \begin{matrix}
        m_{\{1\}}(0) \\
        tN_1(0, \Sigma_{2|(1,3)}) \\
        m_{\{3\}}(0)
    \end{matrix} \right)
    \operatorname{pr}\bigg[ \delta^2>0, \Sigma_{(1,3),(1,3)}^{-1} \left( \begin{matrix}
        \widetilde{\delta}_1 \\
        \widetilde{\delta}_3
    \end{matrix}\right) < 0 \bigg] \\
    &+
    \left( \begin{matrix}
        m_{\{1\}}(0) \\
        m_{\{2\}}(0) \\
        tN_1(0, \Sigma_{3|(1,2)}) \\
    \end{matrix} \right)
    \operatorname{pr}\bigg[ \delta^3>0, \Sigma_{(1,2),(1,2)}^{-1} \left( \begin{matrix}
        \widetilde{\delta}_1 \\
        \widetilde{\delta}_2
    \end{matrix}\right) < 0 \bigg]
    +
    \left(\begin{matrix}
        m_{\{1\}}(0) \\
        m_{\{2\}}(0)\\
        m_{\{3\}}(0)
    \end{matrix} \right)
    \operatorname{pr}[\delta^1<0, \delta^2<0, \delta^3<0],
\end{align*}

\vspace{2mm}
\noindent {\bf The Quasi LR test}
We consider the one-sided hypothesis
\begin{align*}
    H_0: \theta_B = 0 \quad \textrm{vs.} \quad
    H_1: \theta_B \geq 0 \textrm{ and }\theta_{B}\neq 0,
\end{align*}
in the case that $B=\{1,2,3\}$ and $U=\emptyset$.
From Theorem~\ref{theorem:LLRexpansion},
\begin{align*}
    LR 
    &= n
    \bigg( 
    (\widetilde{\delta})^\top W \widetilde{\delta} \cdot I[\widetilde{\delta}_1 > 0, \widetilde{\delta}_2>0, \widetilde{\delta}_3>0]
    +
    (\delta^{12})^\top W_{(1,2),(1,2)} \delta^{12} \cdot
    I[\delta^{12}_1>0, \delta^{12}_2>0, \widetilde{\delta}_3<0] \\
    &+
    (\delta^{13})^\top W_{(1,3),(1,3)} \delta^{13} \cdot
    I[\delta^{13}_1>0, \delta^{13}_3>0, \widetilde{\delta}_2<0]
    +
    (\delta^{23})^\top W_{(2,3),(2,3)} \delta^{23} \cdot
    I[\delta^{23}_2>0, \delta^{23}_3>0, \widetilde{\delta}_1<0] \\
    &+
    W_{11} (\delta^{1})^2 \cdot I\bigg[ \delta^1>0, W_{(2,3)|1} \left( \begin{matrix}
        \widetilde{\delta}_2 \\
        \widetilde{\delta}_3
    \end{matrix}\right) < 0 \bigg]
    +
    W_{22} (\delta^{2})^2 \cdot I\bigg[ \delta^2>0, W_{(1,3)|2} \left( \begin{matrix}
        \widetilde{\delta}_1 \\
        \widetilde{\delta}_3
    \end{matrix}\right) < 0 \bigg] \\
    &+
    W_{33} (\delta^{3})^2 \cdot I\bigg[ \delta^3>0, W_{(1,2)|3} \left( \begin{matrix}
        \widetilde{\delta}_1 \\
        \widetilde{\delta}_2
    \end{matrix}\right) < 0 \bigg] \bigg)
    + 0 \cdot I[\delta^1<0, \delta^2<0, \delta^3<0] + o_p(1).
\end{align*}
In the special case where $\Sigma = W^{-1}$, 
applying 
Corollary \ref{cor:pivotal}
gives the well known result that 
$LR$ is a mixture of 
$\chi^{2}$-distributions
\begin{align*}
    LR
    &\Dcon
    \chi^2_3 \cdot \operatorname{pr}[\widetilde{\delta}_1 > 0, \widetilde{\delta}_2>0, \widetilde{\delta}_3>0] \\
    &+
    \chi^2_2 \cdot \left( \operatorname{pr}[\delta^{12}_1>0, \delta^{12}_2>0, \widetilde{\delta}_3<0] + \operatorname{pr}[\delta^{13}_1>0, \delta^{13}_3>0, \widetilde{\delta}_2<0] + \operatorname{pr}[\delta^{23}_2>0, \delta^{23}_3>0, \widetilde{\delta}_1<0] \right) \\
    &+
    \chi^2_1 \cdot \bigg( \operatorname{pr}\left[\delta^1>0, \Sigma_{(2,3),(2,3)}^{-1} \left( \begin{matrix}
        \widetilde{\delta}_2 \\
        \widetilde{\delta}_3
    \end{matrix}\right) < 0\right] + \operatorname{pr}\left[ \delta^2>0, \Sigma_{(1,3),(1,3)}^{-1} \left( \begin{matrix}
        \widetilde{\delta}_1 \\
        \widetilde{\delta}_3
    \end{matrix}\right) < 0\right] \\
    &+ 
    \operatorname{pr}\left[ \delta^3>0, \Sigma_{(1,2),(1,2)}^{-1} \left( \begin{matrix}
        \widetilde{\delta}_1 \\
        \widetilde{\delta}_2
    \end{matrix}\right) < 0\right] \bigg)
    +
    \chi^2_0 \cdot \operatorname{pr}[\delta^1<0, \delta^2<0, \delta^3<0].
\end{align*}

\vspace{2mm}
\noindent {\bf Example 2. $B=\{1,2\}$ and $U=\{3\}$}
\vspace{2mm}
Consider the case where $p=3$, $B=\{1,2\}, \mathcal{P}(B)=\{\emptyset, \{1\}, \{2\}, \{1,2\}\}$ and $U=\{3\}$.
Suppose the true parameter lies on the boundary (thus $\xi=0$).
This setting is similar to the previous $p=2$ example in Section~\ref{sec:example_p2} except that we here include an unrestricted parameter indexed by $U$.
From Theorem~\ref{theorem:blocks},
\begin{align*}
    \widehat{\delta}
    &= \left( \begin{matrix}
        \widetilde{\delta}_1 \\
        \widetilde{\delta}_2 \\
        \widetilde{\delta}_3
    \end{matrix} \right) I[\widetilde{\delta}_1>0, \widetilde{\delta}_2>0]
    + \left( \begin{matrix}
        \delta^{13}_1 \\
        0 \\
        \delta^{13}_3
    \end{matrix} \right) I[\delta^{13}_1>0, \widetilde{\delta}_2<0] \\
    &+
    \left( \begin{matrix}
        0 \\
        \delta^{23}_2 \\
        \delta^{23}_3
    \end{matrix} \right) I[\delta^{23}_2>0, \widetilde{\delta}_1<0]
    +
    \left( \begin{matrix}
        0 \\
        0 \\
        \delta^{3}
    \end{matrix} \right) I[ \delta^{13}_1 < 0, \delta^{23}_2<0],
\end{align*}
where we define
\begin{align*}
    &\delta^{13} 
    = \left( \begin{matrix}
        \delta^{13}_1 \\
        \delta^{13}_3
    \end{matrix} \right)
    = \left( \begin{matrix}
        \widetilde{\delta}_1 \\
        \widetilde{\delta}_3
    \end{matrix} \right) + 
    \left( \begin{matrix}
        W_{1|3}^{-1} W_{12|3} \\
        W_{3|1}^{-1} W_{32|1}
    \end{matrix}\right) \widetilde{\delta}_2,
    \quad
    \delta^{23} 
    = \left( \begin{matrix}
        \delta^{23}_2 \\
        \delta^{23}_3
    \end{matrix} \right)
    = \left( \begin{matrix}
        \widetilde{\delta}_2 \\
        \widetilde{\delta}_3
    \end{matrix} \right) +
    \left( \begin{matrix}
        W_{2|3}^{-1} W_{21|3} \\
        W_{3|2}^{-1} W_{31|2}
    \end{matrix}\right) \widetilde{\delta}_1, \\
    &\delta^{3} = \widetilde{\delta}_3+W_{33}^{-1}(W_{31}\widetilde{\delta}_1+W_{32} \widetilde{\delta}_2).
\end{align*}
As in the previous example, we use a $3 \times 3$ matrix $T_C$ to derive the asymptotic distribution of $\widehat{\delta}$.
In contrast to the previous examples, we need to account for the block corresponding to $U=\{3\}$.

\begin{itemize}
    \item[(i)] If $C=\{1\}$, then $A_C=\{2\}$, $F_C = \{1,3\}$, and $\sqrt{n}T_{\{1\}} \widetilde{\delta} \Dcon N(0, V^{\{1\}})$ where
    \begin{align*}
        V^{\{1\}} = T_{\{1\}} \Sigma T_{\{1\}}^\top
        =
    \left(
        \begin{matrix}
        v_{11}^{\{1\}} & v_{12}^{\{1\}} & v_{13}^{\{1\}} \\
         v_{21}^{\{1\}} & v_{22}^{\{1\}} & v_{23}^{\{1\}} \\
         v_{31}^{\{1\}} & v_{32}^{\{1\}} & v_{33}^{\{1\}} 
         \end{matrix}
        \right)
        \quad \textrm{where }
        T_{\{1\}} = \left( \begin{matrix}
        1 & W_{1|3}^{-1} W_{12|3} & 0 \\
        0 & W_{2|(1,3)} & 0 \\
        0 & W_{3|1}^{-1} W_{32|1} & 1
        \end{matrix} \right).
    \end{align*}
    \item[(ii)] If $C=\{2\}$, then $A_C = \{1\}$, $F_C = \{2,3\}$, and $\sqrt{n}T_{\{2\}} \widetilde{\delta} \Dcon N(0, V^{\{2\}})$ where
    \begin{align*}
        V^{\{2\}} = T_{\{2\}} \Sigma T_{\{2\}}^\top
        =
    \left(
        \begin{matrix}
        v_{11}^{\{2\}} & v_{12}^{\{2\}} & v_{13}^{\{2\}} \\
         v_{21}^{\{2\}} & v_{22}^{\{2\}} & v_{23}^{\{2\}} \\
         v_{31}^{\{2\}} & v_{32}^{\{2\}} & v_{33}^{\{2\}} 
         \end{matrix}
        \right)
        \quad \textrm{where }
        T_{\{2\}} = \left( \begin{matrix}
        W_{1|(2,3)} & 0 & 0 \\
        W_{2|3}^{-1}W_{21|3} & 1 & 0 \\
        W_{3|2}^{-1}W_{31|2} & 0 & 1
        \end{matrix} \right).
    \end{align*}
    \item[(iii)] If $C= \emptyset$, then $A_C=\{1,2\}$, $F_C = \{3\}$, and $\sqrt{n}T_{\emptyset}\widetilde{\delta} \Dcon N(0, V^{\emptyset})$ where
    \begin{align*}
        V^{\emptyset} = T_{\emptyset} \Sigma T_{\emptyset}^\top
        =
    \left(
        \begin{matrix}
        v_{11}^{\emptyset} & v_{12}^{\emptyset} & v_{13}^{\emptyset} \\
         v_{21}^{\emptyset} & v_{22}^{\emptyset} & v_{23}^{\emptyset} \\
         v_{31}^{\emptyset} & v_{32}^{\emptyset} & v_{33}^{\emptyset} 
         \end{matrix}
        \right)
        \quad \textrm{where }
        T_{\emptyset} = \left( \begin{matrix}
        W_{(1,2)|3} & 0 \\
        W_{33}^{-1}W_{3,(1,2)} & 1
        \end{matrix} \right).
    \end{align*}
\end{itemize}
Lemma~\ref{lemma:SN_TSN} gives the conditional distributions for each $C$.
In consequence, the asymptotic distribution follows from Theorem~\ref{theorem:CLT}
\begin{align*}
    \sqrt{n} \widehat{\delta}
    &\Dcon 
    \left( \begin{matrix}
        tN_2(0, \Sigma_{(1,2),(1,2)}) \\
        \widetilde{SN}_{2,1}(0, \Sigma_{(1,2), 3}/\Sigma_{33}, \Sigma_{33}, \Sigma_{(1,2)|3})
    \end{matrix}\right) \operatorname{pr}[\widetilde{\delta}_1>0, \widetilde{\delta}_2>0] \\
    &+
     \left( \begin{matrix}
        tSN_{1,1}(0, 0, v^{\{1\}}_{21}/ v^{\{1\}}_{11}, v^{\{1\}}_{11}, v^{\{1\}}_{2|1}) \\
        m_{\{2\}}(0) \\
        SN_{2,1}(0, 0, v^{\{1\}}_{23}/v^{\{1\}}_{33}, v^{\{1\}}_{13}/v^{\{1\}}_{33}, v^{\{1\}}_{33}, v^{\{1\}}_{(2,1)|3})
    \end{matrix}\right) \operatorname{pr}[\delta^{13}_1>0, \widetilde{\delta}_2<0] \\
    &+
    \left( \begin{matrix}
    m_{\{1\}}(0) \\
        tSN_{1,1}(0, 0, v^{\{2\}}_{12}/ v^{\{2\}}_{22}, v^{\{2\}}_{22}, v^{\{2\}}_{1|2}) \\
        SN_{2,1}(0, 0, v^{\{2\}}_{13}/v^{\{2\}}_{33}, v^{\{2\}}_{23}/v^{\{2\}}_{33}, v^{\{2\}}_{33}, v^{\{2\}}_{(1,2)|3})
    \end{matrix}\right) \operatorname{pr}[\delta^{23}_2>0, \widetilde{\delta}_1<0] \\
    &+
    \left( \begin{matrix}
    m_{\{1\}}(0) \\
        m_{\{2\}}(0) \\
        \widetilde{SN}_{2,1}(0, -v^{\emptyset}_{(1,2),3}/v^{\emptyset}_{33}, v^{\emptyset}_{33}, v^{\emptyset}_{(1,2)|3})
    \end{matrix}\right) \operatorname{pr}[\delta^{13}_1 < 0, \delta^{23}_2<0].
\end{align*}

\vspace{2mm}
\noindent {\bf Special case $\Sigma=W^{-1}$}
Applying Corollary \ref{theorem:CLT3} gives
\begin{align*}
    \sqrt{n} \widehat{\delta}
    &\Dcon 
    \left( \begin{matrix}
        tN_2(0, \Sigma_{(1,2),(1,2)}) \\
        \widetilde{SN}_{2,1}(0, \Sigma_{(1,2), 3}/\Sigma_{33}, \Sigma_{33}, \Sigma_{(1,2)|3})
    \end{matrix}\right) \operatorname{pr}[\widetilde{\delta}_1>0, \widetilde{\delta}_2>0] \\
    &+
     \left( \begin{matrix}
        tN_{1}(0, \Sigma_{1|2}) \\
        m_{\{2\}}(0) \\
        \widetilde{SN}_{1,1}(0, \Sigma_{13|2}\Sigma_{3|2}^{-1}, \Sigma_{3|2}, \Sigma_{1|(2,3)})
    \end{matrix}\right) \operatorname{pr}[\delta^{13}_1>0, \widetilde{\delta}_2<0] \\
    &+
    \left( \begin{matrix}
    m_{\{1\}}(0) \\
        tN_{1}(0, \Sigma_{2|1}) \\
        \widetilde{SN}_{1,1}(0, \Sigma_{23|1}\Sigma_{3|1}^{-1}, \Sigma_{3|1}, \Sigma_{2|(1,3)})
    \end{matrix}\right) \operatorname{pr}[\delta^{23}_2>0, \widetilde{\delta}_1<0]
    +
    \left( \begin{matrix}
    m_{\{1\}}(0) \\
        m_{\{2\}}(0) \\
       N(0, \Sigma_{3|(1,2)})
    \end{matrix}\right) \operatorname{pr}[\delta^{13}_1 < 0, \delta^{23}_2<0].
\end{align*}
When $C=\emptyset$, we can observe that one-sided skewed normal distribution defined in Definition~\ref{def:SN_Variants} reduces to a normal distribution. 
Hence, in this example,
$\widetilde{SN}_{0, 1}(0, 0, \Sigma_{3|(1,2)},0)$ reduces to $N(0, \Sigma_{3|(1,2)})$.

\vspace{2mm}
\noindent {\bf The Quasi LR test}
We consider the one-sided hypothesis 
\begin{align*}
    H_0: \theta_B = 0 \quad \textrm{vs.} \quad
    H_1: \theta_B \geq 0 \textrm{ and }\theta_{B}\neq 0,
\end{align*}
in the case that $B=\{1,2\}$ and $U=\{3\}$.
Theorem~\ref{theorem:LLRexpansion} gives
\begin{align*}
    LR
    &=
    n   \bigg(
    \left( \begin{matrix}
        \widetilde{\delta}_1 & \widetilde{\delta}_2
    \end{matrix} \right) 
    W_{(1,2)|3}
    \left( \begin{matrix}
        \widetilde{\delta}_1 \\ 
        \widetilde{\delta}_2
    \end{matrix} \right) \cdot
    I[\widetilde{\delta}_1>0, \widetilde{\delta}_2>0]
    + 
    W_{1|3} (\delta_{1}^{13})^2 \cdot I[\delta^{13}_1>0, \widetilde{\delta}_2<0]  \\
    &+
    W_{2|3} (\delta_{2}^{23})^2 \cdot I[\delta^{23}_2>0, \widetilde{\delta}_1<0] \bigg)
    +
    0 \cdot I[\delta^{13}_1<0, \delta^{23}_2<0] + o_p(1).
\end{align*}
In the special case where $\Sigma = W^{-1}$, Corollary~\ref{cor:pivotal} gives
\begin{align*}
    LR 
    &\Dcon \chi^2_2 \cdot \operatorname{pr}[\widetilde{\delta}_1>0, \widetilde{\delta}_2>0] 
    + \chi^2_1 \cdot \left( \operatorname{pr}[\delta^{13}_1>0, \widetilde{\delta}_2<0]+\operatorname{pr}[\delta^{23}_2>0, \widetilde{\delta}_1<0] \right) 
    + \chi^2_0 \cdot \operatorname{pr}[\delta^{13}_1<0, \delta^{23}_2<0].
\end{align*}
Thus, $LR$ is a mixture of $\chi^2$ distributions with $0, 1$, and $2$ degrees of freedom.
It has the same mixture components as in the previous $p=2$ example in Section~\ref{sec:example_p2}, but with different mixture weights due to the presence of $U$.

\end{document}